\documentclass{amsart}
\usepackage{amssymb,latexsym,amsxtra,amscd,ifthen}
\usepackage{amsmath}
\usepackage{amsfonts}
\usepackage{verbatim}
\usepackage[arrow,matrix,cmtip]{xy}

\numberwithin{equation}{section}
\theoremstyle{plain}
\newtheorem{theorem}[equation]{Theorem}

\newtheorem{lemma}[equation]{Lemma}

\newtheorem{corollary}[equation]{Corollary}
\newtheorem{proposition}[equation]{Proposition}

\newtheorem{question}[equation]{Question}

\newtheorem{example}[equation]{Example}

\theoremstyle{definition}
\newtheorem{definition}[equation]{Definition}
\newtheorem{notation}[equation]{Notation}

\newtheorem{assumptions}[equation]{Assumptions}

\newtheorem{remark}[equation]{Remark}
\newtheorem{remarks}[equation]{Remarks}

\newcommand{\fun}{\operatorname{Fun}}

\newcommand{\schemes}{\operatorname{Schemes}}
\newcommand{\rings}{\operatorname{Rings}}
\newcommand{\Morph}{\operatorname{Morph}}
\newcommand{\sets}{\operatorname{Sets}}

\newcommand{\Mod}{\operatorname{-Mod}}

\newcommand{\GT}{\operatorname{GT}}

\newcommand{\bpr}{\begin{proof}}
\newcommand{\epr}{\end{proof}}

\newcommand{\reg}{\operatorname{reg}}

\newcommand{\spec}{\operatorname{Spec}}
\newcommand{\ext}{\operatorname{Ext}}

\newcommand{\lra}{\longrightarrow}

\newcommand{\ra}{\rightarrow}

\newcommand{\mc}{\mathcal}
\newcommand{\mf}{\mathfrak}
\newcommand{\cha}{\operatorname{char}}
\newcommand{\mb}{\mathbb}

\newcommand{\pgl}{\operatorname{PGL}}

\newcommand{\m}{\mf{m}}
\newcommand{\p}{\mf{p}}
\newcommand{\diag}{\operatorname{diag}}

\renewcommand{\hom}{\operatorname{Hom}}

\newcommand{\wt}{\widetilde}
\newcommand{\ptfn}{\operatorname{\mathcal{S}_{\mathrm{pt}}}}
\newcommand{\qptfn}{\operatorname{\mathcal{S}_{\mathrm{q-pt}}}}

\newcommand{\coH}{\operatorname{H}}

\newcommand{\cd}{\operatorname{cd}}
\newcommand{\gld}{\operatorname{gld}}

\newcommand{\rcoh}{\operatorname{coh}}

\newcommand{\catmod}{\operatorname{-mod}}
\newcommand{\catMod}{\operatorname{-Mod}}
\newcommand{\rcatmod}{\operatorname{mod-\negmedspace}}
\newcommand{\rcatMod}{\operatorname{Mod-\negmedspace}}
\newcommand{\Gr}{\operatorname{-Gr}}

\newcommand{\qgr}{\operatorname{-qgr}}
\newcommand{\rGr}{\operatorname{Gr-\negmedspace}}
\newcommand{\rgr}{\operatorname{gr-\negmedspace}}
\newcommand{\rQgr}{\operatorname{Qgr-\negmedspace}}
\newcommand{\rqgr}{\operatorname{qgr-\negmedspace}}

\newcommand{\rTors}{\operatorname{Tors-\negmedspace}}
\newcommand{\rtors}{\operatorname{tors-\negmedspace}}

\newcommand{\id}{\operatorname{id}}
\newcommand{\aut}{\operatorname{Aut}}

\newcommand{\supp}{\operatorname{supp}}

\providecommand{\qedhere}{\qed}
\newboolean{ams2.0}
\ifthenelse{\equal{\qedhere}{\qed}}%
    {\setboolean{ams2.0}{false}}%
    {\setboolean{ams2.0}{true}}%

\DeclareMathOperator{\calHom}{\mathcal{H}\mathit{om}}
\DeclareMathOperator{\calExt}{\mathcal{E}\mathit{xt}}
\DeclareMathOperator{\calTor}{\mathcal{T}\mathit{or}}
\DeclareMathOperator{\Hom}{Hom}
\DeclareMathOperator{\Ext}{Ext}

\DeclareMathOperator{\Aut}{Aut}

\DeclareMathOperator{\character}{char}

\newcommand{\op}{\text{op}}

\newcommand{\ZZ}{{\mathbb Z}}

\newcommand{\NN}{{\mathbb N}}
\newcommand{\PP}{{\mathbb P}}

\newcommand{\calC}{{\mathcal C}}
\newcommand{\calR}{{\mathcal R}}

\newcommand{\calI}{{\mathcal I}}
\newcommand{\calJ}{{\mathcal J}}

\newcommand{\calN}{{\mathcal N}}
\newcommand{\calM}{{\mathcal M}}
\newcommand{\calO}{{\mathcal O}}

\newcommand{\OO}{{\mathcal O}}
\newcommand{\B}{{\mathcal B}}

\begin{document}
\title[Naive Noncommutative Blowing
 Up]{Na{\"\i}ve Noncommutative Blowing Up}
\author{D. S. Keeler}
\address{Department of Mathematics, MIT, 
 Cambridge, MA 02139-4307. } 
 \curraddr{Department of Mathematics,
Miami University, Oxford, OH~45056.}
\email{dskeeler@mit.edu}
\author{D. Rogalski}
\address{Department of Mathematics, University of Washington, 
Seattle, WA 98195-4350.} 
\curraddr{Department of Mathematics, MIT, 
 Cambridge, MA 02139-4307}
 \email{rogalski@math.washington.edu}
\author{J. T. Stafford}
\address{Department of Mathematics, University of Michigan, Ann Arbor,
MI 48109-1109.} 
\email{jts@umich.edu}
 \thanks{The 
first author was supported by an NSF Postdoctoral Research Fellowship.
 The second  author was supported in part by a Clay Liftoff Grant and by the
NSF through the  grant DMS-9801148 and  through
 an NSF Postdoctoral Research Fellowship.
 The third author was supported in part by the
NSF through the  grant DMS-9801148.}

\begin{abstract} Let $B(X,{\mathcal L},\sigma)$ be the twisted homogeneous
coordinate ring of an irreducible  variety $X$ over an algebraically 
closed  field $k$ with $\dim
X\geq 2$.  Assume that $c\in X $ and $\sigma\in \Aut(X)$
 are in sufficiently general position.
   We show that if one follows the commutative prescription for blowing
up $X$ at $c$, but in this noncommutative setting, one obtains a noncommutative
ring $R=R(X,c,{\mathcal L} ,\sigma)$ with surprising properties.  
 In particular:

\begin{enumerate}
\item[(1)] $R$ is always noetherian but never strongly noetherian.
\item[(2)] If $R$ is generated in degree one then the images of the $R$-point modules
in $\rqgr R$ are naturally in (1-1) correspondence with the closed points
of $X$. However, both in $\rqgr R$ and in $\rgr R$, the $R$-point modules are
not parametrized by a projective scheme. 
\item[(3)]  $\rqgr R$ has finite cohomological dimension yet 
$\dim_k {\coH}^1({\mathcal O}_R)= \infty$.
\end{enumerate}
This gives a more geometric approach to results of the second author who 
proved similar results for $X={\mathbb P}^n$ by algebraic methods.
\end{abstract}
 \maketitle 
 \tableofcontents 
 
\clearpage

\section{Introduction}\label{intro}
Noncommutative projective geometry has been very
successful in using the techniques and intuition of classical algebraic
geometry to understand noncommutative connected graded algebras $R=k\oplus
\bigoplus_{n\geq 1}R_n$, over an algebraically closed field $k$. 
In this paper we show
that one of the simplest noncommutative analogues
of blowing up a commutative variety
leads to algebras with a range of interesting properties.
 
One reason why   classical techniques work is that one can construct nontrivial
noncommutative graded rings from commutative data. A typical example is  the
following. Fix an automorphism $\sigma$ of a projective $k$-scheme  $X$ and for
an invertible  sheaf $\mathcal L$ write  
${\mathcal L}_n = {\mathcal L}\otimes_{{\mathcal
O}_X}\sigma^*{\mathcal L}\otimes_{{\mathcal O}_X}\cdots \otimes
(\sigma^{n-1})^* {\mathcal L}$, 
with ${\mathcal L}_0={\mathcal O}_X$. Define the \emph{bimodule algebra}  
${\mathcal B} = {\mathcal B}(X,{\mathcal
L},\sigma)  = \bigoplus_{n\geq 0} {\mathcal L}_n$ with global sections the \emph{
twisted homogeneous coordinate ring} $B(X,{\mathcal L},\sigma) =\bigoplus_{n\geq 0}
{\coH}^0(X, {\mathcal L}_n)$.  If 
${\mathcal L}$ is a $\sigma$-ample 
invertible sheaf, as defined in \eqref{sigma-ample} below, then $B$ is a noetherian ring and the
category $\rqgr B$ of finitely generated right $B$-modules modulo torsion is
equivalent to both $\rqgr \mathcal B$ and to $\rcoh X$, the category of
coherent sheaves on $X$.

In many important examples  (for example, domains  of Gelfand-Kirillov
dimension two generated in degree one \cite{AS}, Artin-Schelter regular rings
of dimension three \cite{ATV} and various higher dimensional   algebras) a
connected graded $k$-algebra $A$ has  such a ring  $B(X,{\mathcal L},\sigma)$
as a factor. One can then use the geometry of $X$ to understand $A$ and show
that it has very pleasant properties. The reason this works is that the  \emph{
point modules} over $A$, cyclic $A$-modules $M=\bigoplus_{i\geq 0}M_i$ with
$\dim_k M_i=1$ for all $i\geq 0$,  are parametrized by  the scheme $X$.
Combined with the fact that the shift functor $\sigma: M\mapsto M[1]_{\geq 0}$
induces an automorphism on the isomorphism classes of $A$-modules, this quickly
leads to the construction of the factor ring  $B(X,{\mathcal L},\sigma)$. It is
important to understand  how generally this technique applies. Specifically, it
had been hoped that these and related  geometric techniques would lead to a
classification of noncommutative surfaces, by which we mean $\rqgr A$
for a connected graded algebra $A$ with Gelfand-Kirillov dimension three.
A general survey of this programme can be found in \cite{SV}.

Recently the second author \cite{Ro} constructed examples of noncommutative
surfaces that do not have such pleasant properties. For example, although these
algebras  $R$ are noetherian, they are never \emph{strongly noetherian} (in other
words there exists a commutative noetherian $k$-algebra $C$  such that
$R\otimes_kC$ is not noetherian) and
their point modules do not appear to be parametrized by a projective scheme.  
The methods of \cite{Ro} are algebraic. The main aim of this paper is to give
an alternative, geometric construction  of these algebras that works  more
generally and helps explain their properties.
 
Our construction uses a natural noncommutative analogue of blowing up a
closed  point $c$ on an irreducible variety $X$. To set this in context,
consider first the classical approach. Write ${\mathcal I}={\mathcal I}_c\subset
{\mathcal O}_X$ for the ideal sheaf corresponding to $c$, let ${\mathcal L}$ be
any invertible sheaf on $X$ and form the sheaf of graded algebras  
${\mathcal
A} = {\mathcal O}_X \oplus \mc{I} \mc{L} \oplus 
\mc{I}^{2}\mc{L}^{\otimes 2}\oplus
\cdots.$ Then $\rqgr{\mathcal A}\simeq \rcoh \widetilde{X}$,  where
$\widetilde{X}$ denotes the blowup of $X$ at $c$.  
The sheaf $\mathcal L$ is
irrelevant to this construction but will be useful later.

Now consider noncommutative analogues of this construction. Using the
definition of $\mc{B}$ as a guide, it is natural to twist the summands of
$\mc{A}$ by powers of $\sigma$. Thus, 
write $\mc{I}_n = \mc{I}\cdot\sigma^*(\mc{I})\cdots(\sigma^{n-1})^*\mc{I}$, set
$\mc{J}_n=\mc{I}_n\mc{L}_n$
 and consider 
${\mathcal R} ={\mathcal R}(X,c,\mc{L},\sigma) = \bigoplus_{n\geq 0} 
\mc{J}_n\subset \mc{B}(X,\mc{L},\sigma)$  
in the place of $\mathcal A$. By analogy with the classical situation, we
define $\rqgr{\mathcal R}$ to be  \emph{the na\"\i ve
noncommutative blowup} of $X$ at~$c$. Although the category  $\rqgr \mc{R}$ 
is independent of $\mc{L}$, the algebra of global sections
$$R=R(X,c,\mc{L},\sigma) ={\coH}^0(X,{\mathcal R})=
\bigoplus_{n\geq 0}{\coH}^0(X,\mc{J}_n)$$ 
obviously does depend upon $\mc{L}$. 
When $\mc{L}$ is very ample, 
this will be the algebra that interests us.

Before stating the main result, we need one more definition. 
 If $c\in X$ is a closed point, 
 the orbit ${\mathcal C}=\{\sigma^n(c): n\in {\mathbb
Z}\}$ is defined 
to be \emph{critically dense} if the Zariski closure of 
any infinite subset ${\mathcal C}'$ of ${\mathcal C}$ equals $X$.
One can regard this as saying that $(c,\sigma)$ is in sufficiently general
position (see Section~\ref{section-critical} or 
\cite[Theorems~14.5 and 14.6]{Ro}).
Under this hypothesis, $R$ has a range of surprising
properties:

\begin{theorem}\label{mainthm} 
Let $X$ be an irreducible  variety with $\dim X\geq 2$ and
$\sigma\in \Aut(X)$. Assume that $\mathcal L$ is a very ample,
$\sigma$-ample invertible sheaf on $X$ and that 
$c\in X$ is a closed point 
 such that ${\mathcal C}=\{\sigma^n(c): n\in {\mathbb
Z}\}$ is critically dense. 

 If $\mc{R}=\mc{R}(X,c,\mc{L},\sigma)$ with global
sections $R=R(X,c,\mc{L},\sigma)$, then:
\begin{enumerate}

\item  $\rqgr R \simeq \rqgr {\mc{R}}$ and $\rqgr \mc{R}$ is independent
of the choice of $\mathcal L$.

\item $R$ is always noetherian.

\item $R$ is never strongly noetherian.

\item  The  simple objects
 in $\rqgr R$ are in (1-1) correspondence with the closed points of $X$. 

\item When $R$ is generated in degree one, the  simple objects in $\rqgr R$ are 
the images of the $R$-point modules. However, in both $\rqgr R$ and $\rgr R$, 
the $R$-point modules are not parametrized by any scheme of locally finite
type.

\item   $\rqgr R$ has finite cohomological dimension. When $X$ 
is smooth, $\rqgr R$ has finite
homological dimension.

\item If ${\coH}^1(R) ={\Ext}^1_{\rqgr R}(R,R)$, 
then $\dim_k {\coH}^1(R)= \infty$.

\item $R$ satisfies $\chi_1$ but does not satisfy $\chi_2$, as 
defined below.
\end{enumerate}
\end{theorem}
This theorem summarizes most of our results and so its proof takes
up most of  the paper.  Specifically, parts~1 and  2 of the theorem 
are proved in Theorem~\ref{general-ampleness} and
their proof takes up most of
Sections~\ref{definitions}--\ref{ampleness-in-general}.
The rest of the paper is then concerned with applying this theorem to get a
deeper understanding of the properties of $\mc{R}$ and $R$. In particular, 
part3 of Theorem~\ref{mainthm} is proved 
in Theorems~\ref{not strong noeth} and \ref{non-strong};
part~4 in Theorem~\ref{GT equiv};
part~5 in Theorem~\ref{non-represent}, Corollary~\ref{non-represent3} and
Remark~\ref{point-class12};
part~6 in Theorem~\ref{fin-cohom-dim} and Corollary~\ref{gldim}; and
parts~7 and 8  in Theorem~\ref{chi_1, not chi_2}.

In the special case where $X={\mathbb P}^n$ and $\mc{L} = \mc{O}(1)$, 
most  parts of the theorem were proved by more algebraic methods in
\cite{Ro}, which in turn was motivated by Jordan's work on  algebras
generated by Eulerian derivatives \cite{Jo}.
The significance of \cite{Ro} was to give counterexamples to a number of open
problems from the literature and Theorem~\ref{mainthm} obviously gives further
 examples. More significantly, it shows that these examples are to be
expected within the geometric framework of noncommutative geometry 
and perhaps also suggests a way of coping with them within, say, the
classification of noncommutative surfaces: if one regards these examples as a
form of noncommutative blowup then one may hope to classify such algebras using
noncommutative analogues of blowing up and down.

Let us explain the significance of individual parts of Theorem~\ref{mainthm}.
Part~1  justifies the idea that $\mathcal R$ is kind of a noncommutative 
blowup of $X$ at $c$. We begin by discussing this aspect of  the theorem since
it illustrates the ways in which  $\rqgr\mathcal R$ is both similar and
dissimilar to $\rcoh \widetilde{X}$.  In the commutative case, if
$\mc{L}=\OO_X$, then ${\mathcal I}_c{\mathcal A}$ corresponds to  an invertible
sheaf on $\widetilde{X}$ and the point $c$ is the  image of the exceptional
divisor  ${\mathcal A}/ {\mathcal I}_c{\mathcal A}\in \rcoh \widetilde{X}$.  In
contrast, an easy computation (Proposition~\ref{blowup2}) shows that  
${\mathcal R}/ {\mathcal I}_c{\mathcal R}$ is a (finite direct sum of copies) 
of a  simple object $\widetilde{c}\in \rqgr \mc{R}$. This can be used to prove
part~4 of the theorem. Ironically, $\widetilde{c}$ does have some properties
that  are more like a divisor than a point; in particular ${\mathcal
I}_c{\mathcal R}$ is still an invertible module in $\rqgr \mc{R}$ (see
Proposition~\ref{invertibility1}). Moreover, if one quotients out the Serre
subcategories ${\mathcal S}$ generated by  the modules corresponding to the
points $\sigma^i(c)$ in the two categories then, as should be expected by
analogy with blowing up a commutative variety, the quotient categories
$(\rqgr\mc{R})/{\mathcal S} $ and $(\rcoh X)/{\mathcal S}$ are equivalent
(Proposition~\ref{factor-equiv}). Remarkably, and in marked contrast to the
commutative situation, the subcategory of torsion sheaves in $\rcoh X$ is also
equivalent to the corresponding subcategory of $\rqgr\mc{R}$, the category of
Goldie torsion modules (see Theorem~\ref{GT equiv}). Thus, in some respects, the
difference between $\rqgr \mc{R}$ and $\mathrm{coh}\, X$ is quite subtle.

There is another version \cite{VB2} of noncommutative blowing up
that has properties much closer to the classical case and has been  
useful in describing noncommutative surfaces (see, for example,
\cite[Section~13]{SV}). Van den Bergh's construction is 
discussed briefly in Section~\ref{sect-blowup}.

The idea of considering strongly noetherian algebras arises in the work of
Artin, Small and Zhang \cite{ASZ,AZ2} who show that many algebras have this
property, and that it has a number of important consequences for an algebra.  
Notably, a strongly noetherian graded $k$-algebra $A$
satisfies \emph{generic flatness} in the following sense:
 for any  finitely generated commutative $k$-algebra
$C$ and any finitely generated  $A\otimes_kC$-module $M$ there exists  $f\in
C\smallsetminus \{0\}$ such that $M[f^{-1}]$ is a flat $C[f^{-1}]$-module
\cite[Theorem~0.1]{ASZ}. In contrast, $R$ fails generic flatness in a rather
dramatic way: 

\begin{proposition}\label{main-thm2}
 {\rm (}Theorem~\ref{not strong noeth}{\rm)}
If $U$  is \emph{any}  open affine subset of $X$,
 then generic flatness fails for the finitely generated $R\otimes {\mathcal
O}_X(U)$-module  ${\mathcal R}(U) =\bigoplus {\mathcal J}_n(U)$.
\end{proposition}

As was noted earlier, point modules have frequently been used to understand
specific classes of noncommutative algebras and \cite{ASZ,AZ2} use generic
flatness  to provide strong structure theorems for these modules (among
others). In particular, if $A$ is a strongly noetherian graded algebra
generated in degree one then the point modules for $A$, both in $\rgr A$ and
$\rqgr A$, are naturally parametrized by a projective scheme (see
\cite[Corollary~E4.5]{AZ2}, respectively Proposition~\ref{shift-aut} below). 
Moreover,  the shift functor $M\mapsto M[1]_{\geq 0}$ induces an automorphism
of this scheme (see Proposition~\ref{shift-aut}, again). 

All three of these
results fail for $R$ (see Section~\ref{section-point}).  This proves
Theorem~\ref{mainthm}(5) and is in marked contrast to part~4 of that result. 
The reason for the dichotomy is that, if one wants to parametrize 
the point modules  in $\rgr R$ or $\rqgr R$
 then, by definition, one needs to parametrize them simultaneously over all 
base spaces; that is, over $R_C=R\otimes_kC$ for all commutative $k$-algebras
$C$. However, Proposition~\ref{main-thm2} can be interpreted as saying
 that $R_{\OO(U)}$ has too few point modules in comparison to its localizations 
$R_{\OO_{X,p}}$ at closed points $p\in X$
for such a parametrization to be possible. (See Theorem~\ref{non-represent}
and Corollary~\ref{non-represent3}  for the precise statement.)

The $\chi$ conditions in part~8 of
Theorem~\ref{mainthm}  are defined as follows. A connected graded algebra $A$
satisfies $\chi_n$ if $\dim {\Ext}^i_{A\operatorname{-Mod}}(k,M)<\infty$ for
all finitely  generated graded $A$-modules $M$ and all $i\leq n$. These
conditions are central to the interplay between $R$ and $\rqgr R$ as described 
in \cite{AZ1}. In particular $\chi_1$ ensures that one can 
(essentially) recover the category
of finitely generated graded right $R$-modules $\rgr R$ from $\rqgr R$. The
higher $\chi$ conditions are related to more subtle properties of $\rqgr R$,
especially the behavior of cohomology.  In fact the failure of $\chi_2$ for $R$
is immediate from the fact that $\mathrm{H}^1(R)$ is infinite dimensional
(Theorem~\ref{mainthm}(7)).  This result is in contrast with a basic theorem
of Serre: $\mathrm{H}^i(Y,\mathcal{F})$ is 
 finite dimensional for any coherent sheaf $\mathcal{F}$ over a projective
 variety $Y$ and any $i>0$ \cite[Theorem~III.5.2]{Ha}. 
 Analogues of this result have 
 also been basic to much of the theory of noncommutative geometry, as developed 
  for example in 
 \cite{AZ1,AZ2},  and so that theory is not available for the study of $R$.

 We would like to thank Michael Artin and Brian Conrad for their help and
 suggestions with this paper, especially with the material from
 Section~\ref{section-point}.  We would also like to thank James Zhang and Paul
 Smith for helpful conversations.

\section{Definitions and background material}\label{definitions}

As was  mentioned in the introduction, we want to work with 
bimodule algebras like ${\mathcal B}=\bigoplus {\mathcal L}_n$ 
and in this section we set up the appropriate notation. Most of this 
comes from   \cite{AV} and \cite{VB1} and the 
reader is referred to those papers  for further  details.

Fix throughout an integral projective scheme $X$ over an algebraically closed
field $k$.   The category of quasicoherent, respectively coherent, sheaves on
$X$ will be written $\mc{O}_X\Mod$, respectively $\mc{O}_X\catmod$. We use the
following notation  for pullbacks: if $\sigma\in \Aut(X)$
is a $k$-automorphism of   $X$,  and $\mc{F}\in \mc{O}_X\catmod$, then
  $\mc{F}^{\sigma}=\sigma^{*}(\mc{F})$. 
We adopt the usual convention that an automorphism $\sigma$ 
acts on functions by $f^\sigma(x) = f(\sigma(x))$, for $x\in X$. 

\begin{definition}
\label{bimod def} 
A \emph{coherent $\mc{O}_X$-bimodule} is a coherent sheaf $\mc{F}$
on  $X \times X$ such that  $Z = \supp{\mc{F}}$ has the property that both
projections  $\rho_1, \rho_2: Z \ra X$ are finite morphisms.
An \emph{$\mc{O}_X$-bimodule} is a  sheaf $\mc{F}$
on  $X \times X$ such that every coherent $X \times X$-subsheaf is a 
coherent $\mc{O}_X$-bimodule.  The left and right $\mc{O}_X$-module
structures associated to $\mc{F}$ are defined 
to be $_{\mc{O}_X} \mc{F} = (\rho_1)_{*} \mc{F}$ and 
$\mc{F}_{\mc{O}_X} = (\rho_2)_{*} \mc{F}$ respectively. 
 \end{definition}

The tensor product of two bimodules is again a bimodule and satisfies 
the expected properties; for example the
tensor product is right exact and associative (see \cite[Section~2]{VB1}). 
In fact we will not be concerned with this generality since all 
the bimodules we
consider   arise from the following construction.

\begin{definition}
\label{standard bimodules}
 Let $\mc{F}\in \mc{O}_X\catmod$ and $\tau,\sigma\in \Aut(X)$. 
Then define an $\mc{O}_X$-bimodule ${}_\tau
\mc{F}_{\sigma}$ by $(\tau, \sigma)_*\mc{F}$ where
 $(\tau, \sigma): X \to X \times X$. 
We usually  write $\mc{F}_{\sigma}$ 
for ${}_1 \mc{F}_{\sigma}$, where $1$ is the identity automorphism. 
\end{definition}

We will see in the next lemma that we only need to consider bimodules
of the form ${}_1 \mc{F}_{\sigma}$.
The reader may check that such a bimodule has left $\mc{O}_X$-module structure
$\mc{F}$ but right $\mc{O}_X$-module structure $\mc{F}^{\sigma^{-1}}$.

When no other bimodule structure  is
given, a coherent sheaf $\mc{F}$ on $X$ will be  assumed to have the bimodule
structure $_1 \mc{F} _1$.  Thus all sheaves become bimodules, and all tensor
products can be thought of as tensor products of bimodules. 
In order to remove ambiguity,  when thinking of a bimodule $\mc{G}$ as a sheaf, we mean the left
$\mc{O}_X$-module structure of $\mc{G}$, unless otherwise stated. In
particular, when we write $\coH^i(X,\mc{G})$ or say that $\mc{G}$ is generated by
its global sections we  are referring to the left structure of $\mc{G}$.
We often write $\Gamma(\mc{G})$ for $\coH^0(X,\mc{G})$.
Working on the left will have notational advantages, but is otherwise not
significant since we have:

\begin{lemma}
\label{tensor bimodules} Let $\mc{F}$, $\mc{G}$ be coherent $\mc{O}_X$-modules, and 
$\sigma, \tau$ automorphisms of $X$.
\begin{enumerate}
\item\label{tensor bimodules 0}
${}_\tau \mc{F}_\sigma \cong {}_1 (\mc{F}^{\tau^{-1}})_{\sigma \tau^{-1}}$.

\item\label{tensor bimodules 1}
$\mc{F}_{\sigma} \otimes \mc{G}_{\tau} \cong (\mc{F} 
\otimes \mc{G}^{\sigma})_{\tau \sigma}$.

\item\label{tensor bimodules 2}
The vector space of 
global sections of ${}_{\OO_X}(_1\mc{G}_\tau)$  is naturally 
isomorphic to that of $(_1\mc{G}_\tau)_{\OO_X}$.
\end{enumerate}
\end{lemma}
\begin{proof} \eqref{tensor bimodules 0}
This follows from the comments after the diagram \cite[(2.2)]{VB1},
taking $X=Y=V=V''$.

 \eqref{tensor bimodules 1} 
This is a special case of \cite[Lemma~2.8(2)]{VB1}, where $X=Y=Z=V=W$.

\eqref{tensor bimodules 2} From the comments before the lemma, 
the right structure of $\mc{G}$ is just
$\mc{G}^{\tau^{-1}}$ and so this has global sections $\coH^0(X,
\mc{G}^{\tau^{-1}}) = \coH^0(X, \mc{G})^{\tau^{-1}}$.
 \end{proof}

We can now define bimodule algebras and their categories of modules.
Since we only need a special case of Van den Bergh's bimodule algebras,
 we will only make our definitions in that special case.

\begin{definition}\label{bimod-alg} 
An \emph{$\mc{O}_X$-bimodule algebra} is an
$\mc{O}_X$-bimodule $\mc{B}$ together with a unit map $1: \mc{O}_X \ra \mc{B}$
and a product map  $\mu: \mc{B} \otimes \mc{B} \ra \mc{B}$ satisfying the usual
axioms. Let $\sigma\in\Aut(X)$. 
 The bimodule algebra $\mc{B}$ is called a
  \emph{graded $(\mc{O}_X, \sigma)$-bimodule
  algebra}
  if:
  \begin{enumerate}
  \item $\mc{B}$ decomposes as a direct sum $\mc{B} =
\bigoplus_{n\geq 0} \mc{B}_n$ of 
$\mc{O}_X$-bimodules  $\mc{B}_n\cong {}_1(\mc{E}_n)_{\sigma^n}$, for some 
$\mc{E}_n\in\ \mc{O}_X\catmod$ with 
$\mc{B}_0={}_1(\mc{O}_X)_1$.
\item The multiplication map satisfies $\mu
(\mc{B}_m \otimes \mc{B}_n) \subseteq \mc{B}_{m+n}$ for all $m,n$ 
 and $1(\mc{O}_X) \subseteq
\mc{B}_0$. Equivalently $\mu$ is defined by   $\OO_X$-module maps
$\mc{E}_n\otimes \mc{E}_m^{\sigma^n}\to \mc{E}_{m+n}$ satisfying the
appropriate associativity conditions.
\end{enumerate}
We will write $\mc{B} =
\bigoplus{}_1(\mc{E}_n)_{\sigma^n}$ throughout the section. 
\end{definition}

\begin{definition} 
\label{graded def} 
Let $\mc{B}$ be a graded $(\mc{O}_X, \sigma)$-algebra.  A 
\emph{right $\mc{B}$-module} $\mc{M}$ is a quasi-coherent
right $\mc{O}_X$-module together with a right $\mc{O}_X$-module  map $\mu:
\mc{M} \otimes \mc{B} \ra \mc{M}$  satisfying the usual axioms. The module
$\mc{M}$ is \emph{graded} if $\mc{M} = \bigoplus_{n \in \mb{Z}} \mc{M}_n$ with 
$\mu(\mc{M}_n\otimes \mc{B}_m ) \subseteq \mc{M}_{m+n}$. The \emph{shift of}
$\mc{M} $ is defined by $\mc{M}[n]=\bigoplus \mc{M}[n]_i$ with 
$\mc{M}[n]_i=\mc{M}_{i+n}$. 

The $\mc{B}$-module $\mc{M}$ is \emph{coherent} (as a  $\mc{B}$-module) if
there is a coherent $\mc{O}_X$-module $\mc{M}_0$  and a surjective map
 $\mc{M}_0 \otimes \mc{B} \ra \mc{M}$ of ungraded 
$\mc{B}$-modules.  Left $\mc{B}$-modules are
defined similarly and  the bimodule algebra $\mc{B}$ is \emph{right (left)
noetherian} if every   right (left) ideal of $\mc{B}$ is coherent.  
For the algebras that interest us, a more natural definition of 
coherence is given in Lemma~\ref{coherent}.
\end{definition}

A priori a graded right
$\mc{B}$-module $\mc{M}=\bigoplus \mc{M}_i$
 is only a right $\mc{O}_X$-module.  One can obviously give 
 the $\mc{M}_i$ various different bimodule structures and 
 we choose the one that is most convenient.
    Specifically, it will cause no loss of generality to
assume that all right $\mc{B}$-modules  have the form 
\begin{equation}\label{nice-modules}
\mc{M} = \bigoplus_{n
\in \mb{Z}} {}_1 (\mc{G}_n)_{\sigma^n}
\text{ for some (left) sheaves }\mc{G}_n\in\mc{O}_X\Mod.
\end{equation}
The advantage of this choice is that the $\mc{B}$-module structure on $\mc{M}$ 
is given by a family of $\mc{O}_X$-module  maps
$\mc{G}_n \otimes \mc{E}_m^{\sigma^n} \to \mc{G}_{n+m}$,
 again satisfying the appropriate associativity conditions.

The graded right $\mc{B}$-modules form an abelian category $\rGr \mc{B}$, with
homomorphisms graded  of degree zero.  Its
subcategory of coherent modules is denoted $\rgr \mc{B}$.   A \emph{bounded}
graded $\mc{B}$-module $\bigoplus \mc{M}_i$ is one such that $\mc{M}_i = 0$ for
all but finitely many $i$. A module $\mc{M} \in \rGr\mc{B}$ is called
\emph{torsion} if every  coherent submodule of $\mc{M}$ is bounded.  Let 
$\rTors \mc{B}$ denote  the full subcategory of $\rGr \mc{B}$ consisting of  
torsion modules, and write $\rQgr \mc{B}$ for  the quotient category $\rGr
\mc{B}/\negthinspace\rTors \mc{B}$. The analogous quotient category of $\rgr \mc{B}$
 will be denoted  $\rqgr \mc{B}$. The
 corresponding  categories of left modules will be  denoted by $\mc{B}\Gr$, etc.

Similar category definitions apply to graded rings.  A graded $k$-algebra 
$A=\bigoplus_{i\geq 0}A_i$ is called \emph{connected graded} if $A_0=k$.
If $A$ is noetherian then  a \emph{torsion} right $A$-module is a graded module
such that every finitely generated submodule is bounded. Write $\rGr A$ for the
category of graded right $A$-modules,  with torsion subcategory $\rTors A$ and 
quotient category $\rGr A/\negthinspace\rTors A$. 
Similarly, write $\rgr A$ for the category
of finitely generated right $A$-modules with  quotient category $\rqgr A=\rgr
A/\negthinspace\rtors A$. 
We denote the natural quotient maps by
\begin{equation}\label{quotient}
\pi_{\mc{B}} : \rGr \mc{B}\to \rQgr\mc{B}\qquad\text{and}\qquad
\pi_A : \rGr A \to \rQgr A
\end{equation}
and write both maps as $\pi$ if no confusion is possible.

We will almost always prove results for right modules over rings or bimodule
algebras. By the next lemma, these results will then have a natural
counterpart on the left.

\begin{definition}\label{opposite bimodule definition}
Let $\psi: X \times X \to X \times X$ be the automorphism given by 
$(x, y) \mapsto (y,x)$.
Let $\mc{B}$ be a  graded $(\mc{O}_X, \sigma)$-algebra. Then 
the \emph{opposite bimodule algebra}
is defined to be  $\mc{B}^\op=\psi_* \mc{B}$.  The unit and product map for $\mc{B}^{\op}$ are 
induced by $\psi_{*}$ from the unit and product of $\mc{B}$.
\end{definition}

\begin{lemma}\label{opposite bimodule lemma}
Let $\mc{B}=\bigoplus{}_1(\mc{E}_n)_{\sigma^n}$ be a  graded 
$(\mc{O}_X, \sigma)$-algebra.
Then $$\mc{B}^\op \cong \bigoplus{}_{\sigma^n}(\mc{E}_n)_1 \cong
\bigoplus {}_1(\mc{E}^{\sigma^{-n}}_n)_{\sigma^{-n}}.$$
Thus $\mc{B}^\op$ is a graded $(\mc{O}_X, \sigma^{-1})$-algebra.
There is a natural category equivalence $\rGr \mc{B} \simeq \mc{B}^\op \Gr$ 
(and similarly for the other module categories).
\end{lemma}
\begin{proof} That $\mc{B}^\op \cong \bigoplus{}_{\sigma^n}(\mc{E}_n)_1$
 follows immediately from
Definition~\ref{standard bimodules} while the second isomorphism 
is just Lemma~\ref{tensor bimodules}\eqref{tensor bimodules 0}.
If $\psi: X \times X \to X \times X$ is the automorphism given by 
$(x, y) \mapsto (y,x)$,
then $\psi_*$ naturally induces the equivalences of categories.
\end{proof}

The notion of coherence for $\mc{B}$-modules should be viewed as an analog of
finite generation, but there is a subtlety here: even if $\mc{B}$ is right
noetherian,  it does not seem to follow that every submodule of a coherent
$\mc{B}$-module  is coherent!  Fortunately, all the bimodule algebras we need
are covered by the following result, so there is no problem.

\begin{proposition} \label{abelian}
Let $\B=\bigoplus \mc{B}_i$  be a right noetherian 
graded $(\OO_X, \sigma)$-bimodule algebra and write $\mc{B}_i=
 {}_1(\mc{E}_i)_{\sigma^i}$ for each $i$.  Assume that
each  $\mc{E}_i$ is a subsheaf of a 
locally free sheaf $\mc{L}_i$ on $X$.  
If $\calM \in \rgr\B$, then every $\B$-submodule
of $\calM$ is coherent. 

Thus $\rgr \B$ is an abelian category and a right
$\B$-module $\calN$ is noetherian if and only if it is coherent.
\end{proposition} 

\begin{proof}
 The only step of substance will be to show that submodules of coherent
$\mc{B}$-modules are coherent. As in the proof of 
\cite[Proposition~3.6(3)]{VB1},  
let $\calC$ denote  the class of $\calM_0 \in \OO_X\catmod$ such that every
submodule of  $\mc{M}_0 \otimes \B$ is coherent. We first  wish to show that
$\calC = \OO_X\catmod$.

Pick a very ample invertible sheaf $\OO_X(1)$ over $X$ and a nonzero global
section $z\in \coH^0(X,\OO_X(1))$. Since $X$ is integral, multiplication
by  $z^n$, for $n\geq 1$, induces an injection $\OO_X\hookrightarrow \OO_X(n)$
and hence an injection $\OO_X(-n)\otimes \mc{L}_i \hookrightarrow \mc{L}_i$. 
As $\mc{E}_i\subseteq \mc{L}_i$ and $\OO_X(-n)$ is locally free, we therefore
get an induced  embedding $\OO_X(-n)\otimes \mc{B}_i \hookrightarrow
\mc{B}_i$.  Thus the  right $\B$-module morphism $ \OO_X(-n)  \otimes
\B  \to \B$ induced by (left) multiplication by $z^n$ is an 
injection for all $n\geq 1$. Therefore,
 any $\B$-submodule of $\OO_X(-n) \otimes \B$ is isomorphic  to a right
$\B$-ideal and hence is coherent. Hence $\OO_X(-n) \in \calC$ for all $n\geq 0$.

As in the proof of \cite[Proposition~3.6(1--2)]{VB1}   it  is easy to see that
quotients and extensions of coherent right $\B$-modules are coherent and so
$\calC=\OO_X\catmod$. By definition, this implies that any $\mc{B}$-submodule
of a coherent  module is coherent.
Therefore, the coherent $\mc{B}$-modules
form an abelian category.   Finally, the proof of
\cite[Proposition~3.6(4)]{VB1}  shows that a $\B$-module is  noetherian if and
only  if it is coherent. 
  \end{proof}

As noted in the introduction we are also interested in the algebra of sections
of a bimodule algebra. This is defined as follows. Fix a graded
$(\OO_X,\sigma)$-algebra $\mc{B}=\bigoplus \mc{B}_n$ with 
$\mc{B}_n = (\mc{E}_n)_{\sigma^n}$ for each $n$. 
Then it is clear that the space of (left) global sections 
$B = \Gamma(\mc{B})=\coH^0(X,\mc{B})$  has a natural graded $k$-algebra
structure $B=\bigoplus \Gamma(\mc{B}_n)$, given by the maps 
$\Gamma(\mc{E}_n) \otimes \Gamma(\mc{E}_m^{\sigma^n}) 
\to \Gamma(\mc{E}_{n+m})$
for all $n, m \geq 0$.  We call $B$ the \emph{section algebra} of $\B$. 
 Similarly, if $\mc{M}$ is a graded right $\B$-module, write 
$\mc{M} = \bigoplus_{n \in \mb{Z}} (\mc{G}_n)_{\sigma^n}$ 
as in \eqref{nice-modules} and  define 
$$\Gamma(\mc{M})=\bigoplus \Gamma(\mc{G}_n)=\bigoplus \coH^0(X,\mc{G}_n).$$
This  is naturally a
right $B$-module, given by   maps
$\Gamma(\mc{G}_n) \otimes \Gamma(\mc{E}_m^{\sigma^n}) \to \Gamma(\mc{G}_{n+m})$
for all $n \in \mb{Z}$, $m \geq 0$.
Thus we get a functor $\Gamma: \rGr \mc{B} \ra \rGr B$.

Conversely,   for a right $B$-module $M$ we define
$  M \otimes_B \mc{B}$ to be the sheafification of the presheaf $V \mapsto
 M \otimes_B \mc{B}(V) $ for open 
 $V\subseteq X$ .  One may check that $  M \otimes_B \mc{B}$ is naturally a
right  $\mc{B}$-module.  Analogously, one may define a left
$\mc{B}$-module  $\mc{B} \otimes_B N$ for any left $B$-module $N$. 
The functor $-\otimes \mc{B}: \rGr B \ra \rGr \mc{B}$ is
 a right adjoint to $\Gamma:  \rGr
\mc{B} \to \rGr B$.

The next definition and theorem from \cite{VB1} provide an important 
situation when one can relate the properties of $\mc{B}$ and $B$.

\begin{definition}
\label{ample def}
Suppose that  $\{ \mc{J}_n \}_{n \in \mb{N}}$ 
is a sequence of $\mc{O}_X$-bimodules.
  Then 
   the sequence is \emph{(right) ample} if,
for any $\mc{M} 
\in \OO_X\catmod$, one has the following:
\begin{enumerate}
\item $\mc{M} \otimes_{\mc{O}_X} \mc{J}_n$ 
is generated by global sections for $n \gg 0$.
\item $\coH^i(\mc{M} \otimes \mc{J}_n) = 0$ for all $i > 0$ and $n \gg 0$.
\end{enumerate}
\end{definition}

\begin{theorem} {\rm (Van den Bergh)}
\label{VdB main theorem}
Let $\mc{B}=\bigoplus \mc{B}_i$ be a graded  $(\mc{O}_X,\sigma)$-algebra.
Assume that  $\mc{B}$ is right noetherian  and that 
  $\{ \mc{B}_n \}_{n \in \mb{N}}$ is a
right ample sequence of $\mc{O}_X$-bimodules such that each $\mc{B}_n$ 
is contained in a locally free left $\mc{O}_X$-module. Then the
section algebra $B = \Gamma(\mc{B})$ is right noetherian, and 
there is an equivalence of categories $\xi : \rqgr \mc{B} \simeq \rqgr B$ 
via the
inverse equivalences $\Gamma(-)$ and $- \otimes_B \mc{B}$. 
\end{theorem}

\begin{proof} Since $\rgr\mc{B}$ is abelian (Proposition~\ref{abelian})  this
follows from the right hand version of \cite[Theorem~5.2]{VB1}, with one
proviso.  The conventions for the right handed version of
 \cite{VB1} require that one defines
$M=\Gamma(\mc{M})$ for $\mc{M}=\bigoplus \mc{M}_i\in \rgr \mc{B}$ by taking the
sections of $\mc{M}_i$ as a right $\OO_X$-module. 
However, the precise construction of $B$ and $M$ is 
not important to that proof; all one requires is that 
 the module structure of $M$ is induced from that of $\mc{M}$.
Thus the proof in \cite{VB1} also proves this theorem. 

Alternatively, it is also not hard to check using 
Lemma~\ref{tensor bimodules}\eqref{tensor bimodules 2}
 that  the left and right section algebras of $\mc{B}$
are actually isomorphic as graded rings. \end{proof}

The first main result of the paper (see Theorem~\ref{general-ampleness} or
part~1 of Theorem~\ref{mainthm}) will be to show that 
this theorem can be applied to our na{\"\i}ve blowups. 
 This will then allow us to identify  $\rqgr R$ with $\rqgr \mc{R}$
 for such algebras.

An important special case of Definition~\ref{ample def} and Theorem~\ref{VdB
main theorem} occurs when $\mc{J}_n = \mc{B}_n = (_1\mc{L}_\sigma)^{\otimes n}$
for an invertible sheaf $\mc{L}$ on $X$. We will usually write 
$\mc{L}_\sigma^{\otimes n}$ for 
$(_1\mc{L}_\sigma)^{\otimes n}$.
It is customary to say that
\begin{equation}\label{sigma-ample} 
\mc{L} \text{ is } \sigma\text{\emph{-ample if }}\
\{ \mc{L}_{\sigma}^{\otimes n} \}_{n \geq 0}\
\text{ is a right ample sequence of 
bimodules.}
\end{equation}
We will always write $\mc{B}=\mc{B}(X,\mc{L},\sigma) =
 \bigoplus \mc{L}_\sigma^{\otimes n}$
with    section algebra $$B=B(X,{\mathcal L},\sigma) =
 \bigoplus_{n\geq 0}
\Gamma(\mc{L}_{\sigma}^{\otimes n}) .$$
This is an equivalent definition of  the  \emph{twisted homogeneous
coordinate ring } of $X$ from the introduction.
The $\sigma$-ampleness condition is subtle since   
given a projective scheme $X$ with automorphism $\sigma$, there may be no 
$\sigma$-ample sheaves.  However, it is known when one such sheaf exists and in that case
all ample invertible sheaves are automatically $\sigma$-ample, on both the left
and the right.  For these and
further  results about $B(X,{\mathcal L},\sigma)$, see \cite{AV,Ke1}.

\section{Rees bimodule algebras}
\label{Rees bimod algs}

 In this section we formally define the algebras
$\mc{R}=\mc{R}(X,c,\mc{L},\sigma)$ and 
$R=\Gamma(\mc{R})$ from Theorem~\ref{mainthm} and  give
conditions under which $\mc{R}$ is noetherian. In the next section  we will
consider when the sequence $\{\mc{R}_n\}$ is ample in the sense of
Definition~\ref{ample def}. Once this has been done, Van den Bergh's Theorem
~\ref{VdB main theorem} can be applied to show that $\mc{R}$ and $R$ are 
noetherian in considerable generality.

The following assumptions will be fixed throughout the section.

\begin{assumptions}\label{global-convention}
Fix an integral projective scheme $X$. Fix
$\sigma\in\Aut(X)$, an invertible sheaf
$\mc{L}$ on $X$ and let $\mc{I}=\mc{I}_c$ denote
the sheaf of ideals defining a closed point $c$ on $X$.
  Assume that $c$ has infinite order 
under $\sigma$ and write $c_i = \sigma^{-i}(c)$ for $i \in \mb{Z}$.
Our  convention on automorphisms
from the beginning of Section~\ref{definitions} means that 
$ \mc{I}^{\sigma^i} = \mc{I}_{c_i}$ 
with quotient $\OO_X/\mc{I}_{c_i}= k(c_i)$ the corresponding 
skyscraper sheaf.
 \end{assumptions}

Mimicking classical blowing up we set 
$$\mc{I}_n =\mc{I} \mc{I}^{\sigma} \dots
\mc{I}^{\sigma^{n-1}},\quad
\mc{L}_n = \mc{L}\otimes \mc{L}^{\sigma} \otimes\dots\otimes 
\mc{L}^{\sigma^{n-1}}\quad\text{and}\quad
\mc{R}_n={}_1(\mc{I}_n\otimes \mc{L}_n)_{\sigma^n},$$
where all tensor products are over $\OO_X$.  From this 
data we define a bimodule algebra
 $$\mc{R} = \mc{R}(X,c,\mc{L},\sigma) = \bigoplus_{n =
0}^{\infty} \mc{R}_n$$ 
with corresponding algebra of sections
$$R = R(X,c,\mc{L},\sigma) =  \Gamma(\mc{R})=\bigoplus \Gamma(\mc{R}_n).$$ 

Note that by Lemma~\ref{opposite bimodule lemma},
\begin{equation}\label{rees opposite}
\mc{R}(X,c,\mc{L},\sigma)^\op =
 \mc{R}(X, \sigma(c), \mc{L}^{\sigma^{-1}}, \sigma^{-1}).
\end{equation}
Thus any results on the right can immediately be transferred to the left.

As the next two easy lemmas show, the hypotheses from \eqref{global-convention}
allow us to replace products by tensor products in the definition
of $\mc{I}_n$ and so there should be no confusion between this and 
the definition of $\mc{L}_n$. 

\begin{lemma}
\label{prod versus tensor}
For $i=1,2$, 
let $\mc{F}_i\in \OO_X\catMod$, and let $Z_i \subseteq X$ be 
the set of closed points of $X$ at 
which $\mc{F}_i$ is not locally free.
Then:
\begin{enumerate}
\item[(1)] If $\calTor$ denotes sheaf Tor then 
$\supp \calTor_j^{\mc{O}_X}(\mc{F}_1, \mc{F}_2) \subseteq Z_1 \cap Z_2$
 for $j > 0$.
 
\item[(2)]  If $\mc{K}  \subseteq \mc{O}_X$ 
  is an ideal sheaf with $\supp(\OO_X/\mc{K})\cap Z_2=\emptyset$, then  
$\mc{K} \otimes \mc{F}_2 \cong \mc{K}\mc{F}_2$. 
\end{enumerate} 
\end{lemma}

\begin{proof}
By \cite[p.~700]{GH}, it suffices to prove this locally, where it is obvious.
\end{proof}

\begin{lemma}
\label{structure J^n}  For 
$n\geq 1$, there are natural isomorphisms
 $\mc{I}_n\cong \mc{I}\otimes \cdots\otimes\mc{I}^{\sigma^{n-1}}$
and $\mc{I}_n\otimes\mc{L}_n\cong \mc{I}_n \mc{L}_n$.

If $\mc{R} = \mc{R}(X,c,\mc{L},\sigma)$, then
the natural homomorphism $\mc{R}_m\otimes \mc{R}_n\to \mc{R}_{m+n}$ is an
isomorphism for all $m,n\geq 0$. In particular, 
$\mc{R}_n  \cong (\mc{R}_1)^{\otimes n}$ as $\OO_X$-bimodules.
\end{lemma}
\begin{proof}
Since  $\OO_X/\mc{I}^{\sigma^i}\cong k(c_i)$ and  
 the points $c_i$ are distinct, the first result  follows 
 from Lemma~\ref{prod versus tensor}(2). The
  second assertion follows from this together with 
 Lemma~\ref{tensor bimodules}.
\end{proof}

The choice of $\mc{L}$ in Assumptions~\ref{global-convention}
is  important to the study of the section algebra
$R$   but, as in the classical case \cite[Lemma~II.7.9]{Ha},
it  is  irrelevant  to the study of
the  bimodule algebra $\mc{R}$.

\begin{proposition}
\label{L no matter} Given two invertible sheaves  $\mc{L}$ and $\mc{L}'$, write 
$\mc{R}=\mc{R}(X,c,\mc{L},\sigma)$, respectively 
$\mc{R}'=\mc{R}(X,c,\mc{L}',\sigma)$.
Then $\mc{R}\Gr \simeq \mc{R}'\Gr$ and $\rGr \mc{R} \simeq \rGr \mc{R}'$. 
\end{proposition} 
\begin{proof}
We prove the result on the right; the left-sided result follows 
from Lemma~\ref{opposite bimodule lemma}.  
It also suffices to prove the result when $\mc{L}'=\OO_X$.
As in \eqref{nice-modules},
we write an arbitrary module $\mc{M}'\in \mc{R}'\Gr$ as $\mc{M}' =
\bigoplus {}_1(\mc{G}_n)_{\sigma^n}$, so that the module structure  is given 
by homomorphisms of sheaves $\alpha_{n,m} : 
\mc{G}_n \otimes \mc{I}_m^{\sigma^n} \to \mc{G}_{n+m}$  for
all $n \in \mb{Z}$ and  $m \geq 0$.

We wish to construct a right $\mc{R}$-module 
$\mc{M} = \bigoplus {}_1(\mc{G}_n \otimes \mc{L}_n)_{\sigma^n}$ from $\mc{M}'$. 
To define the module structure requires maps of sheaves 
\[
\beta_{n,m} : (\mc{G}_n \otimes \mc{L}_n) \otimes (\mc{I}_m \otimes 
\mc{L}_m)^{\sigma^n} \to \mc{G}_{n+m} \otimes \mc{L}_{n+m}.
\]
However,  since 
$\mc{L}_n \otimes (\mc{L}_m)^{\sigma^n} \cong \mc{L}_{m+n}$, 
the maps $\beta_{n,m}$ arise naturally by tensoring the given maps
$\alpha_{n,m}$  with $\mc{L}_{m+n}$.  Checking that
 this does indeed  define a module structure on $\mc{M}$ is straightforward.

This gives a functor $\theta: \rGr \mc{R}' \to \rGr \mc{R}$.  
A similar argument will construct a functor $\psi: \rGr \mc{R} \to \rGr \mc{R}'$  
which sends an $\mc{R}$-module $\bigoplus {}_1(\mc{G}_n)_{\sigma^n}$ 
to the $\mc{R}'$-module $\bigoplus {}_1(\mc{G}_n \otimes
 \mc{L}_n^{-1})_{\sigma^n}$.  It
is obvious that $\psi \theta$ and $\theta \psi$ are naturally isomorphic to 
identity functors, so that $\rGr \mc{R}$ and $\rGr \mc{R}'$ are equivalent.  
\end{proof}

The main result of this section (Proposition~\ref{bimod alg noeth}) determines 
when the bimodule algebra $\mc{R}$ is noetherian.  
The answer will involve the following geometric notion.

\begin{definition}\label{crit-defn}
Let $\mc{C}$ be an infinite set of closed points of an integral scheme $X$. 
 Then
we say that $\mc{C}$ is \emph{critically dense} if every infinite subset
 of $\mc{C}$ has Zariski
closure equal to $X$.
\end{definition}

We start with some easy observations. 

\begin{lemma}
\label{smooth} If  $\mc{C}=\{c_n : n\in{\mathbb Z}\}$ is critically dense then 
$X$ is smooth at every point $c_i \in \mc{C}$.
\end{lemma}
\begin{proof}
If some $c_i$ lies in the non-smooth locus $Z$, then so does every point
$c_j=\sigma^{i-j}(c_i)$. Hence the closure $X$ of $\mc{C}$ is contained in $Z$, 
which is absurd.
\end{proof}

\begin{lemma}\label{bimod-lemma}
Let $\mc{U} ,\mc{V} \subseteq \OO_X$ be comaximal ideal sheaves. 
Assume that $\mc{W}$ is an ideal sheaf such that 
$\mc{U}\cap\mc{V} \subseteq \mc{W}\subseteq \mc{V}$. 
Then $\mc{Z}=\mc{U}+\mc{W}$ is the unique largest sheaf of ideals such that 
$\mc{Z}\mc{V}\subseteq \mc{W}$. Moreover, $\mc{Z}\mc{V}=\mc{Z}\cap \mc{V}
=\mc{W}$.
\end{lemma}

\begin{proof} This follows from a repeated use of comaximality.\end{proof}

The definition of a coherent $\mc{R}$-module from Definition~\ref{graded def}
is not convenient for most applications and so we will use the following
equivalent definition.

\begin{lemma}\label{coherent}
 The following are equivalent for a  module
  $\mc{M}=\bigoplus_{n\in \mathbb Z} \mc{M}_n\in \rGr \mc{R}$:
 \begin{enumerate} 
 \item[(1)] $\mc{M}$ is coherent.
 \item[(2)] Each 
$\mc{M}_n $ is a coherent $\OO_X$-module, with  $\mc{M}_n = 0$
 for $n \ll 0$, and 
the natural map $\mu_n:\mc{M}_n\otimes \mc{R}_1\to \mc{M}_{n+1}$ 
is surjective for   $n\gg 0$.
\end{enumerate}\end{lemma}

\begin{proof} If $\mc{M}$
 is coherent,  there is a surjection 
$\mc{F}\otimes \mc{R} \twoheadrightarrow
 \mc{M}$ for some $\mc{F}\in \OO_X\catmod$. 
Equivalently, for some $a \leq  b$ there is a graded surjection 
 $\Theta: \bigoplus_{m=a}^{b}\mc{F}_m\otimes \mc{R}[-m] 
 \twoheadrightarrow\mc{M}$, where each $\mc{F}_m$ is a copy of $\mc{F}$ 
situated in degree $m$.  Clearly  each $\mc{M}_n$ is then a coherent sheaf and 
$\mc{M}_n = 0$ for $n \ll 0$.
For $n\geq b$ we have a commutative diagram 
 \[ \begin{CD}
 \bigoplus_m\mc{F}_m\otimes \mc{R}_{n-m}\otimes \mc{R}_1 
 @>\theta_1>>\mc{M}_n \otimes \mc{R}_1 \\
 @V\phi_1VV @VV\phi_2V \\
  \bigoplus_m\mc{F}_m\otimes \mc{R}_{n-m+1}
 @>\theta_2>>\mc{M}_{n+1}
\end{CD}
\]
where the $\theta_i$ are induced from $\Theta$ and are therefore surjections,
while $\phi_1$ is the natural isomorphism. Thus, 
$\phi_2\cong \mu_n$ is a surjection. 

For the other direction, we may assume that 
$\mc{M} = \bigoplus_{n =a}^{\infty} \mc{M}_n$ where each $\mc{M}_n$
is coherent, and that $\mu_n$ is a surjection for all $n\geq b$.  
Then $\mc{F}=\bigoplus_{m=a}^b \mc{M}_m$ is a coherent $\mc{O}_X$-module
and there is a surjective map $\mc{F} \otimes \mc{R} \to \mc{M}$. 
\end{proof}

\begin{proposition}
\label{bimod alg noeth}
Keep the hypotheses of \eqref{global-convention}.
The bimodule algebra $\mc{R} = \mc{R}(X,c,\mc{L},\sigma)$ 
is right noetherian if and only if 
$\{c_i\}_{i \geq 0}$ is a critically dense subset 
of $X$, left noetherian if and only if 
$\{c_i\}_{i < 0}$ is critically dense in $X$, and noetherian 
if and only if $\{c_i\}_{i \in \mb{Z}}$ is critically dense in $X$. 

When $\mc{R}$ is not right noetherian, there exists an infinite 
ascending chain of coherent right ideals
of $\mc{R}$ with non-torsion factors.
\end{proposition}
\begin{proof}
By Lemma~\ref{L no matter}, the result is independent of the 
choice of $\mc{L}$ and we choose  $\mc{L} = \mc{O}_X$.
 
Assume first that $\{c_i\}_{i \geq 0}$ is critically dense in $X$.  We need to
show that  every right ideal of $\mc{R}$ is a coherent $\mc{R}$-module in the
sense of  Definition~\ref{graded def}.  By Lemma~\ref{structure J^n}, $\mc{R}_m
\otimes \mc{R}_n \cong \mc{R}_m \mc{R}_n$ for all $m, n \geq 0$,  and so we
may  use products of bimodules  in place of tensor products in 
the proof.   An arbitrary right ideal $\mc{G}$ of $\mc{R}$ is given by a
sequence of  bimodules $\mc{G}_i ={}_1(\mc{H}_i)_{\sigma^{i}}\subseteq
\mc{R}_i$ such that $\mc{G}_i \mc{R}_1 \subseteq \mc{G}_{i+1}$ for all $i
\geq 0$.   By Lemma~\ref{coherent}, $\mc{G}$ will  be coherent if and only if
$\mc{G}_{i} \mc{R}_1  = \mc{G}_{i+1}$ for $i \gg 0$. Equivalently, we are given
that
\begin{equation}\label{key-prop}
\mc{H}_{i} \mc{I}^{\sigma^{i}} \subseteq 
 \mc{H}_{i+1}\subseteq \mc{I}\mc{I}^{\sigma}\cdots\mc{I}^{\sigma^{i}}
\qquad \text{ for all } i \geq 0,
\end{equation}
 and $\mc{G}$ is  coherent if and only if
  $\mc{H}_{i} \mc{I}^{\sigma^{i}} =
\mc{H}_{i+1}$ for $i \gg 0$.

Assume that $\mc{G}\not=0$ and  pick $r$ such that
$\mc{H}_r\not=0$. By critical density, $\mc{H}_r$ can only be contained in
finitely many ideals $\mc{I}_j$ and so there exists $m\geq r$ such that
$\mc{H}_r\not\subseteq \mc{I}^{\sigma^j}=\mc{I}_{c_j}$, for all $j\geq m$.
Set $\mc{U}=\mc{H}_m$ and, for $n> m$, put  $\mc{W} =\mc{H}_n$ and
$\mc{V}=\mc{V}_n=\prod_{i=m}^{n-1} \mc{I}_{c_i}$.
By the choice of $m$ and  \eqref{key-prop},
 $\mc{U} \not\subseteq \mc{I}_{c_i}$ for $i\geq m$ and so $\mc{U}$ and $\mc{V}_n$
are comaximal. Thus \eqref{key-prop} and induction implies that
$\mc{U}\cap\mc{V}_n=\mc{U}\mc{V}_n\subseteq \mc{W}\subseteq \mc{V}_n$. Thus
Lemma~\ref{bimod-lemma} implies that $\mc{Z}_n=\mc{H}_m+\mc{H}_n$ is maximal
with respect to   $\mc{Z}_n\mc{V}_n \subseteq \mc{H}_n$. Since
$$\mc{Z}_n\mc{V}_{n+1}=\mc{Z}_n\mc{V}_n\mc{I}^{\sigma^n} \subseteq
\mc{H}_n\mc{I}^{\sigma^n}\subseteq \mc{H}_{n+1},$$ this implies that
 $\mc{Z}_n\subseteq \mc{Z}_{n+1}$  for all $n\geq m$.

Thus we may pick $n_0\geq m$ such that $\mc{Z}_n=\mc{Z}_{n+1}$ for
 all $n\geq n_0$. For all such $n$, Lemma~\ref{bimod-lemma}
 implies that 
$\mc{H}_{n+1}=\mc{Z}_n\mc{V}_{n+1}=\mc{Z}_n\mc{V}_n\mc{I}^{\sigma^n}
=\mc{H}_{n}\mc{I}^{\sigma^n}.$
Thus, $\mc{G}$ is coherent and $\mc{R}$ 
is right noetherian.

Conversely, suppose that $\{c_i\}_{i \geq 0}$ is not 
critically dense.  Then there exists an infinite set $A \subseteq \mb{N}$ such
 that 
the Zariski closure of the set of points $\{c_i\}_{i \in A}$ is a
reduced closed 
subscheme $Y \subsetneq X$ with defining ideal say
$\mc{I}_Y = \bigcap_{i \in A} \mc{I}^{\sigma^i}$.  
Set  $\mc{H}_n = \mc{I}_Y \cap \mc{I}_n$ 
for  $n \geq 0$, and observe that 
$\mc{G} = \bigoplus \mc{G}_n = \bigoplus {}_1(\mc{H}_n)_{\sigma^n}$ is 
a right ideal of $\mc{R}$.  
Write  $\mf{m}_{c_n}$ for
the maximal ideal in the local ring $\mc{O}_{X, c_n}$.  Looking 
locally at the point $c_n$, for $n\in A$,
 we have $(\mc{H}_n \mc{I}^{\sigma^n})_{c_n} = 
(\mc{I}_Y)_{c_n}\mf{m}_{c_n}$ while $(\mc{H}_{n+1})_{c_n} = (\mc{I}_Y)_{c_n}$.  
But $(\mc{I}_Y)_{c_n} \neq 0$ and so, by Nakayama's lemma,
$(\mc{I}_Y)_{c_n}\mf{m}_{c_n} \neq (\mc{I}_Y)_{c_n}$. Thus
$\mc{H}_n \mc{I}^{\sigma^n} 
\neq \mc{H}_{n+1}$ for any $n\in A$. Therefore, 
 the right ideal $\mc{G}$  is not coherent and  $\mc{R}$ 
fails to be right noetherian.  

 By \eqref{rees opposite}, 
 $\mc{R}^\op \cong \mc{R}(X, \sigma(c),\mc{L}^{\sigma^{-1}},\sigma^{-1})$
 and so  the  result on the left follows from the one on the right. The 
result for noetherian algebras is then obvious.

To prove the final statement, assume that $\mc{R}$ 
is not right noetherian, and let 
$\mc{G} = \bigoplus \mc{G}_n$ be the non-coherent right ideal
defined above.  
Set $\mc{M}^j = \sum_{0 \leq i \leq j} \mc{G}_i \mc{R}$; thus 
$(\mc{M}^j)_n = ((\mc{I}_Y\cap \mc{I}_j)\mc{I}^{\sigma^{j+1}}\cdots
 \mc{I}^{\sigma^{n-1}})_{\sigma^n}$ for $n\geq j$. This gives a chain of
coherent right ideals $\mc{M}^0 \subseteq \mc{M}^1 \subseteq \dots \subset
\mc{R}$.
Looking locally at a point $c_m\in A$ one finds that 
$\mc{M}^m_i \subsetneq \mc{M}^{m+1}_i$ for 
 all $m\in A$ and all $i\geq m$.  Thus the 
 subsequence $\{ \mc{M}^m : m \in A\}$ 
gives the desired 
ascending chain of right ideals of $\mc{R}$.
\end{proof}

\section{Ampleness}\label{ampleness-in-general}
We maintain the hypotheses from \eqref{global-convention}.  
The main aim of this section (Theorem~\ref{general-ampleness})
 is to prove, in considerable generality,
that the sequence of bimodules
$\{ \mc{R}_n = (\mc{I}_n\otimes \mc{L}_n)_{\sigma^n} \}$ 
is ample in the 
 sense of Definition~\ref{ample def}. 
Combined with the results of the last two sections this will 
prove parts~1 and 2 of Theorem~\ref{mainthm}.

\begin{theorem}
\label{general-ampleness}
Assume that $X$ is an integral projective scheme such that 
$X\not\cong\PP^1$ and fix $\sigma\in\Aut(X)$. 
Let $c\in X$ be a closed point for which
$\{c_i\}_{i\in{\mathbb Z}}$ is critically dense in $X$. 
Suppose that $\mc{L}$ is a very ample and $\sigma$-ample 
invertible sheaf on $X$.

Set  $\mc{R}=\mc{R}(X,c,\mc{L},\sigma)$. 
Then the section algebra $R = \Gamma(\mc{R})$ is  noetherian and
there is an equivalence of categories $\xi:\rqgr \mc{R} \simeq \rqgr R$ via the
inverse equivalences $\Gamma(-)$ and $-\otimes_R\mc{R}$. Similarly,
 $\mc{R}\qgr\simeq
R\qgr$.  
\end{theorem}

We prove this theorem through a series of reductions that give increasingly  
simple criteria for the ampleness of the sequence $\{ \mc{R}_n \}$.  Note that
only the left sheaf structure of the bimodules $\mc{R}_n$ really matters in 
the definition of right ampleness and 
 we write $\mc{J}_n=\mc{I}_n\otimes \mc{L}_n$ for that left sheaf. 
 By Lemma~\ref{structure J^n} we can and will identify $\mc{J}_n$ with
 $\mc{I}_n\mc{L}_n\subset \mc{L}_n$ for any $n\geq 1$.

\begin{lemma}
\label{reduction 2}  The following are equivalent:
\begin{enumerate}
\item[(1)]
The sequence $\{ \mc{R}_n \}$ is ample.
\item[(2)] The sheaf
$\mc{L}$ is $\sigma$-ample and  $\coH^1(\mc{M} \otimes \mc{R}_n) = 0$ 
for all invertible sheaves $\mc{M}$  and all $n \gg 0$.
\end{enumerate}\end{lemma}
\begin{proof}
Suppose that (1) holds and let $\mc{M}$ be an arbitrary invertible sheaf. 
Since $\mc{M}$ is  flat, the sequence 
\[
0 \to \mc{M} \otimes \mc{J}_n \to \mc{M} \otimes \mc{L}_n 
\to \mc{M} \otimes \mc{L}_n/\mc{J}_n \to 0
\]
is exact. Fix $i>0$. 
 Since $\mc{L}_n/\mc{J}_n$ is supported on a finite set of points, 
so is $\mc{M} \otimes \mc{L}_n/\mc{J}_n$. Therefore, 
$\coH^i(\mc{M} \otimes \mc{L}_n/\mc{J}_n) = 0$ while,   by hypothesis,
$\coH^i(\mc{M} \otimes \mc{J}_n) = 0$ for $ n \gg 0$. Thus
$\coH^i(\mc{M} \otimes \mc{L}_n) = 0 $ for $n \gg 0$ 
 and so \cite[Proposition~3.4]{AV} implies that  $\mc{L}$ is $\sigma$-ample.
  The claim regarding $\coH^1$ 
is a special case of Definition~\ref{ample def}(2).

Now suppose that (2) holds and 
let $\mc{F}$ be an arbitrary coherent sheaf. There is an exact sequence
\[
0 \to \calTor_1(\mc{F}, \mc{L}_n/\mc{J}_n) \to \mc{F} \otimes \mc{J}_n
\to \mc{F} \otimes \mc{L}_n \to \mc{F} \otimes \mc{L}_n/\mc{J}_n \to 0.
\]
By Lemma~\ref{prod versus tensor}(1),
$\calTor_1(\mc{F}, \mc{L}_n/\mc{J}_n)$ and $\mc{F} \otimes \mc{L}_n/\mc{J}_n$
are supported on finitely many points and hence their higher cohomology
vanishes. Thus $H^i(\mc{F} \otimes \mc{J}_n) \cong H^i(\mc{F} \otimes \mc{L}_n)$
for $i > 1$. Since $\{\mc{L}_n\}$ is an ample sequence, 
 these groups vanish for $n \gg 0$.

It remains to consider $H^1(\mc{F} \otimes \mc{J}_n)$.
There exists a short exact sequence $0 \to \mc{K} \to \mc{E} \to \mc{F} \to 0$
where $\mc{E}$ is a direct sum of invertible sheaves. Tensoring with $\mc{J}_n$
gives the exact sequence
\[
0 \to \calTor_1(\mc{F}, \mc{J}_n) \to \mc{K} \otimes \mc{J}_n
\to \mc{E} \otimes \mc{J}_n \to  \mc{F} \otimes \mc{J}_n \to 0.
\]
Once again, $\calTor_1(\mc{F}, \mc{J}_n)$ is finitely supported and so 
 $\coH^i(\calTor_1(\mc{F}, \mc{J}_n)) = 0$ 
for $ i > 0$.
By the last paragraph,
 $\coH^i(\mc{K} \otimes \mc{J}_n) = 0,$ for $ i > 1$ and $ n \gg 0$ while,
by hypothesis, $\coH^1(\mc{E} \otimes \mc{J}_n) = 0$ for $ n \gg 0$.
Thus $\coH^i(\mc{F} \otimes \mc{J}_n) = 0,$ for $ i > 0$ and $ n \gg 0$.
Hence Definition~\ref{ample def}(2) holds for $\mc{F}$.

Let $\mc{O}(1)$ be an arbitrary  very ample invertible sheaf on $X$.
Apply the conclusion of the last paragraph   to 
$\mc{F}\otimes \mc{O}(-i)$, for $1\leq i \leq \dim X$. This implies that, for
some $n_0\geq 0$, one has 
$${\coH}^i(X,\, \mc{F}\otimes\calO(-i)\otimes \calJ_{n})=0
\qquad \text{ for }  i>0\ \text{and}\  n \geq n_0.$$
By  \cite[Proposition~1, p.~307]{Kl} this implies that
$\mc{F} \otimes \calJ_n$ is generated by its sections for all $n\geq n_0$.  
Thus Definition~\ref{ample def}(1) holds for $\mc{F}$.
\end{proof}

Our second reduction gives a criterion for ampleness based on the following
notion  of separating points. Let $\mc{N}$ be an invertible sheaf on $X$ with
global sections $V = \coH^0(X,\mc{N})$.  We say $\mc{N}$
\emph{separates a set of closed points} $d_1, d_2, \dots, d_m\in X$ if, for
$1\leq i\leq m$, there exists $\alpha_i \in V$ such that $\alpha_i(d_j) =
0$ for $j \neq i$ but  $\alpha_i(d_i) \neq 0$.
By a slight abuse of notation we write this as $\alpha_i(d_j)=\delta_{ij}$
for  $1 \leq i,j \leq m$.

\begin{lemma}\label{separating}
Let $\calN$ be a very ample invertible sheaf on $X$ 
and set $1+h=\dim_k{\coH}^0(X,\calN)$. Let  
$\eta :X\hookrightarrow \mathbb P^h$
be the immersion corresponding to a basis of $\coH^0(X, \mc{N})$. Then:
\begin{enumerate}
\item 
$\mc{N}$ separates any pair of distinct points
of $X$.
\item $\mc{N}$ separates a triple $a,b,c$ of distinct closed points of
 $X$ if and only if the points 
$\eta(a),\eta(b),\eta(c)$ are not collinear in $\mathbb P^h$. 
\end{enumerate}
\end{lemma}

\begin{proof}
(1) This is just  \cite[Proposition~II.7.3]{Ha}.

(2) 
Each nonzero element $\alpha \in V = \coH^0(\mc{N})$ corresponds to a
hyperplane $Y \subseteq \mb{P}^h$ such that $\eta^{-1}(Y) = \{ x \in X \mid
\alpha(x) = 0 \}$. Now $\eta(a), \eta(b), \eta(c)$ are collinear in
$\mb{P}^h$ if and only if every hyperplane of $\mb{P}^h$ which contains
$\eta(a)$ and $\eta(b)$ also contains $\eta(c)$, if and only if  every $\alpha
\in V$ such that $\alpha(a) = \alpha(b) = 0$ also satisfies
 $\alpha(c) = 0$.  Since
the ordering of $a, b, c$ in this argument is immaterial this proves the
result. 
\end{proof}

\begin{lemma}
\label{reduction 3}
Assume that $\mc{L}$ is a $\sigma$-ample invertible sheaf on $X$. In order 
to prove that $\{ \mc{R}_n \}$ is an ample sequence, it is enough to show 
that, for any   $N \geq 0$, there is some $n_0 \geq 0$ such that 
the invertible sheaf $\mc{L}_n$ separates 
the points $c_{-N}, c_{-N+1}, \dots, c_{n-1}$ for all $n \geq n_0$.
\end{lemma}

\begin{proof} 
 Recall from Section~\ref{definitions} 
 that our convention is that $f^\sigma(x)=f(\sigma(x))$. 

 We first need to be careful about how  the points separated by 
$\mc{L}$ are related to the points separated by ${\mc L}^{\sigma^n}$.
So, suppose that an invertible sheaf $\mc{N}$ separates points $d_1,\dots, d_r$.
For some $1\leq i\leq d$, assume that
 $v\in \coH^0(\mc{N}) $ satisfies $v(d_j)=\delta_{ij}$ for $1\leq j\leq d$.
Then $w=v^{\sigma^m}\in \coH^0(\mc{N})^{\sigma^m} \cong 
\coH^0(\mc{N}^{\sigma^m})$ satisfies 
$w(\sigma^{-m}(d_j))=v(d_j)=\delta_{ij}$ and so $\mc{N}^{\sigma^m}$ separates 
the points $\sigma^{-m}(d_j)$. In particular, if 
$\mc{N}$ separates $c_{-m},\dots,c_r$, then $\mc{N}^{\sigma^m}$ separates 
$c_0,\dots,c_{m+r}$, since $c_i=\sigma^{-i}(c_0)$.
 Tensoring by a very ample sheaf $\mc{P}$
also preserves the property of separating a set of points. Indeed, 
if $v\in \coH^0(\mc{N})$ is defined as above, pick 
any   $u\in \coH^0(\mc{P})$ such that 
$u(d_i)\not= 0$.    Then 
$v\otimes u\in \coH^0(\mc{N}\otimes \mc{P})$
 certainly satisfies $(v\otimes u)(d_j)=\delta_{ij}$. Thus 
 $\mc{N}\otimes\mc{P}$ separates the points $d_1,\dots,d_r$.

We now turn to the proof of the lemma. 
 Fix an   invertible sheaf $\mc{M}$ on $X$.
Since $\mc{L}$ is $\sigma$-ample, 
there exists $m \geq 0$ such that $\mc{M} \otimes \mc{L}_m$ 
is very ample \cite[Proposition~2.3]{Ke1}.  Taking $N = m$ in the hypothesis 
shows that $\mc{L}_{n-m}$ separates $c_{-m}, c_{-m+1}, \dots, c_{n-m-1}$
for $n\gg 0$. By the first paragraph, 
$(\mc{L}_{n-m})^{\sigma^m}$   separates $c_0, c_1, \dots, c_{n-1}$.
Thus
$(\mc{M} \otimes \mc{L}_m) \otimes (\mc{L}_{n-m})^{\sigma^m} =
 \mc{M} \otimes \mc{L}_n$  also
separates $c_0, c_1, \dots, c_{n-1}$ for $n \gg 0$.

Write $\mc{N} = \mc{M} \otimes \mc{L}_n$ and 
consider the exact sequence 
\begin{equation}\label{reduction33}
0 \to  \mc{I}_n \otimes  \mc{N} \to \mc{N} \to 
\mc{N}/\mc{I}_n \mc{N} \to 0.
\end{equation}
Since $\mc{O}_X/\mc{I}_n \cong \bigoplus_{i =0}^{n-1} k(c_i)$ and  
  $\mc{N}$ is invertible, there are isomorphisms 
\[
\mc{N}/\mc{I}_n \mc{N}\  \cong \ \bigoplus_{i = 0}^{n-1} 
\mc{N}/\mc{I}^{\sigma^i}{\hskip-1pt}\mc{N}\  \cong \ 
\bigoplus_{i = 0}^{n-1} k(c_i).  
\]
Thus  $\alpha(c_j) = 0$ for some $\alpha \in \coH^0(\mc{N})$ 
if and only if $\psi_j(\alpha) = 0$, where $\psi_j$ is the map
\[
\psi_j: \coH^0(\mc{N})\ \overset{\theta}{\longrightarrow}\
  \coH^0(\mc{N}/\mc{I}_n \mc{N})\ = \
\bigoplus_{i = 0}^{n-1} \coH^0(k(c_i))\ \overset{\pi_j}{\longrightarrow}\
 \coH^0(k(c_j))
\ \cong\ k.
\]
The fact that
$\mc{N}$ separates the set of points $\{c_0, c_1, \dots, c_{n-1}\}$ ensures that 
the map  $\theta$ 
is surjective.  
But  $\mc{L}$ is $\sigma$-ample and so  $\coH^1(\mc{N}) = 0$
 for $n \gg 0$. From the long exact sequence 
$$\coH^0(\mc{N})\ {\buildrel{\theta}\over{\longrightarrow}}\
 \coH^0(\mc{N}/\mc{I}_n\mc{N})\longrightarrow
\coH^1(\mc{I}_n\otimes \mc{N}) \longrightarrow\coH^1(\mc{N}),$$ arising
 from \eqref{reduction33},
this implies that $\coH^1(\mc{I}_n \otimes \mc{N}) = 0$
for $n \gg 0$.  Thus $\coH^1(\mc{M} \otimes \mc{R}_n) = 0$ for $n \gg 0$ and 
the result follows from
 Lemma~\ref{reduction 2}.
\end{proof}

We can now  prove in considerable generality
that $\{ \mc{R}_n \}$ is ample.
 
\begin{proposition}
\label{not-general-ampleness} Assume that \eqref{global-convention}
holds.
Let $\mc{L}$ be a very ample and $\sigma$-ample sheaf on $X$
and assume that the Zariski closure of $\{c_i: i\geq 0\}$ is not isomorphic to
$\PP^1$.
Then the sequence $\{ \mc{R}_n \}$ is an ample sequence of bimodules.
\end{proposition}

\begin{proof}  Since $\sigma^m$ is an automorphism of $X$, the assumption on 
Zariski closures also ensures that the Zariski closure of $\{c_i: i\geq m\}$ 
is also  not isomorphic to $\PP^1$, for any $m\in{\mathbb Z}$. Consequently, 
for any closed immersion $\eta: X\hookrightarrow \PP^t$, the Zariski closure 
of $\{\eta(c_i) : i\geq m\}$ in $\PP^t$ cannot be contained 
in a line in $\PP^t$. It is this consequence of our hypothesis 
that will actually be used in the proof. In particular, 
if $\mc{M}$ is any very ample invertible sheaf, then  Lemma~\ref{separating}(2)
 implies that there exist three
points from the set $\{c_i\}_{i \geq m}$ which are separated by $\mc{M}$.

The strategy of the proof is to apply 
 Lemma~\ref{reduction 3}; thus for $N\geq 0$ we need 
 to show that, if $n$ is sufficiently large,
then for every $-N \leq k \leq n -1$
we can find $w\in \coH^0(\mc{L}_n)$
such that $w(c_i)=\delta_{ik}$  for  
$-N\leq i\leq n-1$.
The way we do this is to use the tensor product structure 
$\mc{L}_n=\mc{L}\otimes\cdots\otimes\mc{L}^{\sigma^{n-1}}$ and, for each
$0\leq i<n$, find some $w_i\in \coH^0(\mc{L}^{\sigma^i})$ such that 
$w_i(c_k)\not=0$, but $w_i(c_j)=0$ for ``sufficiently many'' $j\not=k$.
Then $w=w_0\otimes\cdots \otimes w_{n-1}$ will have the desired property.

Fix $N\geq 0$.  Pick three distinct points $d_{01}, d_{02}, d_{03}$ from
the set  $\{c_\ell  : \ell \geq -N \}$ that are separated by $\mc{L}$.  
Then pick 
$m_1$ such that $$\{d_{01}, d_{02}, d_{03} \} \subset \{c_{-N}, c_{-N+1}, \dots,
 c_{m_1 -1} \}.$$ Since  $\mc{L}^{\sigma}$ is  also very ample 
 we may choose  three distinct points $d_{11}, d_{12}, d_{13}$ from
$\{c_\ell : \ell \geq m_1\}$ that  are separated by $\mc{L}^{\sigma}$.
  Continuing inductively, 
   we may pick points $\{ d_{ji} : 0 \leq j \leq N-1\ \text{ and }\ 1 \leq i
\leq 3\}$ such that  the $d_{ji}$ are all distinct and such that
$\mc{L}^{\sigma^j}$ separates  $d_{j1}, d_{j2}, d_{j3}$ for each $j$.
Note that this implies that, for any $j$ and an arbitrary point $d\in X$, 
we can find  $u,v\in \{1,2,3\}$ such that
$\mc{L}^{\sigma^j}$ separates the points  $\{d,d_{ju},d_{jv}\} $.

Now choose $n$ large enough so that all of the $\{d_{ji}\}$ are contained in 
the set $\{c_\ell : -N \leq \ell \leq n  \}$ (where we have added an 
extra point for notational convenience).  Fix $-N \leq k \leq n $. 
 For each 
$0 \leq j \leq N-1$, the conclusion of the last paragraph implies that  
we can  find two of the $d_{ji}$, say $d_{j1}$ and $d_{j2}$, and some
$w_j \in {\coH}^0(\mc{L}^{\sigma^j})$ 
such that $w_j(c_k)\not=0$ but $w_j(d_{j1})=0=w_j(d_{j2})$.
Now let $\{e_i\}$ denote some enumeration of the 
remaining $n- N$ points; more precisely  write 
\[
\{c_\ell : -N \leq \ell \leq n \}=
\{ c_k\} \cup \{ d_{j1},d_{j2} : 0 \leq j \leq N-1\}
\cup \{e_1,\dots, e_{n- N}\}.
\]
Since each $\mc{L}^{\sigma^j}$ is very ample,
  Lemma~\ref{separating}(1) implies that, for 
$N \leq j \leq n -1$, there exists 
$w_j \in {\coH}^0(\mc{L}^{\sigma^j})$
such that $w_j(c_k)\not=0$ but
$w_j(e_{j-N+1})=0$.

Thus $w = w_0 \otimes w_1 \otimes \dots \otimes w_{n-1}$ is an element of 
$\coH^0(\mc{L} \otimes \mc{L}^{\sigma} \dots \otimes 
\mc{L}^{\sigma^{n-1}}) = \coH^0(\mc{L}_n)$
with the property that $w(c_i) =\delta_{ik}$ for  
$-N \leq i \leq n $.  Since $k$ is arbitrary this implies 
 that $\mc{L}_n$ separates
the set of points $c_{-N}, c_{-N+1}, \dots, c_{n }$. Thus   
  Lemma~\ref{reduction 3} can be applied to prove 
the ampleness of the sequence of bimodules $\{ \mc{R}_n \}$.
\end{proof}

\begin{proof}[Proof of Theorem~\ref{general-ampleness}] 
We may assume that $\dim X\geq 1$.  Suppose that 
$\{c_i : i\in {\mathbb Z}\}$ is critically dense.
By Proposition~\ref{bimod alg
noeth} $\mc{R}$ is right noetherian and by 
Proposition~\ref{abelian} $\rgr \mc{R}$ is an abelian category.
  Since $X \not \cong \mb{P}^1$, the hypotheses of 
Proposition~\ref{not-general-ampleness} are satisfied and that result
 implies that 
  $\{ \mc{R}_n \}$ is an ample sequence of bimodules.   
Thus all of the hypotheses of Theorem~\ref{VdB main theorem} are satisfied 
and so $R$ is right noetherian with $\rqgr \mc{R} \simeq \rqgr R$. 

As was noted in \eqref{rees opposite},
$\mc{R}^\op\cong \mc{R}(X,\sigma(c),\mc{L}^{\sigma^{-1}},\sigma^{-1})$ and,
by \cite[Corollary~5.1]{Ke1}, $\mc{L}^{\sigma^{-1}}$ is $\sigma^{-1}$-ample. 
Thus the claims on the left follow from 
those on the right. This completes the proof of 
  Theorem~\ref{general-ampleness}.
\end{proof}

\begin{remark}\label{remark-better} 
Proposition~\ref{not-general-ampleness} and hence 
 Theorem~\ref{general-ampleness} are  not the best results
possible. Using a  similar but more complicated argument we can prove that
$\{\mc{R}_n\}$ is an ample sequence as long as $\mc{L}$ is $\sigma$-ample, ample
and generated by global sections, and  the set $\{c_i\}$ is not all contained
in a dimension $1$ subscheme of $X$.  We do not know  what are
 the weakest possible
assumptions on $\mc{L}$ and the $\{c_i\}$ under which the proposition will hold.
 \end{remark}
 
The results we have proved this far also give
 a partial converse to Theorem~\ref{general-ampleness}. 

\begin{proposition}\label{not-noeth} Keep the assumptions of
\eqref{global-convention}.
Suppose that $\mc{C} = \{c_i\}_{i \geq 0}$ is not critically dense, 
but that the Zariski closure of $\mc{C}$ is not isomorphic to $\PP^1$.
Let $\mc{L}$ be very ample 
and $\sigma$-ample.
Then neither  $\mc{R}=\mc{R}(X,c,\mc{L},\sigma)$ nor 
$R = \Gamma(\mc{R})$ is right noetherian.
\end{proposition}
\begin{proof}
By Theorem~\ref{bimod alg noeth}, $\mc{R}$ is not right noetherian.  By 
Proposition~\ref{not-general-ampleness} $\{ \mc{R}_n \}$ is, however, a 
right-ample sequence of bimodules.

The final assertion  of 
Theorem~\ref{bimod alg noeth} provides an infinite proper 
ascending chain $\mc{M}^0 \subsetneq \mc{M}^1 \subsetneq \dots$ of coherent 
right $\mc{R}$-ideals such that the 
 factors $\mc{M}^{n+1}/\mc{M}^n$ are not in $\rtors \mc{R}$.
Since $\{ \mc{R}_n \}$ is an ample sequence, and  these are coherent ideals, 
the proof of \cite[Theorem~5.2, Step 1]{VB1} shows that, for each $n$,
$\mc{M}^n_i$ is generated by its global sections for  $i \gg 0$.
Writing $M^n = \Gamma(\mc{M}^n)$, this forces $M^n_i \subsetneq M^{n+1}_i$ 
for all $i \gg 0$,
and so $M^0 \subsetneq M^1 \subsetneq \dots$ is also
 a proper ascending chain of 
right $R$-ideals.  
Thus $R$ is not right noetherian. 
\end{proof}
 
 There is one minor case of  Theorem~\ref{general-ampleness} that  will not
be of interest in the sequel. This is when $X$ is a curve or a point.
The latter case is completely trivial, so assume that $X$ is a curve.
  In this case, since
$X$ has an infinite automorphism, it is either rational or elliptic. In
the former case it must also be singular.  In either case, since $c$ is a
smooth point (Lemma~\ref{smooth}), both  $\mc{I}=\mc{I}_c$ and 
$\mc{N}=\mc{I}\otimes\mc{L}$ are  invertible sheaves. Thus   $R$ is
nothing more than the twisted homogeneous coordinate ring 
$B(X,\mc{N},\sigma)$, as defined at the end of Section~\ref{definitions}.
 Note that, as $\mc{L}$ is very ample, the fact that 
 $X\not\cong {\mathbb P}^1$ implies that $\mc{L}$ must have at least $3$ global
 sections. An easy exercise 
then implies that $\mc{N}$ is ample. Hence
Theorem~\ref{general-ampleness} is just a very special case of 
\cite[Theorems~1.3 and 1.4]{AV} and $\mc{R}$ does not have any unusual
properties. In contrast,  the theorem does not hold for $X=\PP^1$ since one
can take $\mc{L}=\mc{I}^{-1}$.

The aim of the rest of the paper is to obtain a deeper understanding of 
the algebra $\mc{R}(X,c,\mc{L},\sigma)$ and its section ring $R=\Gamma(\mc{R})$
under the assumptions of Theorem~\ref{general-ampleness}. By
the last paragraph we are only interested in the case when
$\dim X\geq 2$. Thus for the rest of the paper we will  make the
following assumptions:

\begin{assumptions}\label{global-convention2} 
Let $X$ be a integral projective 
scheme of dimension $d\geq 2$.  
Fix $\sigma\in \Aut(X)$ and a very ample, 
$\sigma$-ample invertible sheaf $\mc{L}$. 
Finally assume that  $c\in X$ is chosen
so that $\mc{C} = \{ c_i :i \in \mb{Z}\} =
 \{ \sigma^{-i}(c) : i \in \mb{Z}\}$ 
is critically dense in $X$.  

We will always write $\mc{R} = \mc{R}(X,c,\mc{L},\sigma)$ and 
 $R=\Gamma(\mc{R}) = R(X,c,\mc{L},\sigma)$. 
By Theorem~\ref{general-ampleness}, $R$ is noetherian with 
 $\rqgr R\simeq \rqgr \mc{R}$.
\end{assumptions}

It is often useful to work with connected graded rings that are generated in
degree one and we end the section by giving two ways in which this may be
achieved; either by replacing $R$ by a large Veronese ring or by assuming that
the invertible sheaf $\mc{L}$ is ``sufficiently ample.''

\begin{proposition}
\label{veronese} 
Keep the hypotheses of \eqref{global-convention2}.
  Then  the Veronese ring 
$R^{(p)} = \bigoplus_{i \geq 0} R_{pi}$ is generated 
in degree $1$ for some $p\gg 0$.
\end{proposition}

\begin{remark} The ring $R^{(p)}$ is not quite an algebra of the same
form as $R$; explicitly $R^{(p)}$ is the global sections of the bimodule 
algebra  $ \bigoplus_{i \geq 0} ((\mc{I}_p \otimes
\mc{L}_p)_{\sigma^p})^{\otimes i}$ where $\mc{I}_p$ is the ideal sheaf defining
the finite set of points $\{c_0, c_1, \dots, c_{p-1}\}$. \end{remark}

\begin{proof} Set $\mc{J}_t=\mc{I}_t\otimes \mc{L}_t$ for $t\geq 1$.
By Proposition~\ref{not-general-ampleness} we may chose $r\geq 1$ such that 
$\mc{J}_r$ is generated by its global sections. 
Thus, there exists a short exact sequence
\[
0 \to \mc{V} \to \mathrm{H}^0(X,\mc{J}_r)\otimes \OO_X \to \mc{J}_r \to 0,
\]
of $\OO_X$-modules, for some sheaf $\mc{V}$.
Since  $\{\mc{J}_n\}$ and hence $\{\mc{J}^{\sigma^r}_n\}$ 
 is an ample sequence, there exists 
$n_0$ such that
$\mathrm{H}^1(X, \mc{V}\otimes \mc{J}^{\sigma^r}_n) = 0$ for $n\geq n_0$.
Hence, if one tensors the displayed exact sequence on the right with
$\mc{J}_{nr}^{\sigma^r}$ for $nr\geq n_0$ and takes global sections, 
one obtains the exact sequence
$$
\mathrm{H}^0(X,\mc{J}_r)\otimes \mathrm{H}^0(X,\mc{J}_{nr}^{\sigma^r})
\to \mathrm{H}^0(X,\mc{J}_{(n+1)r})
\to \mathrm{H}^1(X, \mc{V}\otimes \mc{J}_{nr}^{\sigma^r}).
$$
By construction, this final term is zero, so the exact sequence
is nothing more that the statement that the natural map $R_r\otimes R_{nr}
\to R_{(n+1)r}$ is a surjection.

By induction, for all $j\geq 1$ and all $n\geq  n_0$, we get
$R_{jr}R_{nr}=R_{(n+j)r}$. This implies that $R^{(nr)}$
is generated in degree one; that is by $R_{nr}$.
\end{proof}

In the next result, we write $\mc{L}_n^m = (\mc{L}_n)^{\otimes m}
\cong \mc{N}_n$ for $\mc{N}=\mc{L}^{\otimes m}$.

\begin{proposition}\label{degree 1}
Keep the hypotheses of 
\eqref{global-convention2}. 
Then there exists $M \in \NN$ such that, for $m \geq M$:
\begin{enumerate}
\item\label{deg11} $\mc{I}_n \otimes \mc{L}_n^m$ is 
generated by its global sections for all $n \geq 1$.
\item\label{deg12} $R(X,c,\mc{L}^m,\sigma)$ is generated in 
degree $1$. 
\end{enumerate} 
\end{proposition}

\begin{proof} (1) 
By Proposition~\ref{not-general-ampleness} 
there exists $n_0$ such that $\mc{I}_n \otimes \mc{L}_n$ is
generated by its global sections for $n \geq n_0$. Since $\mc{L}_n$ is already
globally generated, the sheaves $\mc{I}_n \otimes \mc{L}_n^m$
are globally generated for $n \geq n_0$ and $m \geq 1$. On the other hand,
 $\mc{L}_n$
is ample for all $n \geq 1$ and so there exists $m_0$ such that $\mc{I}_n \otimes
\mc{L}_n^m$ is globally generated for all $1\leq n\leq n_0$ and $ m \geq m_0$.
 Combining these   observations proves  (1).

(2) 
 The proof will use the following notion of Castelnuovo-Mumford regularity:
 An $\OO_X$-module $\mc{F}$ is \emph{$r$-regular with 
 respect to a very ample invertible sheaf $\mc{H}$} if
 $\mathrm{H}^i(X, \mc{F}\otimes \mc{H}^{(r-i)})=0$ for $i>0$.
 In this proof, all regularities will be taken with respect to 
 $\mc{H}= \mc{L}^\sigma$.
 The minimum $r$ such that $\mc{F}$ is $r$-regular is denoted $\reg\mc{F}$
 and called the \emph{regularity} of $\mc{F}$.
 Set $r = \max\{1, \reg \mc{O}_X \}$.
 
  For any $m\geq m_0$, part~1 provides a short exact sequence 
 \begin{equation}\label{degree 12}
0 \to \mc{K}_m \to \mathrm{H}^0(\mc{I} \otimes \mc{L}^m) \otimes \mc{O}_X \to 
\mc{I} \otimes \mc{L}^m \to 0,
\end{equation}
 for some sheaf $\mc{K}_m$. 
Since $\mc{L}$ is very ample, there exists $m_1\geq m_0$ such that
 $\mc{I} \otimes \mc{L}^m$ is $0$-regular with respect to
$\mc{L}^\sigma$ for all $m \geq m_1$. 
By \cite[Lemma~3.1]{AK} this implies that
 $\mc{K}_m$ is $r$-regular, independently of $m \geq m_1$.
 
 We   want to find similar upper bounds on the regularity of  other 
  sheaves.
 Since $\{ (\mc{I}_n \otimes \mc{L}_n)^\sigma \}$ is an ample sequence,
the regularity of $(\mc{I}_n \otimes \mc{L}_n)^\sigma$ is bounded above, 
independently of  $n \geq 1$.
Moreover,  by the vanishing theorem of \cite[Theorem~5.1]{Fj},
 there is a universal upper bound
on the regularity of \emph{any} ample invertible sheaf. Thus,
 $\reg((\mc{L}_{n-1}^{m-1})^{\sigma^2})$
is bounded above, independently of $n \geq 1$ and $ m \geq 2$.
Finally,  if $\mc{F}, \mc{G}$ are locally free except in a 
subscheme of dimension $\leq 2$,
then \cite[\S2]{Ke2} shows that 
$
\reg(\mc{F} \otimes \mc{G} ) \leq \reg \mc{F} + \reg \mc{G} +t$, where 
$t= (r-1)(\dim X -1).
$
Note that each $\mc{I}_n$, and hence each $\mc{K}_m$, is locally free except in
a subscheme of dimension $0$. 
Combining these observations shows that 
 $\reg(\mc{K}_m \otimes (\mc{I}_n \otimes \mc{L}_n)^\sigma
 \otimes (\mc{L}_{n-1}^{m-1})^{\sigma^2})$
is bounded above, independently of $n \geq 1$ and $m \geq m_1$.

Thus there exists 
$M\geq m_1$ such that, for all $m\geq m_1$,
$$\begin{array}{rl}
M & \geq  \
 \reg(\mc{K}_m \otimes (\mc{I}_n \otimes \mc{L}_n)^\sigma 
 \otimes (\mc{L}_{n-1}^{m-1})^{\sigma^2}) \\
 \noalign{\vskip 5pt}
&\qquad =\ \reg(\mc{K}_m \otimes (\mc{I}_n \otimes \mc{L}_n)^\sigma 
\otimes (\mc{L}_{n-1}^{m-1})^{\sigma^2} \otimes
(\mc{L}^\sigma)^{m-1}) +(m-1)
 \\
 \noalign{\vskip 5pt}
&\qquad\qquad= \
 \reg(\mc{K}_m \otimes (\mc{I}_n \otimes \mc{L}_n^m)^\sigma) +(m-1).
\end{array}$$
In other words, 
$\reg(\mc{K}_m \otimes (\mc{I}_n \otimes \mc{L}_n^m)^\sigma)  \leq 1$ 
for all $m\geq M$.
By \cite[Proposition~1, p.307]{Kl} this implies that 
$H^1(\mc{K}_{m} \otimes (\mc{I}_n \otimes \mc{L}_n^{ m})^\sigma) = 0$.

Now tensor   \eqref{degree 12}  
  with $(\mc{I}_n \otimes \mc{L}_n^{ m})^\sigma$ and note that 
the resulting sequence is exact, by Lemma~\ref{prod versus tensor}. 
Taking cohomology gives the  exact sequence
\begin{multline*}
\mathrm{H}^0(X, \mathcal{I}\otimes
\mathcal{L}^m)\otimes \mathrm{H}^0(X, (\mathcal{I}_n\otimes
\mathcal{L}^m_n)^\sigma)
\ {\buildrel{\theta}\over{\longrightarrow}}\
 \mathrm{H}^0(X, \mathcal{I}_{n+1}\otimes
\mathcal{L}^m_{n+1}) \longrightarrow \\
\noalign{\vskip 5pt}
 \longrightarrow  H^1(\mc{K}_{m} \otimes 
(\mc{I}_n \otimes \mc{L}_n^{ m})^\sigma).
\end{multline*}
By the conclusion of the last paragraph, this  final term is zero and
hence   $\theta$ 
is surjective for $n \geq 1$ and $ m \geq M$. This is equivalent to the assertion of
part~2.
\end{proof}

\section{Na{\"\i}ve noncommutative blowing up}\label{sect-blowup}

The hypotheses of Assumptions~\ref{global-convention2} will 
remain in force throughout this section.
In the introduction we asserted that one should regard the bimodule algebra
$\mc{R}$ as a sort of noncommutative blowup of $X$ at the point $c$. In this
section we justify and expand upon those comments, showing that they are easy
consequences of the basic construction. One should note, however, that
modules get slightly shifted and so 
it may be more natural to think of $\mc{R}$ as the blowup of $c_{-1}$, or 
perhaps better yet as the blowup of the 
entire orbit $\{c_n\}_{n \in \mb{Z}}$. We discuss this in more detail after 
Proposition~\ref{invertibility21}.
 
As was noted in the introduction, one way to form the blowup $\widetilde{X}$
of  $X$ at the closed point $c\in X$ is to use the identity $
\OO_{\widetilde{X}}\catmod = \rqgr \mc{A}$, 
where  $\mc{A}=\bigoplus {\mathcal I}_c^n$.
Since this bimodule algebra equals $ \mc{R}(X,c,\OO_X,\id)$, it is
natural to define $\rqgr \mc{R}(X,c,\OO_X,\sigma)$ to be   \emph{the na\"\i ve
noncommutative blowup} of $X$ at $c$. By Proposition~\ref{L no matter}, 
  $\rqgr \mc{R}(X,c,\OO_X,\sigma)\simeq
 \rqgr \mc{R}(X,c,\mc{L},\sigma)$ for any $\mc{L}$ and so 
 Theorem~\ref{general-ampleness} can be restated as:
 
 \begin{corollary}\label{blowup1}  Keep the hypotheses of
\eqref{global-convention2}.
 Then $\rqgr R(X,c,\mc{L},\sigma)$ is a na\"\i ve
 noncommutative blowup of $X$ at $c$. 
  \qed\end{corollary}

Perhaps the biggest difference between the classical and na{\"\i}ve
noncommutative blowups 
is in the properties of the inverse image of the smooth 
point that has been blown up.
In the commutative case, of course, one gets a divisor. However, in our case
we get just (a finite sum of copies of) one point, where we define
 a \emph{point} in $\rqgr \mc{R}$ to be a simple object in that category.

To explain this we need a further definition. 
As in \eqref{global-convention}, 
the skyscraper sheaf at a closed point $x\in X$ is written
$k(x)=\OO_X/\mc{I}_x$. Then one has a natural 
graded right $\mc{R}$-module 
 \begin{equation}\label{exceptional}
 \overline{x}= k(x)\oplus k(x)_{\sigma} \oplus k(x)_{\sigma^2} \cdots.
 \end{equation}
 The image of $\overline{x}$ in $\rqgr\mc{R}$ will   be written 
 $ \widetilde{x}$.
 It is an easy exercise to see that  $ \widetilde{x}$  is a simple object in 
 $\rqgr \mc{R}$ and we call it a \emph{closed point} in $\rqgr \mc{R}$.
 As will be seen in Theorem~\ref{GT equiv}, all points
  in $\rqgr \mc{R}$ are closed points, so the 
 notation is reasonable.

\begin{proposition}\label{blowup2}
 Keep the  hypotheses of
\eqref{global-convention2} and let $x\in X$ be a closed point. 
 Let $\widetilde{\rho}:  \OO_X\catmod\to \rqgr\mc{R}$  denote the blowup map
defined by $\rho: 
{\mathcal M}\mapsto {\mathcal M}\otimes_{{\mathcal O}_X}\mc{R}\in \rgr \mc{R}$. 
Then:
\begin{enumerate}
\item[(1)] If $x\not=c_j$ for $j\geq 0$, then $\widetilde{\rho}(k(x)) $ 
is the closed point $\widetilde{x}$ in $\rqgr\mc{R}$.
\item[(2)] If $x=c_j$ for $j\geq 0$, then $\widetilde{\rho}(k(x)) $ is a 
direct sum of $d=\dim X$ copies of 
the closed point $\widetilde{x}$  in $\rqgr \mc{R}$.
\end{enumerate}
\end{proposition}
\begin{remark} When $x=c_n$ for $n\in {\mathbb Z}$, we call $\widetilde{c}_n$ 
an \emph{exceptional point}.  In Theorem~\ref{GT equiv} we will 
slightly modify the functor $-\otimes \mc{R}$ in order to remove the 
direct sum that appears in part~2 of the proposition. This will show that  the
subcategory of torsion sheaves in  $ \OO_X\catMod$ is 
equivalent to a corresponding subcategory of $\rqgr \mc{R}$.
\end{remark}

\begin{proof} 
In the notation of Section~\ref{Rees bimod algs},
 $\rho(k(x))$ has the following structure as a left
$\OO_X$-module:
$$\rho(k(x)) = 
\bigoplus \left({\mathcal I}_n/{\mathcal I_x}{\mathcal I_n}\right)\otimes
 {\mathcal L}_n
\cong \bigoplus \left({\mathcal I}_n/{\mathcal I_x}{\mathcal I_n}\right).$$
If $x\notin\{c_j : j\geq 0\}$, then $\mc{I}_x$ and $\mc{I}_n$ are
comaximal for all $n \geq 0$ with
${\mathcal I}_n/{\mathcal I_x}{\mathcal I_n}
\cong {\mathcal I}_n/{\mathcal I_x}\cap{\mathcal I_n}
 \cong {\mathcal O}_X/{\mathcal I}_x\cong k(x).$
Thus $\widetilde{\rho}(k(x))=\widetilde{x}$.

On the other hand, if $x=c_j$ for some $j\geq 0$ then, for 
$n>j$, one has 
 ${\mathcal I}_{c_j}{\mathcal I}_n
={\mathcal I}_{c_j}^2{\mathcal I}^*_{n}$, where 
${\mathcal I}^*_n = \prod_{i=0}^{n-1}\{\mc{I}_{c_i} \mid i\not= j\}.$
By Lemma~\ref{smooth}, each $c_i$ is a smooth point of $X$ and so
 $${\mathcal I}_n/{\mathcal I}_{c_j}{\mathcal I}_n
\cong \left({\mathcal I}_{c_j}\cap{\mathcal I}^*_n\right) /\left(
{\mathcal I}_{c_j}^2\cap{\mathcal I}^*_n \right)
 \cong {\mathcal I}_{c_j}/{\mathcal I}_{c_j}^2\cong \bigoplus_{r=1}^d k(c_j).$$
Thus $\widetilde{\rho}(k(x))$ is 
   the direct sum of $d$ copies of the  
 point $\widetilde{c}_j$. 
 \end{proof}

In the commutative case one way to see that the exceptional divisor really is 
a divisor is as follows: Let $\rho:\widetilde{X}\to X$ be the blowup of  a
smooth variety $X$ at a point $c$ and  write $\widetilde{{\mathcal I}}_c =
\rho^{-1}({\mathcal I}_c)\cdot {\mathcal O}_{\widetilde{X}}$ for  the inverse
image ideal sheaf. Then \cite[Proposition~II.7.13]{Ha} shows that, under the
identification  $\OO_{\widetilde{X}}\catmod=\rqgr \mc{A}$ for
$\mc{A}=\bigoplus {\mathcal I}_c^n$, the sheaf $\widetilde{\mathcal I}_c$ is
simply the twisting sheaf ${\mathcal O}_{\widetilde{X}}(1)
 =\bigoplus_{n\geq 1}{\mathcal I}_c^{n+1}$. 
 As such it is invertible and so corresponds to a
(Cartier)  divisor.

A similar argument works in our situation with the exception of the last
sentence: there is no correspondence between invertible objects and divisors.
To explain this it is enough to work with $\mc{L}=\OO_X$ 
and $\mc{S}=\mc{R}(X,c,\OO_X,\sigma)$
and so for most of this section we will work with that bimodule algebra. 
As we will see, the exceptional point $\widetilde{c}_{-1}$ of
Proposition~\ref{blowup2} is indeed equal to $\mc{S}/\mc{K}$
where $\mc{K}$ is naturally isomorphic to 
the shift $\mc{S}[1]$. Similar results hold for each $\widetilde{c}_n$.

It is natural to write the shifts $\mc{S}[m]$ in the form 
$\bigoplus \mc{F}_n \otimes \mc{L}_{\sigma}^n$, for some $\OO_X$-modules 
$\mc{F}_n$ but, as the next lemma shows, one has to be careful about 
 the powers of $\sigma$ appearing in the $\mc{F}_n$.  
 
\begin{lemma}
\label{shifts}  Suppose that $\mc{N}$ is a right module 
over $\mc{R}=\mc{R}(X,c,\mc{L},\sigma)$, for $\mc{L}$ arbitrary,
that can be written in the form 
 $\mc{N} = \bigoplus \mc{F}_n \otimes  \mc{L}_\sigma^{\otimes n}$ for some 
 $\OO_X$-modules $\mc{F}_n$ with the trivial bimodule
structure. If $m \in \mathbb Z$, then
$\mc{N}[m] \cong \bigoplus \mc{G}_n \otimes  \mc{L}_\sigma^{\otimes n}$ where:
 
$\mc{G}_n = (\mc{F}_{n+m}\otimes \mc{L}_m)^{\sigma^{-m}}$
(with the trivial bimodule structure) if $m>0$ and

$\mc{G}_n = (\mc{F}_{n+m})^{\sigma^{-m}} \otimes\mc{L}_{-m}^{-1}$ if $m<0$.
\end{lemma}

\begin{proof} The result is well-known for the bimodule algebra
$\mc{B} = \mc{B}(X,\mc{L},\sigma)$, as defined at the 
end of Section~\ref{definitions}. In particular, when $m \geq 0$ the
result for $\mc{B}$ follows immediately from \cite[(3.1)]{SV} and 
the same argument works for $\mc{R}$. A simple computation then 
gives the required  formula for $m\leq 0$. 
\end{proof}

Using this lemma we find that, at least in $\rqgr\mc{S}$,
 the shift $\mc{S}[n]$ is isomorphic to the following right 
 $\mc{S}$-module $\widetilde{\mathcal K}(n) $.
 
\begin{definition}\label{invert-defn1} Let $\mc{S}=\mc{R}(X,c,\OO,\sigma)
\cong \bigoplus_{r\geq 0} \mc{I}\mc{I}^\sigma\cdots\mc{I}^{\sigma^{r-1}}
 \otimes \OO_{\sigma}^{\otimes r} $ for $\OO=\OO_X$.
For  $n\in \mathbb Z$,  define $\OO_X$-bimodules

$$
\widetilde{\mathcal K}(n) =
\bigoplus_{r\geq 1+|n|} {\mathcal I}^{\sigma^{-n}}
{\mathcal I}^{\sigma^{-n+1}}
\cdots{\mathcal I}^{\sigma^{r-1}} \otimes \mc{O}_{\sigma}^{\otimes r}
\cong 
\bigoplus_{r\geq 1+|n|} \left({\mathcal I}^{\sigma^{-n}}
\cdots{\mathcal I}^{\sigma^{r-1}}\right)_{\sigma^r}$$
 and 

$$
\widetilde{\mathcal K}(n)^* = 
\bigoplus_{r\geq 1+|n|}  {\mathcal I}{\mathcal I}^{\sigma}
\cdots{\mathcal I}^{\sigma^{r-n-1}} \otimes \mc{O}_{\sigma}^{\otimes r} 
\cong
\bigoplus_{r\geq 1+|n|}  \left({\mathcal I}
\cdots{\mathcal I}^{\sigma^{r-n-1}}\right)_{\sigma^r}$$
 \end{definition}

The next result describes the bimodule structure of 
 $\widetilde{\mathcal K}(n)$ and requires the following definitions.
  A module $\mc{M}\in \rqgr \mc{B}$ is called  \emph{an invertible 
$(\mc{A},\mc{B})$-bimodule in qgr} if it is the image under the quotient
functor $\pi_{\mc{B}}$  of an $(\mc{A},\mc{B})$-bimodule $\mc{M}$ and there
exists a $(\mc{B},\mc{A})$-bimodule $\mc{M}'$ such that, up to torsion,
$\mc{M}\otimes_{\mc{B}}\mc{M}' \cong\mc{A}$ and $\mc{M}'\otimes_{\mc{A}}\mc{M}
\cong \mc{B}$. The module $\mc{M}'$ is called the \emph{inverse} of $\mc{M}$.

\begin{proposition}\label{invertibility1}
 Set $\mc{S} = \mc{R}(X, c,\OO,\sigma) $ and $\mc{S}'=\mc{R}(X,
c_{-n},\OO,\sigma).$ Then  $\widetilde{\mathcal K}(n)$ is 
 an invertible $\left(\mc{S}',\mc{S}\right)$-bimodule in qgr with inverse 
  $\widetilde{\mathcal K}(n)^*$.
 \end{proposition}

 \begin{proof} By Lemma~\ref{structure J^n},
  $\mc{S}_1 \otimes \widetilde{\mc{K}}(n)^*_r =
  \mc{I}\otimes {}_1\mc{O}_\sigma\otimes \widetilde{\mc{K}}(n)^*_r 
 =\widetilde{\mc{K}}(n)^*_{r+1}$, for any $r$ such that 
  $\widetilde{\mc{K}}(n)^*_r\not=0$. Thus
   $ \widetilde{\mc{K}}(n)^*$ is a left $\mc{S}$-module and, 
 similarly, $ \widetilde{\mc{K}}(n)$ is a right $\mc{S}$-module.
 If
 $\mc{H}=\mc{I}^{\sigma^{-n}}$,  then $\mc{S}'$ is  the bimodule algebra
 $\bigoplus_{r\geq 0}\mc{H}\mc{H}^\sigma\cdots
 \mc{H}^{\sigma^{r-1}}\otimes \mc{O}_{\sigma}^{\otimes r}.$ From 
 its definition one finds that 
$$ \widetilde{\mathcal K}(n) =
\bigoplus_{r\geq 1+|n|} \mc{H}\mc{H}^\sigma \cdots
 \mc{H}^{\sigma^{n+r-1}}\otimes \mc{O}_\sigma^{\otimes r},$$  which
is obviously a left $\mc{S}'$-module. A similar argument shows that 
$ \widetilde{\mathcal K}(n)^*$ is a right $\mc{S}'$-module and the 
left and right  actions are clearly
compatible.

Now consider $\widetilde{\mathcal K}(n)^*\otimes_{\mc{S}'}
\widetilde{\mathcal K}(n)$. For $t$ sufficiently large, Lemma~\ref{tensor
bimodules} implies that its
  $t^{\mathrm{th}}$ summand 
is 
$$
\sum_i 
\left(\mc{I}\cdots\mc{I}^{\sigma^{i-n-1}}
\right)_{\sigma^i}
\left(\mc{I}^{\sigma^{-n}}\cdots\mc{I}^{\sigma^{t-i-1}} 
\right)_{\sigma^{t-i}} 
= \sum_i\left(\mc{I}\cdots\mc{I}^{\sigma^{t-1}}
\right)_{\sigma^t}
 = \mc{S}_t,
$$
as required. 
  The proof that 
  $\widetilde{\mathcal K}(n)\otimes_{\mc{S}}
\widetilde{\mathcal K}(n)^* \cong \mc{S}'$ in $\rqgr\mc{S}'$
   is essentially the same.
  \end{proof}
  
The point in constructing the
 $\widetilde{\mc{K}}(n)$ was to show that 
   the exceptional point $\widetilde{c}_{-1}$ could be written as 
 $\mc{S}/\widetilde{\mc{K}}(1)$.
In fact we have a more general result: 
 
  \begin{proposition}\label{invertibility21}
 In $\rqgr\mc{S}$ one has 
 $\widetilde{c}_{-n-1} \cong \widetilde{\mc{K}}(n)/\widetilde{\mc{K}}(n+1)$
for all $n\in \mathbb Z$.
  \end{proposition}
  
\begin{proof} Since the ${\mathcal I}^{\sigma^i}$ are comaximal,
   ${\mathcal I}^{\sigma^{-n}}
\cdots{\mathcal I}^{\sigma^{-1}}\big/{\mathcal I}^{\sigma^{-n-1}}
\cdots{\mathcal I}^{\sigma^{-1}} \cong 
{\mathcal O}_X/{\mathcal I}^{\sigma^{-n-1}}$,
for all $n$.
 Thus in $\rqgr\mc{S}$ we have 
\begin{align*}
\qquad\qquad\widetilde{\mc{K}}(n)/\widetilde{\mc{K}}(n+1)
\ \cong\ \bigoplus_{r\geq 2 + |n|}
  \left(\OO_X/\mc{I}^{\sigma^{-n-1}}\right)_{\sigma^r} 
  \ =\ \widetilde{c}_{-n-1}.\qquad\qquad{}\qedhere\end{align*}
\end{proof}

We remarked in the introduction to this section that there is some ambiguity in
what is actually been blown up in the passage from $\rcatmod\OO_X$ to $\rqgr
\mc{S}$. Given the way the construction works, we feel the correct
interpretation is that we have just na{\"\i}vely blown up the point $c$ but, in
the process, we automatically blow up the full orbit $\{c_i : i\in \mathbb
Z\}$.  This is illustrated by Proposition~\ref{invertibility21}:  all the
modules  $\widetilde{c}_n\in \rqgr \mc{S} $ can be written as a factor  of two
invertible bimodules and so, in this respect, they  are more similar to a
divisor than to a point. Another way of viewing the same result (which can also
be proved directly from Lemma~\ref{shifts}) is that $\widetilde{c}_n =
\widetilde{c}[-n]$. 

\begin{remark} In the discussion above we have concentrated on  modules over
$\mc{S}$ since this most naturally corresponds to the commutative 
description of blowing up.  In fact, one can also define the commutative
blow-up in terms of $\bigoplus \mc{I}^n\mc{L}^{\otimes n}$ (see
\cite[Lemma~II.7.9]{Ha}) and so one should expect that these results for 
$\mc{S}$ have natural analogues for   $\mc{R}=\mc{R}(X,c,\mc{L},\sigma)$.
This is true. We leave it to the reader to check that  
Definition~\ref{invert-defn1}--Proposition~\ref{invertibility21} all remain true
as results about $\mc{R}$-modules if one simply replaces  $\OO$ by $\mc{L}$
and $\mc{S}$ by $\mc{R}$ in
their statements and adjusts the proofs accordingly. Note that the new module
$\widetilde{\mc{K}}(n)$ will not be isomorphic to $\mc{R}[n]$ in $\rqgr \mc{R}$.
  \end{remark}

There is another version of noncommutative blowing up, due to Van den Bergh
\cite{VB2}, that does have the    properties expected of blowing up  and his
theory  has  been very useful in describing noncommutative surfaces. While Van
den Bergh's construction is notationally   very similar to ours, it is actually
rather different. For one thing, it only works in the following situation: one
has a  connected graded ring $S$ of Gelfand-Kirillov dimension $3$ that
surjects onto a twisted homogeneous coordinate ring  $B=B(E,{\mathcal N},
\tau)$, where $E$ is a curve, and the aim is to blow up a point $c$ on that 
curve. This requires the construction of something  analogous to the graded
algebra $\bigoplus \mc{I}_c^n$. However, since  $B$ is noncommutative,
$\bigoplus \mc{I}_c^n$ does not carry a natural algebra structure and so one
has to be more subtle.   Rather than use a category like $\rqgr S$, Van den
Bergh works in the category of left exact functors from $\rqgr S$ to itself and
so, in particular, $\pi(S)$ is replaced by the identity functor on $\rqgr S$.
It is then nontrivial to show that the blowup has the expected properties. 
In particular, the inverse image of $c$  does look like a divisor. The
details can be found in \cite{VB2}  and a brief introduction to the
construction and its applications are described in \cite[Section~13]{SV}.

\section{$\mc{R}$-modules and equivalences of categories}
\label{R-modules}

The hypotheses from Assumptions~\ref{global-convention2} will remain in force
throughout this section. One nice consequence of critical density is that
it forces modules over $\mc{R}=\mc{R}(X,c,\mc{L},\sigma)$ and
$R=\Gamma(\mc{R})$ to have a very pleasant structure; indeed in many cases they
are just induced from $\OO_X$-modules. 
This will be used in this section to give various equivalences of categories,
notably that the category of  coherent 
torsion $\OO_X$-modules is equivalent to the subcategory 
of Goldie torsion modules in $\rqgr \mc{R}$, as defined below. 
This gives the promised improvement of Proposition~\ref{blowup2}. 
We also give a natural analogue of the standard fact that,
for a blowup 
$\rho:\widetilde{X}\to X$ at a smooth point $x$, the schemes
$X\smallsetminus \{x\}$ and $\widetilde{X}\smallsetminus \rho^{-1}(x)$ 
are isomorphic. For these  results it would be sufficient, by
Proposition~\ref{L no matter}, to work with just $\mc{R}(X,c,\OO_X,\sigma)$,
but we will work in the general case since this will enable us to draw
conclusions about $R$-modules.

If $A$ is a noetherian graded domain, a graded
$A$-module $M$ is called  \emph{Goldie torsion} (to distinguish this from the
notion of torsion already defined)  if every homogeneous element of $M$ is
killed by some nonzero homogeneous  element of $A$. Equivalently, 
$M$ is a direct limit of modules of the form $(A/I)[n]$, for nonzero 
graded right ideals $I$. The latter notion passes to all the categories $Q$ 
we consider; for example a right $\mc{R}$-module is Goldie torsion if 
it is a direct limit of modules of the 
form $(\mc{R}/\mc{K})[n]$ for nonzero right ideals $\mc{K}$ of $\mc{R}$. 
Of course, Goldie torsion $\OO_X$-modules are just the torsion $\OO_X$-modules,
as in \cite[Exercise~II.6.12]{Ha}.
We write $\GT Q$ for the full  subcategory of
  Goldie torsion modules
 in~$Q$. 

We start by giving some technical results on the structure of Goldie torsion
modules.
If $\mc{N}=\bigoplus \mc{N}_n\in \rGr\mc{R}$, 
recall  from \eqref{nice-modules} that we may write 
each $\mc{N}_n$ as 
an $\mc{O}_X$-bimodule of  the form $(\mc{G}_n)_{\sigma^n}$. 
It is often convenient to write   $(\mc{G}_n)_{\sigma^n} =
\mc{F}_n \otimes \mc{L}_{\sigma}^{\otimes n}$,
  where $\mc{F}_n = {}_1(\mc{F}_n)_1$ 
 has trivial bimodule structure and, as usual, 
 $\mc{L}_\sigma^{\otimes n}=(_1\mc{L}_\sigma)^{\otimes n}$.
    This has the advantage that 
 the module structure of
$\mc{N}$ is now given by  maps of (left) sheaves 
 $\mc{F}_n\otimes \mc{I}^{\sigma^n} \to \mc{F}_{n+1}$ for all $n$.

\begin{lemma}
\label{Goldie tors facts}
{\rm{(1)}} If 
$\mc{N} = \bigoplus \mc{F}_n \otimes \mc{L}_{\sigma}^{\otimes n}\in \GT\rgr
\mc{R}$, then    there exists a single module 
$\mc{F}\in \GT \OO_X\catmod$ such that
$\mc{F}_n = \mc{F}$ for all
 $n \gg 0$.  

{\rm{(2)}}
Conversely,
if $\mc{F}\in \GT \OO_X\catmod$, then 
$\bigoplus_{n=0}^{\infty} \mc{F}\otimes 
 \mc{L}_{\sigma}^{\otimes n}\in \GT\rgr \mc{R}$.

\end{lemma}
\begin{proof} (1) Clearly
 $\mc{F}_n$ is  Goldie torsion for $n \gg 0$. 
Since $\mc{N}$ is coherent, by Lemma~\ref{coherent} there is a surjection 
$\mc{F}_n \otimes \mc{I}^{\sigma^n} \twoheadrightarrow \mc{F}_{n+1}$ for 
$n \gg 0$ and so the  supports satisfy 
$\supp \mc{F}_{n+1} \subseteq \supp \mc{F}_n$ for all $n\gg 0$. 
Since $\mc{C} = \{c_i \}_{i \geq 0}$ is critically dense, 
$\mc{C} \cap\, \supp \mc{F}_n$ is finite for each $n$ such that $\mc{F}_n$
is torsion.
Thus, $c_n \not \in \supp \mc{F}_n$ for $n \gg 0$.
By Lemma~\ref{prod versus tensor}(2) and the fact that
$\mc{F}_n/\mc{F}_n\mc{I}^{\sigma^n}$ is supported on $(\supp \mc{F}_n)\cap
\{c_{n}\}=\emptyset$, one has $\mc{F}_n\otimes \mc{I}^{\sigma^n}\cong
\mc{F}_n  \mc{I}^{\sigma^n}=\mc{F}_n$ for all $n\gg 0$.
 Thus, we obtain a surjection 
 $\mc{F}_n \twoheadrightarrow \mc{F}_{n+1}$
 for all $n\geq n_0$.
 Since the $\mc{F}_n$ are noetherian this forces
 $\mc{F}_n \cong \mc{F}_{n+1}$ for $n \gg n_0$.

(2) Since $Y = \supp \mc{F}$ is a proper closed subset of $X$ and the points
$\{ c_i \}_{i \geq 0}$ are critically dense,  
$Y \cap \{c_i\}_{i \geq n} = \emptyset$ for $n \gg 0$.  So, just as in part~1, 
$\mc{F}\otimes  \mc{I}^{\sigma^n}\cong 
\mc{F} \mc{I}^{\sigma^n} \cong \mc{F}$ for $n \gg 0$. Thus, if 
$\mc{M} = \bigoplus_{n=0}^\infty \mc{F}\otimes  \mc{L}_{\sigma}^{\otimes n}$,
then there is a surjection
$\mc{M}_n \otimes \mc{R} \twoheadrightarrow \mc{M}_{\geq n}$ for
$n \gg 0$ and $\mc{M}$ is coherent by Lemma~\ref{coherent}.  
By construction, $\mc{M}$ is a Goldie torsion $\mc{R}$-module. 
\end{proof}

\begin{lemma}\label{goldie-subfactors}
If $\mc{M}\in\rgr \mc{R}$, then $\mc{M}$ has a finite filtration of
submodules $0=\mc{M}^0\subset\mc{M}^1\subset \cdots\subset \mc{M}^r=\mc{M}$
such that the factors $\mc{M}^i/\mc{M}^{i-1}$ 
are equal to either a shift $\mc{R}[j]$ of $\mc{R}$ or 
to a Goldie torsion module.
\end{lemma}
\begin{proof} 
This is similar to the analogous result for finitely generated 
modules over rings and
the proof is left to the reader.
\end{proof}

The point of the next lemma is that, at least  in large 
degree, the structure of an $\mc{R}$-module may be written using
 products 
instead of tensor products.
\begin{lemma}
\label{product versus tensor}
Let $\mc{N}  = 
\bigoplus \mc{F}_n\otimes \mc{L}_{\sigma}^{\otimes n}
\in \rgr \mc{R}$.  Then, under the natural map, 
 $$\mc{F}_n\otimes \mc{I}^{\sigma^n}\cong  \mc{F}_n\mc{I}^{\sigma^n}
=\mc{F}_{n+1}\qquad \text{for }\ n\gg
0.$$
\end{lemma}

\begin{proof}
Let $0 \to \mc{H}' \to \mc{H} \to \mc{H}'' \to 0$ be an exact sequence
 of sheaves on $X$, and consider the commutative diagram
\[
\xymatrix{ \mc{H}' \otimes \mc{I}^{\sigma^n} \ar[r] \ar[d]^{\theta_1}
	& \mc{H} \otimes \mc{I}^{\sigma^n} \ar[r] \ar[d]^{\theta_2}
	& \mc{H}'' \otimes \mc{I}^{\sigma^n} \ar[r] \ar[d]^{\theta_3} & 0 \\
		\mc{H}'\mc{I}^{\sigma^n} \ar[r]^{\phi} 
	& \mc{H}\mc{I}^{\sigma^n} \ar[r]
	& \mc{H}''\mc{I}^{\sigma^n} }
\]

\noindent
where the $\theta_i$ are the natural surjections. 
Although the bottom row is not in general exact, the map $\phi$
is injective. A diagram chase then shows that, if $\theta_1$ 
and $\theta_3$ are isomorphisms, then $\theta_2$ is an isomorphism. 

If  $\mc{N}$ is a Goldie
torsion module, then the proof 
of Lemma~\ref{Goldie tors facts} shows that
 $\mc{F}_n \otimes \mc{I}^{\sigma^n} 
\cong \mc{F}_n \mc{I}^{\sigma^n}$ for $n \gg 0$. 
Alternatively, suppose that  $\mc{N}=\mc{R}[m]$
 is a shift of $\mc{R}$.  Lemma~\ref{shifts} implies that 
 $\mc{F}_n = \mc{I}^{\sigma^{-m}}\cdots\mc{I}^{\sigma^{n-1}}
 \otimes \mc{L}_{|m|}^\alpha$
 for  $n>|m|$ and the appropriate  $\alpha$.  
  Lemma~\ref{prod versus tensor}(2) therefore implies 
   that the map $\mc{F}_n \otimes \mc{I}^{\sigma^n} \to \mc{F}_n
\mc{I}^{\sigma^n}$ must  be  an isomorphism for $n > |m|$.

If $\mc{N}$ is an arbitrary coherent $\mc{R}$-module, 
then by Lemma~\ref{goldie-subfactors}  it has a finite
filtration by shifts of $\mc{R}$ and Goldie torsion modules. It therefore
follows  by induction from the last two paragraphs that $\mc{F}_n \otimes
\mc{I}^{\sigma^n} \to \mc{F}_n \mc{I}^{\sigma^n}$  is an isomorphism for $n \gg
0$. Since $\mc{N}$ is coherent, the induced maps  $\phi_n: \mc{F}_n \otimes
\mc{I}^{\sigma^n}=  \mc{F}_n \mc{I}^{\sigma^n} \to \mc{F}_{n+1}$  defining  the
module structure of $\mc{N}$  are  surjections for all $n \geq n_0$. 

It remains to prove that the surjections 
$\phi_n: \mc{F}_n \mc{I}^{\sigma^n} \to \mc{F}_{n+1}$ are isomorphisms
 for $n\gg n_0$. Pulling back to $\mc{F}_{n_0} $, we may write 
$\mc{F}_{n} = \mc{A}_n/\mc{B}_n$, for submodules 
$\mc{B}_n\subseteq \mc{A}_n\subseteq
\mc{F}_{n_0}$. Since $\mc{F}_{n+1}$ is a homomorphic image of 
$\mc{F}_n\mc{I}^{\sigma^n} = (\mc{A}_n  \mc{I}^{\sigma^n} + \mc{B}_n)/\mc{B}_n$, we find that 
$\mc{B}_{n+1} \supseteq \mc{B}_n$ for each $n\geq n_0$.
 Since $\mc{F}_{n_0}$ is noetherian, $\mc{B}_n=\mc{B}_{n+1} $ for all 
 $n\gg n_0$, and hence 
$ \mc{F}_n \mc{I}^{\sigma^n} \cong \mc{F}_{n+1}$ for all such $n$.
\end{proof}

As might be expected, the fact that $\mc{R}$-modules have a nice form 
is also reflected in the homomorphism groups. The next lemma collects the
relevant facts.  Recall from \eqref{quotient} that  the natural map from
 $\rGr\mc{R}$ to $\rQgr\mc{R}$ is denoted $\pi$.

\begin{lemma}
\label{Hom in Qgr}
Let  $\mc{N} = \bigoplus \mc{G}_n \otimes 
\mc{L}_{\sigma}^{\otimes n}\in \rGr \calR$.  
Then:

{\rm (1)} There is a natural isomorphism 
$$
\Hom_{\rQgr \mc{R}}(\pi(\mc{R}), \pi(\mc{N})) \cong \lim_{n \to \infty} 
\Hom_{\mc{O}_X}(\mc{I}_n, \mc{G}_n),$$
  where the limit is induced from the 
multiplication map; specifically it 
sends $\theta\in \Hom_{\mc{O}_X}(\mc{I}_n, \mc{G}_n)$ to the map
$$\theta\otimes \mc{I}^{\sigma^n} : \mc{I}_{n+1}\cong \mc{I}_n\otimes
\mc{I}^{\sigma^n}\to \mc{G}_n\otimes \mc{I}^{\sigma^n}\to \mc{G}_{n+1}.$$

{\rm (2)}
Suppose that $\mc{M} = \bigoplus \mc{F}_n \otimes\mc{L}_{\sigma}^{\otimes n}$ 
and $\mc{N}$ are coherent Goldie torsion 
and, by Lemma~\ref{Goldie tors facts}, 
write $\mc{F}_n = \mc{F}$  and $\mc{G}_n = \mc{G}$ for $n \geq n_0$.  
Then there is a natural isomorphism
$$
\Hom_{\rQgr \mc{R}}(\pi(\mc{M}), \pi(\mc{N})) 
\cong \Hom_{\mc{O}_X}(\mc{F}, \mc{G}).
$$
\end{lemma}

\begin{proof}
  The definition of homomorphisms in quotient categories implies that 
\begin{equation}\label{hom-in-quotient}
\Hom_{\rQgr \mc{R}}(\pi(\mc{M}),\pi(\mc{N})) = \lim_{n\to \infty}
\Hom_{\rGr \mc{R}}(\mc{M}_{\geq n},\mc{N}),
\end{equation} 
whenever $\mc{M}$ is coherent (see, for example,
\cite[p.31]{VB2}).
 On the other hand, we claim that there are natural vector space maps
\begin{equation}\label{hom-in-quotient2}
\hom_{\rGr \calR}(\mc{M}_{\geq n}, \mc{N}) 
\overset{\phi_n}{\lra} \hom_{\OO_X}
(\mc{M}_n, \mc{N}_n) \overset{\rho_n}{\lra} \hom_{\mc{O}_X}(\mc{F}_n, \mc{G}_n).
\end{equation} 
Indeed, if $f \in \hom_{\rGr \calR}(\mc{M}_{\geq n}, \mc{N})$, then $f$
is a morphism of right $\mc{O}_X$-modules, so 
we may define $\phi_n(f)$ to be the restriction of $f$ to $\mc{M}_n$.
The map $\rho_n$ is the natural isomorphism obtained by tensoring with
the invertible bimodule $(\mc{L}_{\sigma}^{\otimes n})^{-1}.$ 

(1)  The result will  follow   from 
 \eqref{hom-in-quotient} and  \eqref{hom-in-quotient2}
 once we prove that the map $\rho_n\circ\phi_n$ is an isomorphism for $n\gg 0$.
 In this case, $\mc{F}_n=\mc{I}_n $ and, by Lemma~\ref{structure J^n},
 $\mc{I}_{n+r}\cong \mc{I}_n\otimes \mc{I}_r^{\sigma^{n}}$, for any $n,r\geq 0$.
 Thus, if $g\in \hom_{\mc{O}_X}(\mc{I}_n, \mc{G}_n)$,
 then $g$ induces a unique map 
 $$\mc{R}_{n+r} \cong
  \mc{I}_n\otimes (\mc{I}_r^{\sigma^{n}}\otimes\mc{L}_\sigma^{\otimes (n+r)})
 \to \mc{G}_n\otimes (\mc{I}_r^{\sigma^{n}}\otimes\mc{L}_\sigma^{\otimes (n+r)})
 \to \mc{G}_{n+r}  \otimes\mc{L}_\sigma^{\otimes (n+r)}=\mc{N}_{n+r},$$
 for any $r\geq 0$. This clearly defines an 
 $\mc{R}$-module map $f\in \hom_{\OO_X}
(\mc{R}_{\geq n}, \mc{N})$ such that 
 $\rho_n\phi_n(f)=g$. Thus  $\rho_n\phi_n$ is in fact an isomorphism  for 
 all $n \geq 0$.

(2) The proof is essentially the same as that of part~1. 
Any element 
$g \in \hom_{\mc{O}_X}(\mc{F} , \mc{G})$
determines a unique map
 $\mc{F} \otimes \mc{L}_\sigma^n \to \mc{G}\otimes\mc{L}_\sigma^n $
 for   $n\geq n_0$. For such an $n$
 this  defines  an $\mc{R}$-module map 
 $f: \mc{M}_{\geq n}\to \mc{N}_{\geq n}$ with $\rho_n\phi_n(f)=g$.
Thus $\rho_n\phi_n$ is an isomorphism for $n\geq n_0$
and the result  follows   from 
 \eqref{hom-in-quotient} and  \eqref{hom-in-quotient2}.
 \end{proof}

It is now easy to define an 
 equivalence of categories between $\GT \rqgr R$ and $\GT \OO_X\catmod$.

\begin{theorem}
\label{GT equiv} Keep the hypotheses from \eqref{global-convention2}.
Then there are equivalences of categories 
$$\GT \rQgr R\ \simeq\ \GT \rQgr \mc{R} \ \simeq\ \GT \OO_X\catMod,$$
which restrict to equivalences
$\GT \rqgr R\simeq \GT \rqgr \mc{R} \simeq \GT \OO_X\catmod.$
This equivalence is given by mapping $\mc{F}\in \GT \OO_X\catMod$
to $\pi\left(\bigoplus \mc{F}\otimes \mc{L}^{\otimes n}_\sigma\right)
\in \rQgr \mc{R}$.
\end{theorem}

\begin{remark}\label{GT equiv1} For any closed point $x\in X$, 
this equivalence sends $k(x) \in \OO_X\catmod$ to 
$\widetilde{x}\in\GT \rqgr \mc{R}$ and so 
it does give the promised refinement of Proposition~\ref{blowup2}.
Since the simple objects in $\OO_X\catMod$ are precisely
these modules $k(x)$, this 
 also proves part~4 of Theorem~\ref{mainthm}.
\end{remark}

\begin{proof}  The equivalence of categories Theorem~\ref{general-ampleness}
clearly restricts to an equivalence $\GT \rQgr R\simeq \GT \rQgr \mc{R}$, so
only the second equivalence needs proving.

Define $\theta: \GT \mc{O}_X\catmod \to \GT \rqgr \mc{R}$  by $\mc{F} \mapsto
\pi(\bigoplus \mc{F}\otimes \mc{L}_{\sigma}^{\otimes n})$.   That $\theta$ 
lands in $\GT \rqgr \mc{R}$  rather than $\GT \rQgr \mc{R}$  follows from
Lemma~\ref{Goldie tors facts}(2), and $\theta$ is  clearly functorial
since $-\otimes (\bigoplus \mc{L}_{\sigma}^{\otimes n})$ and $\pi$ are functors.
Lemma~\ref{Goldie tors facts}(1) shows that  $\theta$ is surjective on
objects.  Finally, it follows  from Lemma~\ref{Hom in Qgr}(2) that $\theta$ is
full and faithful on morphisms,  so that $\theta$ is an equivalence.  

By \cite[Theorem~1.1.1]{SV}, 
  $\GT \rQgr\mc{R}$  is   the closure of 
   $\GT \rqgr\mc{R}$ under direct limits  and similarly for 
$\GT\OO_X\catMod$.   Thus, the equivalence $\GT \mc{O}_X\catmod
\simeq \GT \rqgr \mc{R}$ extends  to an equivalence $\GT
\mc{O}_X\Mod \simeq \GT \rQgr \mc{R}$. Since direct limits commute with tensor
products, the equivalence does still have the specified form.
 \end{proof}

A standard fact in geometry is that if $\rho: \widetilde{X}\to X$ 
is the blowup of $X$ at a smooth point $x$, then $X\smallsetminus\{x\}$ 
is isomorphic to $\widetilde{X}\smallsetminus \rho^{-1}(x)$
\cite[Proposition~II.7.13]{Ha}. The final result of this section 
proves the analogous result for $\rqgr\mc{R}$. 
As may be expected from the results of Section~\ref{sect-blowup}
we have to remove all of the $c_i$ from $X$ rather than just one point. 
Thus we define $C_X$ to be the smallest localizing subcategory of $\mc{O}_X\catMod$ 
containing all of the $\{k(c_i)\}_{i \in \mb{Z}}$.  
The reader may check that the 
objects in this subcategory are exactly those
quasicoherent sheaves which are supported at the set of 
points $\{c_i\}_{i \in \mb{Z}}$ and so  $C_X$ is a
subcategory of $\GT \mc{O}_X\catMod$.  

Similarly,  write 
$C_\mc{R}$ for the localizing subcategory  of $\rQgr\mc{R}$ 
generated by the modules $\widetilde{c}_n$ for $n\in \mathbb Z$.
By Remark~\ref{GT equiv1}, the equivalence $\GT \mc{O}_X\catMod\simeq 
\GT \rQgr \mc{R}$ restricts to an equivalence  
$C_X \simeq C_{\mc{R}}$.

\begin{proposition}
\label{factor-equiv} Keep the hypotheses from \eqref{global-convention2}.
There is an equivalence of  categories 
$\OO_X \catMod/C_X  \simeq \rQgr \mc{R}/C_{\mc{R}}.$
\end{proposition}

\begin{proof} 
Given a sheaf $\mc{F} \in \OO_X \catMod$, write $\overline{\mc{F}}$ 
for the corresponding 
object in $\OO_X \catMod/C_X$.  Set $D_X = C_X \cap \OO_X \catmod$ and 
$D_{\mc{R}} = C_{\mc{R}} \cap \rqgr\mc{R}$. 

As in Proposition~\ref{blowup2}, define a map
$\theta' : \OO_X\catmod \to \rqgr\mc{R}$ by 
 $ \mc{F} \mapsto \pi(\mc{F}\otimes \mc{R})$. 
   By that proposition, $\theta'(k(c_i))$ is a direct sum of copies of 
   $\wt{c}_i$,    and so $\theta'$
induces  a functor $\theta: \OO_X\catmod/D_X  \to \rqgr \mc{R}/D_{\mc{R}}$. 
Conversely, let
$\mc{M} = \bigoplus \mc{F}_n \otimes \mc{L}_{\sigma}^{\otimes n} 
\in \rgr \mc{R}$. Then Lemma~\ref{product versus tensor} implies
 that, for some $\omega$,
 \begin{equation}\label{factor-equ}
 \mc{F}_n\otimes  \mc{I}^{\sigma^n} \cong \mc{F}_n \mc{I}^{\sigma^n} \cong 
 \mc{F}_{n+1}\qquad \text{for}\ n\geq \omega.
 \end{equation} 
  Thus,  
$\overline{\mc{F}}_n \cong \overline{\mc{F}}_{n+1}$ for all $n \geq  \omega$
and so the 
rule $\psi':\mc{M}\mapsto \overline{\mc{F}}_{\omega}$ defines a functor from  
$\rqgr \mc{R}$ to  $\OO_X\catmod/D_X$.  
For all $i$, the object $\wt{c}_i$ maps to $0$ and so there is an  
induced functor $\psi: \rqgr \mc{R}/D_{\mc{R}} \to \OO_X \catmod/D_X$.  

We need to check that these are inverse functors. First, since 
$\mc{I}_n$ and $\OO_X$ are isomorphic modulo $D_X $, it is clear that 
$\psi'\theta'(\mc{F}) = \overline{\mc{F}} $ for $\mc{F}\in \OO_X\catmod$.
On the other hand, 
if $\mc{M} 
\in \rgr \mc{R}$ and   $\omega$  satisfy \eqref{factor-equ},
then  $\theta\psi'(\mc{M})$ is the image of the module 
$\mc{M}' =\mc{F}_{\omega}\otimes \mc{R} 
\in \rgr \mc{R}$. Therefore, 
$\mc{M}'_{\geq \omega} = \mc{F}_{\omega}\otimes \mc{I}_{\omega}\otimes
\mc{L}_\sigma^{\otimes \omega}\otimes \mc{R}$
whereas $\mc{M}_{\geq \omega} = \mc{F}_{\omega}\otimes 
\mc{L}_\sigma^{\otimes \omega}\otimes \mc{R}$, by the choice of $\omega$.
Thus, there is a natural map 
$\alpha: \mc{M}'_{\geq \omega}\to \mc{M}_{\geq \omega}$ whose 
cokernel  is a homomorphic image of 
$\mc{F}_{\omega}\otimes \OO_X/\mc{I}_{\omega}\otimes 
\mc{L}_\sigma^{\otimes \omega}\otimes \mc{R}$, which is clearly contained in
$\theta'(D_{X})\subseteq D_{\mc{R}}$.
 Similarly, Lemma~\ref{prod versus tensor} implies that
  $\mathrm{Ker}(\alpha)  \subseteq D_{\mc{R}}$ and so 
  $\mc{M}=\mc{M}'$ in $\rgr\mc{R}/D_{\mc{R}}$. 
 
Therefore, both $\theta\psi$ and $\psi\theta$ are naturally
isomorphic to the identity and we have an equivalence
$\OO_X \catmod/D_X  \simeq \rqgr
\mc{R}/D_{\mc{R}}$.  It now follows formally, as in the proof of
Theorem~\ref{GT equiv}, that this induces the desired equivalence $\OO_X
\catMod/C_X  \simeq \rQgr \mc{R}/C_{\mc{R}}$. \end{proof}

\section{The chi conditions}
\label{chi} 

The hypotheses from Assumptions~\ref{global-convention2} remain in force
throughout the current section.  The aim is to prove that
$R=R(X,c,\mc{L},\sigma)$ satisfies the $\chi_1$ condition but not the $\chi_2$
condition, as defined in the introduction. In the process we will prove that 
${\coH}^1(\pi(R))=\ext^1_{\rQgr R}(\pi(R),\pi(R))$
  is infinite dimensional, thereby proving parts~7 and 8 of
Theorem~\ref{mainthm}.  
  An easy consequence (see Corollary~\ref{artin-zhang3})
will be that $\rqgr R$ satisfies the ampleness condition of Artin and Zhang
 \cite{AZ1}.
  For a  detailed discussion of these various conditions the reader is referred
to  \cite{AZ1,SV}.

\begin{theorem}
\label{chi_1, not chi_2} Keep the hypotheses from
\eqref{global-convention2}. Then, on both the left and the  right:
\begin{enumerate}
\item 
$\chi_1$ holds for $R$.
\item 
$\chi_2$ fails for $R$. Indeed $\ext_{\rcatMod R}^2(k,R)$
 is infinite dimensional.
\item
 ${\coH}^1(\pi(R))=\ext^1_{\rQgr R}(\pi(R),\pi(R))$
  is infinite dimensional.
  \end{enumerate}
\end{theorem}

Before beginning the proof of the theorem we need two easy technical lemmas.
The proofs will use the derived functors $\calExt^i$ of $\calHom$, the 
Sheaf Hom on $\mc{O}_X$-modules in the sense of \cite[Section~III.6]{Ha}.

\begin{lemma}
\label{exts from points}
Let $\mc{F} \in \mc{O}_X\catmod$ and fix $i \in \mb{Z}$.  Then:  
\begin{enumerate}
\item If $c_i \not \in \supp \mc{F}$, then $\Hom_{\OO_X}(k(c_i), \mc{F}) = 0 = 
\ext_{\OO_X}^1(k(c_i), \mc{F})$. 
\item Suppose that  $\mc{F} \subseteq \mc{M}$, where $\mc{M}$ is locally free
with  $c_i \not \in \supp \mc{M}/\mc{F}$. Then $\ext_{\OO_X}^1(k(c_i),
 \mc{F}) = 0$.
 \item If $\mc{F}\in \GT \OO_X\catmod$ then 
 $\ext^j_{\OO_X}(\mc{I}_n,\mc{F})=0$ for all $j\geq \dim X$ and all $n \geq 0$.
\end{enumerate}
\end{lemma}

\begin{proof} (1,2)
Let $\mc{O}(1)$ be a very ample invertible sheaf on $X$. Then 
$k(c_i) \otimes \mc{O}(n) \cong k(c_i)$ for all $n\in \mathbb Z$ and so,
by \cite[Propositions~III.6.7 and  III.6.9]{Ha}, it suffices to prove that 
$\calExt^j(k(c_i),\mc{F})= 0$ (where $j=0,1$ in part~1 and $j=1$ in part~2).
By \cite[Proposition~III.6.8]{Ha} the question is now a local one;
for any $x\in X$ one has
$\calExt^i(k(c_i),\mc{F})_x = \ext^i_{\OO_{X,x}}(k(c_i)_x,\,\mc{F}_x).$

Part~1 now follows from the fact that, for any closed point $x\in X$, either 
$k(c_i)$ or $\mc{F}$ is zero at $x$. Hence so are the $\calExt$ groups.

 In order to prove part~2, note that, by part~1,
 $\calHom(k(c_i), \mc{M}/\mc{F}) = 0$.  
So by the long exact sequence in $\calExt$ it is enough to show that 
$\calExt^1(k(c_i), \mc{M}) = 0$.  
But $\dim X\geq 2$ and, by Lemma~\ref{smooth},
$X$ is smooth and hence Cohen-Macaulay at $c_i$.  
Thus $\calExt^1(k(c_i), \mc{M})_{c_i} = 0$. At any other  
closed point $c_i\not=y\in X$, one has $k(c_i)_y=0$ and hence
 $\calExt^1(k(c_i), \mc{M})_{y} = 0$. 

(3) Pick $j\geq d=\dim X$ and consider the exact sequence 
$$\ext^j(\OO_X, \mc{F}) \to \ext^j(\mc{I}_n, \mc{F}) \to 
\ext^{j+1}(\bigoplus_{i=0}^{n-1}k(c_i),\mc{F}).$$
Since $\mc{F}$ is torsion, $\mc{F}$
 is supported on a subscheme 
of dimension $< d$ and so  $\ext^j(\OO_X, \mc{F})=0,$ by 
\cite[Lemma~III.2.10]{Ha}. It therefore  suffices
to show that $\ext^{m}(k(c_i),\mc{F}) = 0$, for all $m>d$.
 It again suffices to prove this
 locally, and since $k(c_i)$ is supported at $c_i$
  we only need to look locally at  $c_i$.  
Since $X$ is smooth at
  $c_i$,  $\OO_{X,c_i}$ has global dimension 
  $\leq d$ and hence $\calExt^m(k(c_i),\mc{F})_{c_i}= 0$. 
\end{proof}

\begin{lemma} \label{h1-thm2} 
As an $\calO_X$-module,
$\calHom_{\calO_X}(\calI_n, \mc{O}_X/\calI_n)
\cong \bigoplus_{i=0}^{n-1} k(c_i)^{d}$, where $d = \dim X$.  Thus $\dim_k 
\Hom_{\calO_X}(\calI_n, \mc{O}_X/\calI_n) = nd$.
\end{lemma}
\begin{proof} 
Looking locally, one calculates that 
$\calHom_{\calO_X}(\calI_n, \mc{O}_X/\calI_n) \cong \calI_n/\calI_n^2$.
But $\mc{O}_X/\mc{I}_n \cong \bigoplus_{i=0}^{n-1} k(c_i)$.
Since $X$ is smooth locally at each point $c_i$,
the first statement follows.  The second assertion follows because
$\Hom_{\calO_X}(\calI_n, \mc{O}_X/\calI_n)$ is isomorphic to
 the global sections of 
$\calHom_{\calO_X}(\calI_n, \mc{O}_X/\calI_n)$.
\end{proof}

\begin{proof}[Proof of Theorem~\ref{chi_1, not chi_2}]
We prove this result on the right only; the left-sided results then
 follow by appealing to \eqref{rees opposite}.  

(1) We begin by making some reductions to the problem. 
One characterization of the $\chi_1$ 
condition for $R$ is that, for all $N \in \rgr R$, the natural map 
\begin{equation}\label{chi-first}
N \to \bigoplus_m \Hom_{\rQgr R}(\pi(R), \pi(N)[m]) 
\end{equation}
should have right bounded cokernel \cite[Proposition~3.14]{AZ1}.  
This condition clearly holds for $N$ if and only if it holds for a shift
$N[r]$. Since 
 $N$ has a filtration by shifts of $R$ and Goldie torsion
modules, it suffices to prove the condition in those two
 cases.
 
 We convert \eqref{chi-first} into a statement about the $\mc{R}$-module
 $\mc{N}=(N\otimes_R\mc{R})$. By the equivalences of categories,
  Theorem~\ref{general-ampleness}, 
   $\xi^{-1}\circ\pi_R(N)=\pi_{\mc{R}}(\mc{N})$
  and $N=\Gamma(\mc{N})$ in high degree.   Thus 
we may rephrase  $\chi_1$
as requiring that  the natural map 
\begin{equation}
\label{max tors ext} \Gamma({\mc{N}}) \to \bigoplus_m \Hom_{\rQgr 
\mc{R}}(\pi(\mc{R}), \pi(\mc{N})[m])
\end{equation}
has a right bounded cokernel. 
Since $R\otimes_R\mc{R}=\mc{R}$ and
 $N\otimes_R\mc{R}$ is Goldie torsion for any Goldie torsion module $N$,
 it will also suffice
  to show that  \eqref{max tors ext} has right bounded cokernel 
when $\mc{N}$ is either $\mc{R}$ or a Goldie torsion module.

Write 
 $\mc{N} = 
\bigoplus  \mc{F}_n\otimes  \mc{L}_{\sigma}^{\otimes n}$ and 
$\mc{R} = \bigoplus  \mc{I}_n \otimes \mc{L}_{\sigma}^{\otimes n}$. Fix $m\gg 0$ and
write $\mc{N}[m] = \bigoplus \mc{G}_n \otimes 
\mc{L}_{\sigma}^{\otimes n}$; thus
$\mc{G}_n \cong (\mc{F}_{n+m} \otimes \mc{L}_m)^{\sigma^{-m}}$
 by Lemma~\ref{shifts}.
By  Lemma~\ref{Hom in Qgr}(1) we have   an isomorphism
$$
\Hom_{\rQgr \mc{R}}(\pi(\mc{R}), \pi(\mc{N})[m]) \cong \lim_{n \to \infty} 
\Hom_{\mc{O}_X}(\mc{I}_n, \mc{G}_n).
$$
Write  \begin{equation}
\label{nth map} \psi_n: \hom_{\mc{O}_X}(\mc{I}_n, \mc{G}_n) \to 
\hom_{\mc{O}_X}(\mc{I}_{n+1}, \mc{G}_{n+1})
\end{equation}
for the $n^{\mathrm{th}}$ map in this direct limit and observe that
the  zeroth term  is nothing more  than 
$\hom_{\mc{O}_X}(\mc{I}_{0}, \mc{G}_{0})=\coH^0(X,\mc{G}_0)
\cong\coH^0(X,\mc{F}_m\otimes \mc{L}_m)=
 \Gamma(\mc{N})_m$. 
Thus if we can show, for $m \gg 0$, that the map $\psi_n$ is
an isomorphism for all $n \geq 0$, then we will have shown 
that \eqref{max tors ext} has right bounded cokernel and proved the result.

As  $m$  is sufficiently large, Lemma~\ref{product versus tensor}
implies that 
$\mc{G}_{n+1} = \mc{G}_n \mc{I}^{\sigma^n}$ for all $n \geq 0$.
In particular,  $\mc{I}_{n+1} \subseteq \mc{I}_n$, 
$\mc{G}_{n+1} \subseteq \mc{G}_n$, and so Lemma~\ref{Hom in Qgr}(1) implies
that $\psi_n$ is just the restriction map
for all $n \geq 0$.
It is clear then that $\psi_n$ fits into the following commutative
diagram, where the rows and columns are parts of the long exact sequences 
in $\Hom$ induced from the aforementioned inclusions:
{\small \[
\xymatrix@C=2.2em{ \Hom(k(c_n),\mc{G}_{n+1}) \ar[r] \ar[d] &
\Hom(\mc{I}_n,\mc{G}_{n+1}) \ar[r]^{\theta_{n+1}} \ar[d] &
\Hom(\mc{I}_{n+1},\mc{G}_{n+1}) \ar[r] \ar[d]^{\rho_n} &
\ext^1(k(c_n),\mc{G}_{n+1}) \ar[d] \\
	\Hom(k(c_n),\mc{G}_n) \ar[r] &
\Hom(\mc{I}_n,\mc{G}_n) \ar[r]_{\theta_n} \ar[ru]_{\psi_n} &
	\Hom(\mc{I}_{n+1},\mc{G}_n) \ar[r] &
	\ext^1(k(c_n),\mc{G}_n) } 
\]}

Assume first that $\mc{N} = \mc{R}$ and fix $n \geq 0$.  Then 
$\mc{G}_n = (\mc{I}_{n+m} \otimes \mc{L}_m)^{\sigma^{-m}}$, and so 
$\supp (\mc{L}_m^{\sigma^{-m}}/\mc{G}_n)=\{c_{-m},\dots,c_{n-1}\}
\not\ni c_n$. By Lemma~\ref{exts from points}(2) this implies that 
  $\ext^1(k(c_n), \mc{G}_n) = 0$.  Also $\Hom(k(c_n), \mc{G}_n) = 0$ since 
$\mc{G}_n$  is a subsheaf of an invertible sheaf.  Thus $\theta_n$ is 
an isomorphism.  Since $\rho_n$ is injective, the maps 
$\psi_n$ and $\rho_n$ are also isomorphisms for all $n \geq 0$.   
For later reference, note that this step  works for any $m \geq 0$.

Now assume that $\mc{N}$ is Goldie torsion.  For $n \geq n_0\gg 0$,
Lemma~\ref{Goldie tors facts} implies that $\mc{F}_n = \mc{F}_{n+1}=\mc{F}$,
 say. We may assume that $m\geq n_0$ and hence that 
  $\mc{G}_n=\mc{G}_{n+1}=\mc{G}
=(\mc{F}\otimes\mc{L}_m)^{\sigma^{-m}}$ for all $n\geq 0$. 
 Since $\mc{F}$ is a Goldie torsion sheaf, its support
is a proper closed subset of $X$  and so  
 $\supp \mc{F} \cap \{c_i\}_{i \in \mb{Z}}$
is finite, say contained in $\{c_j : |j|\leq r\}$. Thus  
$\supp \mc{G} \cap \{c_i\}_{i \in \mb{Z}}\subseteq \{c_j : j\leq r-m\}$. As 
$m$ is sufficiently large we may also  assume that $m>r$ and hence 
that $\supp \mc{G} \cap \{ c_i : i \geq 0\}= \emptyset$.
  Lemma~\ref{exts from points}(1) therefore implies that 
$\Hom(k(c_n), \mc{G}) = 0 = \ext^1(k(c_n), \mc{G})$ for all 
$n \geq 0$. The commutative diagram is now just the statement that
$\psi_n = \theta_n=\theta_{n+1}$ is an 
isomorphism for all $n \geq 0$. 

Thus, in either case $\psi_n$ is an isomorphism for $n \geq 0$ 
and so \eqref{max tors ext} 
does have the required bounded cokernel.

(2) In order to prove that $\chi_2$ fails on the right for the $R$-module 
$R$, it   suffices
to show  that 	$\dim_k\ext^1_{\rQgr R}(\pi(R),\pi(R))=\infty$
 \cite[Theorem~7.4]{AZ1}.  Thus (2) follows from (3).

(3)  
By Theorem~\ref{general-ampleness}, again,
 this is equivalent to showing  that 
$$
\dim_k\ext^1_{\rQgr \calR}(\pi(\calR) , \pi(\calR))=\infty.
$$
Write $\mc{B} = \bigoplus_{n \geq 0} \mc{L}_{\sigma}^{\otimes n}$, 
which we think of as  
a right $\mc{R}$-module.
The long exact sequence in $\Hom$ induced from the inclusion $\pi(\mc{R})\subset
\pi(\mc{B})$ provides the following exact sequence:
\begin{multline}\label{chitwo}
\qquad\Hom_{\rQgr \mc{R}}(\pi(\mc{R}), \pi(\mc{B}))
 \overset{\phi_1}{\longrightarrow}
\Hom_{\rQgr \mc{R}}(\pi(\mc{R}), \pi(\mc{B}/\mc{R})) \\
\overset{\phi_2}{\longrightarrow} 
\ext^1_{\rQgr \mc{R}}(\pi(\mc{R}), \pi(\mc{R})). \qquad
\end{multline}
We need to understand the first two terms in the sequence.

 Lemma~\ref{Hom in Qgr}(1) implies that 
 $\Hom_{\rQgr \mc{R}}(\pi(\mc{R}), \pi(\mc{B}))
\cong {\displaystyle \lim_{n\to \infty}}
 \hom_{\mc{O}_X}(\mc{I}_n, \mc{O}_X)$. 
By Lemma~\ref{exts from points}(2),
$\ext^1(k(c_i), \mc{O}_X) = 0$  and so 
 the natural map $\Hom(\mc{O}_X, \mc{O}_X) \to \Hom(\mc{I}_n, \mc{O}_X)$ is an 
isomorphism for all $n \geq 0$. Thus
$\Hom_{\rQgr \mc{R}}(\pi(\mc{R}), \pi(\mc{B}))
\cong 
 \hom_{\mc{O}_X}(\mc{O}_X, \mc{O}_X) = k$.

On the other hand,  $\Hom_{\rQgr \mc{R}}(\pi(\mc{R}), \pi(\mc{B}/\mc{R}))\cong 
{\displaystyle \lim_{n\to \infty}}
 \hom_{\mc{O}_X}(\mc{I}_n, \mc{O}_X/\mc{I}_n)$. 
As
$\Hom_{\mc{O}_X}(k(c_n), \mc{O}_X/\mc{I}_n) = 0$, it follows from the same 
commutative diagram as 
in part~1 that all of the maps $\psi_n$ in this direct limit are injective.
However, by Lemma~\ref{h1-thm2},
${\displaystyle \lim_{n\to \infty}}
\dim_k \hom_{\mc{O}_X}(\mc{I}_n, \mc{O}_X/\mc{I}_n) = \infty$ and so 
$$\dim_k \lim_{n\to \infty}
 \hom_{\mc{O}_X}(\mc{I}_n, \mc{O}_X/\mc{I}_n) = \infty.$$ 
 Thus the map $\phi_2$ in \eqref{chitwo}
has an  infinite dimensional cokernel, and we are done.
\end{proof}

In the proof of Theorem~\ref{chi_1, not chi_2} we were careful to 
show that, for  $\mc{N} = \mc{R}$,
 \eqref{nth map} is an isomorphism for all $m \geq 0$ and hence that
 \eqref{max tors ext} is an isomorphism.
 This gives the following corollary.
\begin{corollary}\label{R-equals-sections}
$R \cong 
\bigoplus_{m \geq 0} \Hom_{\rQgr \mc{R}}(\pi(\mc{R}), \pi(\mc{R})[m])$.\qed
\end{corollary}

One of the main aims of \cite{AZ1} was to show that, up to a finite dimensional
vector space, one can recover a connected graded noetherian ring $S$ from
$\rqgr S$ whenever $S$ satisfies $\chi_1$. Conversely, given an appropriate
category $\mc{C}$ with an ample shift functor then one  can construct a ``ring
of sections'' $S$ with $\rqgr S\simeq \mc{C}$.  It is almost immediate from
Theorem~\ref{chi_1, not chi_2} that these results do apply to $R$. The reader
should note, however, that  the definitions of ampleness and rings of sections
from \cite{AZ1}  are not the same as the ones used in this paper. Instead,
they are defined as follows:

\begin{definition}\label{artin-zhang2}
 Let $S$ be a noetherian 
connected graded $k$-algebra.  Then the shift functor 
$s:\mc{P}\to\mc{P}[1] $ on $\rqgr S$ is
called \emph{Artin-Zhang ample} provided:
\begin{enumerate}
\item[(a)] For all $\mc{M}\in \rqgr S$ there
 is an epimorphism  $\bigoplus_{i=0}^m 
\pi(S)[-\ell_i]\twoheadrightarrow \mc{M}$ for some $\ell_i\geq 0.$

\item[(b)] For every epimorphism $f : \mc{M} \twoheadrightarrow \mc{N}$
in $\rqgr S$, and for all $n\gg 0$,
the induced map $\hom_{\rqgr S}( \pi(S), \mc{M}[n]) \to 
\hom_{\rqgr S}( \pi(S), \mc{N}[n])$ is surjective.
\end{enumerate}
\end{definition}

One also has a natural functor, which we write as $\Gamma_{AZ}$ to distinguish
it from the global sections functor
$\Gamma$, from $\rqgr S\to \rGr S$ defined by
\begin{equation}\label{artin-zhang21}
\Gamma_{AZ}(\mc{M}) = \bigoplus_{n\geq 0} \hom_{\rqgr S}( \pi(S), \mc{M}[n]).
\end{equation}
Combining the results of this section with those from \cite{AZ1} we obtain:

\begin{corollary}\label{artin-zhang3}
 Keep the hypotheses from \eqref{global-convention2}.
Then:
\begin{enumerate}
\item[(1)]  $\Gamma_{AZ}(\pi_R(R)) = R$. 

\item[(2)]   The shift functor $s$ is Artin-Zhang ample in $\rqgr R$.
\end{enumerate}
\end{corollary}

\begin{proof} (1) Using 
Theorem~\ref{general-ampleness}, this is just a restatement of
Corollary~\ref{R-equals-sections}.

(2) By \cite[Theorem~4.5]{AZ1} and part~1, 
this is equivalent to $R$ satisfying $\chi_1$.
\end{proof}

As a final application of Theorem~\ref{chi_1, not chi_2} we
 note that, by \cite[Theorem~4.2]{YZ}, Theorem~\ref{chi_1, not chi_2}
implies that $R$ does not have a balanced dualizing complex, in the sense of
Yekutieli \cite{Ye}. By \cite{Jg} and Theorem~\ref{fin-cohom-dim}, below, it
does however have a dualizing complex in the weaker sense of \cite{Jg}.

\section{Homological and cohomological dimensions}
\label{cat equivs}

The hypotheses of Assumptions~\ref{global-convention2} remain in force 
in this section.  We will continue to study the homological properties
of  $\rQgr R\simeq \rQgr\mc{R}$ by showing that this category always has finite
cohomological dimension and even has finite homological dimension when $X$ is
smooth. This proves part~6 of Theorem~\ref{mainthm} from the introduction.
 In some sense this result is not surprising: it is not difficult
to reduce both results to a consideration of Goldie torsion modules, in
which case Theorem~\ref{GT equiv} can be applied. 

We recall the definitions of these concepts: 
\begin{definition}
\label{cohomdim-index} 
 The \emph{global dimension} of $\rQgr \mc{R}$ (or $\rQgr R$) is 
 defined to be
\[
\gld( \mc{R}) = \sup \{i \mid \ext^i_{\rQgr \mc{R}}
(\mc{M}, \mc{N}) \neq 0\ \text{for some}\ \mc{M}, 
\mc{N} \in \rQgr \mc{R} \}. 
\]
The \emph{cohomological dimension} of $\rQgr \mc{R}$ (or $\rQgr R$) is  
defined to be 
$$\cd(\rQgr\mc{R}) = \sup \{\cd(\mc{N})\ \mid \ \mc{N}\in \rQgr \mc{R}\},$$
where $\cd(\mc{N}) = 
 \sup \{i \mid \coH^i(\mc{N}) \neq 0.\} $
\end{definition}
 Cohomological dimension is  just defined for $\rqgr R$ in \cite{AZ1} but,
 by \cite[Proposition~7.2(4)]{AZ1}, this is equivalent to our definition.

\begin{theorem}\label{fin-cohom-dim} Keep the hypotheses from 
\eqref{global-convention2}. Then 
$\rQgr \mc{R}$ has finite cohomological dimension.
Indeed,  $\dim X-1 \leq \cd(\rQgr \mc{R})\leq \dim X$.
\end{theorem}

\begin{proof} 
In order to give the upper bound,  what we need to prove is that 
$\coH^i(\mc{M}) = \ext^i_{\rQgr\mc{R}}(\pi(\mc{R}),\mc{M})=0$
 for all $i>d=\dim X\geq 2$ 
and all $\mc{M}\in \rQgr \mc{R}$. We will drop $\pi$ in the proof.
Let $\mc{E}$ denote the injective hull of $\mc{M}$ in $\rQgr\mc{R}$.
By Theorem~\ref{general-ampleness}  and \cite[(7.1.4)]{AZ1},
injective hulls in $\rQgr\mc{R}\simeq\rQgr R$ are induced from those in 
$\rGr R$.  This implies that $\mc{E}/\mc{M}\in \GT\rQgr\mc{R}$, 
since the analogous statement holds in $\rGr R$. Since  
  $\ext^i_{\rQgr \mc{R}}(\mc{R},\mc{M}) =
\ext^{i-1}_{\rQgr \mc{R}}(\mc{R},\mc{E}/\mc{M})$, 
it suffices to prove that 
$\ext^{j}_{\rQgr \mc{R}}(\mc{R},\mc{N})=0$
for $j\geq \dim X$ and $\mc{N}\in\GT \rQgr\mc{R}$.
  By \cite[Proposition~7.2(4)]{AZ1}, we may also 
assume that $\mc{N}\in \GT \rqgr\mc{R}$.

By Lemma~\ref{Goldie tors facts}, we may write 
$\mc{N} = \bigoplus_{n\geq 0} \mc{F}\otimes \mc{L}^{\otimes n}_\sigma$,
 for some $\mc{F}\in \GT \OO_X\catmod.$ Now take an injective resolution
 $\mc{F}\to \mc{E}^{\bullet}$  of 
 $\mc{F} $ in $\OO_X\catmod$ and observe that, since 
 $\mc{F}\in \GT\OO_X\catmod$, so is each $\mc{E}^\ell$. Thus, 
 by Theorem~\ref{GT equiv}
 this induces an injective resolution  
 $ \mc{N}\to \bigoplus_{n \geq 0} \mc{E}^{\bullet}\otimes\mc{L}^{\otimes n}_\sigma$  of 
 $\mc{N} $ in $\rQgr \mc{R}.$
 Thus, for all $j\geq \dim X$ we have

 $$\begin{array}{rll}
\quad \ext^j_{\rQgr \mc{R}}(\mc{R},\,\mc{N}) &=
 h^j\hom_{\rQgr \mc{R}}(\mc{R},\, 
 \mc{E}^{\bullet}\otimes\mc{L}_\sigma^{\otimes n}) \\
\noalign{\vskip 7pt} &=
 h^j{\displaystyle\lim_{n\to \infty}}\hom_{\OO_X}(\mc{I}_n,\, 
 \mc{E}^{\bullet})&\qquad\text{by Lemma~\ref{Hom in Qgr}(1)} \\
\noalign{\vskip 7pt} &=
 {\displaystyle\lim_{n\to \infty}}\, h^j\hom_{\OO_X}(\mc{I}_n,\, 
 \mc{E}^{\bullet})&
\\
 \noalign{\vskip 7pt}&=
 {\displaystyle\lim_{n\to \infty}}\ext^j_{\OO_X}(\mc{I}_n,\, 
 \mc{F}) = 0 &\qquad \text{by Lemma~\ref{exts from points}(3)}.
 \end{array}$$

Thus, $\cd(\rQgr \mc{R})\leq \dim X$.

For the other direction, by Proposition~\ref{L no matter}
it suffices to prove the result for $\mc{R}=\mc{R}(X,c,\OO_X,\sigma)$.
By Lemma~\ref{smooth}, $X$ is smooth at the point $x=c_{-1}$
and so $\OO_{X,x}$ has global dimension equal to $d=\dim X$. 
Thus
$\ext^d_{\OO_{X,x}}(k(x),k(x))\not=0$ (use, for example,
\cite[Exercise~9.25 and Corollary~9.55]{Rt}).
Using the arguments at the beginning of the proof of  
Lemma~\ref{exts from points} this implies that 
$\ext^d_{\OO_{X}}(k(x),k(x))\not=0$.

Now notice that $\GT \mc{O}_X\catMod$ is a localizing 
subcategory of $\mc{O}_X\catMod$ 
which is closed under essential extensions and hence injective hulls. 
 The same is true
of the subcategory $\GT \rQgr \mc{R}$ of $\rQgr \mc{R}$.  It follows 
that $\ext$ groups
involving only Goldie torsion objects may be calculated entirely 
within the Goldie torsion subcategories.
Then by Theorem~\ref{GT equiv} and Remark~\ref{GT equiv1} we have  
$\ext^d_{\rQgr \mc{R}}(\widetilde{x},\widetilde{x}) 
\cong \ext^d_{\OO_{X}}(k(x),k(x)) \neq 0$. 
Since we have taken $\mc{L}=\OO_X$, Definition~\ref{invert-defn1} and
Proposition~\ref{invertibility21} 
 provides a short exact sequence 
$0\to \mc{R}[+1]\to\mc{R}\to \widetilde{x}\to 0$. 
This induces the exact sequence
$$\ext^{d-1}_{\rQgr \mc{R}}(\mc{R}[+1], \widetilde{x})\to
\ext^{d}_{\rQgr \mc{R}}(\widetilde{x}, \widetilde{x})\to
\ext^{d}_{\rQgr \mc{R}}(\mc{R}, \widetilde{x}).$$

Since $\wt{x}$ is Goldie torsion,  
$\ext^{d}_{\rQgr \mc{R}}(\mc{R}, \widetilde{x}) = 0$ by the first part of the
proof.  Therefore, 
$\coH^{d-1}(\widetilde{x}[-1])\cong 
\ext^{d-1}_{\rQgr \mc{R}}(\mc{R}[+1], \widetilde{x})
\not=0$.
\end{proof}

\begin{corollary}\label{gldim}
If $X$ is smooth, then 
 $\rQgr \mc{R}$ has finite global dimension.
 Indeed,  $\dim X\leq \gld(\mc{R})\leq 1+\dim X$.
\end{corollary}

\begin{proof} In the proof of 
 Theorem~\ref{fin-cohom-dim} we show that $\ext^d_{\rQgr
 \mc{R}}(\widetilde{c}_{-1},\widetilde{c}_{-1})\not=0$, which gives the
 required lower bound.
 
 In order to prove the upper bound, we
  need to prove that $\ext^i_{\rQgr \mc{R}}(\mc{M}, \mc{N}) = 0$
for all $\mc{M}, \mc{N} \in \rQgr \mc{R}$ and all $i \geq \dim X + 1$. The
first step is  to reduce to the case of noetherian objects.  One can assume
that $\mc{M}$ is noetherian by the usual  module-theoretic argument
\cite[Theorem~9.12]{Rt} provided one replaces  Baer's criterion by
\cite[Lemma~1, p.136]{Gr}. Now that  $\mc{M}$ is noetherian, the functor
$\ext^i_{\rQgr \mc{R}}(\mc{M}, -)$ commutes with direct limits
\cite[Proposition~7.2(4)]{AZ1}  and so we may assume that $\mc{N}$ is also
noetherian.

So, assume that $\mc{M}$ and $\mc{N}$ are noetherian.
By Lemma~\ref{goldie-subfactors}, we may assume that 
 $\mc{M}$ is either a shift of  $\mc{R}$ or  a 
 Goldie torsion module. If $\mc{M}=\mc{R}[r]$, then 
$\ext^i_{\rQgr \mc{R}}(\mc{M},\mc{N}) 
=\ext^i_{\rQgr \mc{R}}(\mc{R},\mc{N}[-r])=0$,
 by  Theorem~\ref{fin-cohom-dim}.  Thus we may assume 
that $\mc{M}\in \GT\rqgr \mc{R}$.

As in the proof of Theorem~\ref{fin-cohom-dim}, it  
suffices to prove that 
 $\ext^j_{\rQgr \mc{R}}(\mc{M},\mc{N})=0$ for 
 all $\mc{N}\in \GT\rqgr \mc{R}$ 
 and all $j >\dim X$.
By Theorem~\ref{GT equiv}, we may write 
$\mc{M}= \bigoplus_{n \geq 0} \mc{F} \otimes \mc{L}_{\sigma}^{\otimes n}$ 
and $\mc{N} =  \bigoplus_{n \geq 0} \mc{G} 
\otimes \mc{L}_{\sigma}^{\otimes n}$ for some
$\mc{F},\mc{G}\in \GT\OO_X\catmod$.  Moreover, the   argument given in 
the proof of Theorem~\ref{fin-cohom-dim} shows that we may
calculate $\ext$ inside the Goldie torsion subcategories and so 
$\ext^j_{\rQgr \mc{R}}(\mc{M},\mc{N})=  \ext^j_{\OO_X}(\mc{F}, \mc{G}).$
As $X$ is smooth,  $\OO_X\catmod$ has homological dimension $\dim X$
\cite[Exercise~III.6.5 and Proposition~III.6.11A]{Ha}. 
 Thus $ \ext^j_{\OO_X}(\mc{F},  \mc{G})= 0$ for all $j>\dim X$
 and $\gld(\mc{R})\leq 1+\dim X$. 
 \end{proof}

In both Theorem~\ref{fin-cohom-dim}
and Corollary~\ref{gldim}
we conjecture that the correct dimension is $\dim X$.

Curiously, there seems to be no known example of a noetherian connected graded
 algebra that
does {\it not} have finite cohomological dimension.  We presume that such
examples do exist, if only because the standard proofs that commutative
varieties have finite cohomological dimension clearly do not work in a
noncommutative setting.

\section{Generic flatness}\label{section-genflat}

The hypotheses from Assumptions~\ref{global-convention2} will be assumed
throughout this section. Recall that a $k$-algebra $A$ is  \emph{strongly
noetherian} if $A \otimes_k C$  is noetherian for all commutative noetherian
$k$-algebras $C$.  The strong noetherian condition was introduced in
\cite{ASZ,AZ2} where it is shown that strongly noetherian graded algebras have
a remarkable number of nice properties. Notably, \cite[Theorem~0.1]{ASZ}
shows that they satisfy generic 
flatness, as defined in the introduction. In this section  we show that 
$R=R(X,c,\mc{L},\sigma)$ fails  generic flatness in a quite dramatic way:
for \emph{any} open affine subset $V\subseteq
X$ it fails for the algebras $C=\OO_X(V)\subset A=R\otimes_kC$  and the module
$M=\mc{R}(V)$. We will also  use this to 
construct an explicit noetherian ring $C$ such that $R\otimes_kC$ is not
noetherian.

\begin{lemma}
\label{sections} Let $V$ be an open affine subset of $X$  and 
set $C=\calO_X(V)$. Then $M = \calR(V)$ is 
a finitely generated 
$(C, R)$-bimodule.  Equivalently, (after identifying $C^{\op}$ with
$C$)
$M$ is a finitely generated right $R \otimes_k C$-module.
\end{lemma}

\begin{proof} 
Trivially $M$ is a graded left $C$-module.  Write 
$\mc{J}_n = \mc{I}_n \otimes \mc{L}_n$ for $n \geq 0$; thus  
 $\mc{R}_n={}_1(\mc{J}_n)_{\sigma^n}$. 
Then the right $R$-module structure
on $M$ is the natural one induced by the maps 
\[
\mc{J}_n(V) \otimes \mc{J}_m(X) \overset{1 \otimes (\sigma^n)^*}{\lra} 
\mc{J}_n(V) \otimes \mc{J}_m^{\sigma^n}(X) \lra
 \mc{J}_n(V) \otimes \mc{J}_m^{\sigma^n}(V) \lra \mc{J}_{n+m}(V)
\]
for $m, n \geq 0$. 
The commutativity of the two actions is clear.

By Theorem~\ref{general-ampleness}, there exists $n_0$ such that $\calR_n$ is
generated by its sections
for all $n \geq n_0$. Thus the 
natural map $R_{\geq n_0} \otimes C \to M_{\geq n_0}$ 
is a surjective $R \otimes C$-module homomorphism.
Since $R$ is noetherian, $R_{\geq n_0}$ is a finitely generated right
$R$-module and hence
$M_{\geq n_0}$ is a finitely generated $R \otimes C$-module.
Finally, as  $\bigoplus_{i=0}^{n_0-1}M_i$ is a finitely generated left
$C$-module,  $M$ is indeed a finitely generated
$R \otimes C$-module.
\end{proof}

\begin{theorem} 
\label{not strong noeth} Keep the hypotheses of
\eqref{global-convention2}. Let $V$ be any open affine subset of $X$ 
 and write $C = \mc{O}_X(V)$ and
 $M = \mc{R}(V)$.  
Then $M$ is a finitely generated right $R\otimes_k C$-module which is not 
generically flat over $C$.

 Thus, $R$ is neither strongly right noetherian nor strongly left noetherian.  
\end{theorem}

\begin{remarks}\label{not strong noeth1} (1)  One can make the theorem more
precise by exactly determining the maximal ideals  $\p$ of $C$ at which $M$ is
not flat. Indeed,  let $\m_x$ denote the maximal ideal of $C$ corresponding to
the  closed point $x\in X$.  Then $M$ fails to be flat at precisely the maximal
ideals  $\m_{c_i}$ for $i\in P=\{n \geq 0 : c_n\in V\}$. Moreover,   $M_t$ will
not be flat at such an $\m_{c_i}$  whenever $t>i$. Note that, as   
$\{ c_i \}_{i\geq 0}$  is critically dense, ${\mathbb N}\smallsetminus P$ 
is a finite set.

(2) If one wishes to work more scheme-theoretically, then 
Theorem~\ref{not strong noeth} generalizes naturally to one describing
$\mc{R}$  as a sheaf of right modules over  $R\otimes_k\OO_X$.

(3) The theorem  proves Theorem~\ref{mainthm}(3) and 
Proposition~\ref{main-thm2} from the introduction.
\end{remarks}

\begin{proof} 
If $R$ is strongly right noetherian, then so is $R\otimes_k C$, since
$C$ is a finitely generated $k$-algebra. 
Thus once we prove that $M$ is not generically flat over $C$, 
it follows from \cite[Theorem~0.1]{ASZ} that
 $R$ cannot be strongly right noetherian.
It then follows  from \eqref{rees opposite} 
that $R$ is not strongly left noetherian.

 By  Lemma~\ref{sections} $M$ is   a finitely generated $R \otimes C$-module. 
Consider the short exact sequence of sheaves  
$$ 
0\longrightarrow \mc{I}_n \otimes \mc{L}_n  \longrightarrow 
\mc{L}_n  \longrightarrow 
\bigoplus_{i=0}^{n-1} k(c_i)  \longrightarrow 0.    
$$
Localizing this sequence at the open set $V$ gives
\begin{equation}\label{exact33}
 0  \longrightarrow M_n  \longrightarrow \mc{L}_n(V) 
 \longrightarrow \bigoplus_{i=0}^{n-1} k(c_i)(V)  \longrightarrow 0.
\end{equation}
Since $\mc{L}_n$ is locally free of rank 1, 
for any $\p \in \spec C$ the module $(\mc{L}_n(V))_\p$ will be 
isomorphic to $C_\p$.

Let $P$ be the set defined by Remark~\ref{not strong noeth1}
 and let $\p$ be any prime ideal  in $\spec C$ that is
not equal to   $\m_i=\m_{c_i}$ for $i\in P$. In particular, 
$\p\not=\m_j$, for  $j\geq 0$. Thus, for any $n$, 
localizing \eqref{exact33} 
at $\p$ shows that
$(M_n)_\p \cong C_\p$ and hence that $M_\p$ is a flat $C_\p$-module.

On the other hand, if $\p = \m_i$ for $i\in P$ and $n-1 \geq i$, 
then localizing \eqref{exact33}  at $\p$ gives
the exact sequence
$$ 0 \longrightarrow(M_n)_{\m_i} \longrightarrow C_{\m_i} 
\longrightarrow C_{\m_i}/\m_i C_{\m_i} \longrightarrow 0.
$$
Thus $(M_n)_{\m_i}\cong \m_i C_{\m_i}$.
  Suppose that $\m_iC_{\m_i}$ is a flat $C_{\m_i}$-module.  Then it is
free, and hence principal.  By the Krull
principal ideal theorem,  $\m_iC_{\m_i}$ must have height one, contradicting
 the fact that $\dim C_{\m_i} \geq 2$.  Thus  $(M_n)_{\m_i}$ is never
 flat for  $n \gg 0$ and so 
$M_{\m_i}$ is not a flat $C_{\m_i}$-module  when $i\in P$.

Now let  $f\in C$ be non-zero. Then $C_f=\OO_X(W)$ for some open affine set
$W\subseteq V$. Since the conclusion of the last paragraph is independent of the
choice of $V$, this implies that $M_f = \mc{R}(W)$ cannot be 
a flat $C_f$-module.

Finally, we have shown that $M_\p$ is 
a flat $C_\p$-module if and only if $\p \not= \m_i$
for $i\in P$. This justifies the assertions in 
Remark~\ref{not strong noeth1}(1).
\end{proof}

A natural question raised by Theorem~\ref{not strong noeth} is to identify a
commutative noetherian ring $C$ for which $R\otimes C$ is not noetherian.  
This was achieved in \cite{Ro} for the special case $X = \mb{P}^t$, but without
an explanation for why the particular ring $C$ given there, an infinite affine
blowup of affine space, was a natural choice.  We show next how the results of
\cite{ASZ} lead inevitably to similar choices of $C$ in general. To set this
up, we need to discuss some commutative  constructions from \cite{ASZ}.

Let $C$ be a finitely generated commutative $k$-algebra that is a domain
with fraction field $F$, and
let  $c$ be a closed smooth point of $\spec C$ with associated maximal
ideal $\p$.   The \emph{affine blowup} of $\spec C$ at $c$ is $\spec C'$, where
the ring $C'$ is  formed as follows:  write
$\p=\sum_{i=0}^rx_iC$ with $x_0 \not \in \p^2$ and 
define $C' = C[x_1 x_0^{-1}, x_2  x_0^{-1}, \dots, x_r x_0^{-1}]\subset F$.
Here $x_0$ will be called the denominator for the blowup.
  Now suppose we are given an infinite sequence  of
distinct smooth closed points $c_1, c_2, \dots$ of $\spec C$  corresponding to maximal
ideals  $\p_1, \p_2, \dots$, and for each $i \geq 1$ suppose that we can find 
a denominator 
$x^{(i)}=x^{(i)}_0$ such that
\begin{equation}
\label{denominators} x^{(i)} \in \p_i \setminus {\p_i}^2\
 \text{but}\ x^{(i)} \not \in 
\p_k\ \text{for}\ k \neq i. 
\end{equation}
Then we can successively blow up each of the points $c_i$, giving 
a sequence of rings
$
C \subseteq C_1 \subseteq C_2 \subseteq \dots
$
The ring $\wt{C} = \bigcup C_i$ is called the \emph{infinite affine blowup} of
$C$ at  the points $c_i$ (with respect to the particular choice of
denominators $x^{(i)}$). 
  Let $\rho:  \spec \wt{C} \to \spec C$ be the blowup map,
induced by the inclusion $C \to \wt{C}$.  For each $i \geq 1$ the ideal $\eta_i
= \p_i \wt{C}$ is a prime ideal of $\wt{C}$  that corresponds to the
exceptional divisor $\rho^{-1}(c_i)$ of $\spec \wt{C}$ \cite[Lemma~1.3]{ASZ}.

Returning to our specific scheme $X$, 
let $V$ be an open affine subset of $X$
and set $C = \mc{O}_X(V)$. Define $P$ by Remark~\ref{not strong noeth1} and 
 write $\p_i=\m_i$ for  the
maximal ideal of $C$ corresponding to the point $c_i$ for $i\in P$.
Recall that the $c_i$
are smooth by Lemma~\ref{smooth}.  
Because the  sequence $\{c_i : i\in P\}$ is critically dense, it is
possible to choose an infinite   subsequence 
$\mc{C}'=\{c_{i} : i\in I\subseteq P\}$ which have 
 denominators  satisfying \eqref{denominators}
\cite[Proposition~1.6]{ASZ}, and so the  infinite blowup $\wt{C}$ of $C$ at
that subsequence $\mc{C}'$ is well defined. By   \cite[Theorem~1.5]{ASZ} 
 $\wt{C}$ is   noetherian.  If
 we invert $\{f \in \wt{C}\smallsetminus \bigcup_{i\in I}\p_i\wt{C}\}$,
  we obtain a further localization $D$ of $\wt{C}$ that  is a
Dedekind domain \cite[Proposition~2.8]{ASZ}.

\begin{theorem} \label{non-strong}
Assume that \eqref{global-convention2} holds and set
$R=R(X,c,\mc{L},\sigma)$ and let $\wt{C}$ and $D$ be defined as above. Then
$R \otimes_k \wt{C}$ and $R \otimes_k D$ are not noetherian rings.
\end{theorem}

\begin{proof} Set $M =\mc{R}(V)$.   By Remark~\ref{not strong noeth1}(1) $M_\p$
is not flat for exactly the maximal ideals $\p_i \in C$ which  correspond to
the points in $\{c_i : i\in P\}$. Following the proof of
\cite[Theorem~2.3]{ASZ}, we see that   $D$ has maximal ideals  $\{\mu_j =\p_jD 
: j\in I\}$ and that $M \otimes_C  D_{\mu_j}$ is never a flat
$D_{\mu_j}$-module.  
 Since  each $D_{\mu_j}$ is a DVR, this means that $M \otimes_C 
D_{\mu_j}$ is not even a torsion-free   $D_{\mu_j}$-module
\cite[Lemma~3.3(4)]{ASZ}. Now the points $\mu_j$ are  critically dense in
$\spec D$.  So if $M \otimes_C D$ were a noetherian $R
\otimes_k D$-module,  then  \cite[Lemma~3.3(2)]{ASZ} would imply that  $M
\otimes_C D$ would be torsion-free at all but finitely many  closed points of
$\spec D$. This contradicts our construction and proves that $M \otimes_C D$ is
not noetherian. However it is a finitely generated module over  $R\otimes_k
D\cong R\otimes_k C\otimes_C D$ simply because $M$  is a finitely generated
$R\otimes_k C$-module. Thus  $R \otimes_k D$ is not  noetherian which, since
$D$ is a localization of $\wt{C}$, implies that $R  \otimes_k \wt{C}$ is also 
not noetherian. \end{proof}

\section{Point modules}\label{section-point}
The hypotheses from
Assumptions~\ref{global-convention2} will remain in force throughout the
section. Here we discuss 
 the point modules for $R=R(X,c,\mc{L},\sigma)$
  and show in particular that they are not
parametrized by a projective scheme and that the shift functor 
does not induce an automorphism on the set of point modules.  
This is in marked contrast to the behaviour of strongly noetherian 
  rings and so provides further proofs of
the fact that $R$ is not strongly noetherian. 

To set this  in context, we begin by reviewing the work of \cite{AZ2}
and answering a special case of one their conjectures. Let $S=\bigoplus_{n\geq
0} S_n$ be a connected graded $k$-algebra. Given a commutative
$k$-algebra $C$,   write $S_C = S\otimes_kC=\bigoplus (S_n\otimes_kC)$,
regarded as a graded $C$-algebra. Fix a finitely generated graded $S$-module
$P$ and a Hilbert function  $h: \mathbb{N} \to \mathbb{N}$. For a finitely
generated commutative $k$-algebra $C$,  write $\mathcal{P}_h(C)$ for the set of
isomorphism classes of  graded factors $V$ of $P_C=P\otimes_SS_C$ with the
property that  each $V_n$ is a flat $C$-module of constant rank $h(n)$.  Given
a map of finitely generated algebras $C \to D$,  one gets a map from
$\mc{P}_h(C)$ to $\mc{P}_h(D)$ via $M \mapsto M \otimes_C D$. 
Following \cite[Section~E4]{AZ2}, for an
arbitrary commutative $k$-algebra $C$ we let
$\mathcal{P}_h(C)$ denote the direct limit of the sets $\mc{P}_h(C')$ as $C'$
varies over the finitely generated subalgebras of $C$.  In this way $\mc{P}_h$
becomes a functor from  $(\rings)$ to $(\sets)$, where $(\rings)$ denotes the
category of commutative  $k$-algebras. 

One of the main aims of \cite{AZ2} was to show that for any $h$ the functor 
$\mathcal{P}_h$ can be represented by a projective scheme.  The relevant
definitions  are just like the commutative case: given a scheme $Y$ over $k$, 
its \emph{functor of points}  $h_{Y}: (\rings) \to (\sets)$ is defined by $C
\mapsto \Morph(\spec C, Y)$, where $\Morph$ denotes morphisms in the category
of $k$-schemes.  By Yoneda's lemma this  induces an embedding of the category
$(\schemes)$ as a full subcategory  of the category $(\fun)$  of covariant
functors from  $(\rings)$ to $(\sets)$ (see \cite[Section~VI]{EH} for more
details.) A functor $F \in (\fun)$ that is equal to $h_{Y}$ for some scheme 
$Y$ is said to be \emph{represented} by $Y$.

For strongly noetherian algebras one has the following result:

\begin{theorem}\label{AZ-param} {\rm(}\cite[Theorem~E4.3]{AZ2}{\rm)}
Let $S$ be a strongly noetherian connected graded $k$-algebra
and fix a finitely generated graded 
$S$-module $P$ and a Hilbert function $h$.
Then the functor $\mathcal{P}_h$ is represented by a projective scheme $Y$.
\end{theorem}

The paper \cite{AZ2} also proves an analogous but weaker result for
parametrizing objects in $\rqgr S$ with a given Hilbert series.  The
definitions are as follows: Given $P$ and $h$ as above, for any finitely
generated commutative $k$-algebra $C$ we let $\mathcal{P}_h^{\mathrm{qgr}}(C)$
denote the set of isomorphism classes of graded factors $V$ of $P_C$ with the
property that  $V_n$ is a flat $C$-module of constant rank $h(n)$ for all $n\gg
0$.  As before, the definition  is extended to all $C \in (\rings)$ by taking
limits, and thus $\mc{P}_h^{\mathrm{qgr}}$ becomes a functor from $(\rings) \to
(\sets)$.\footnote{It is not clear what is the best version of flatness to use
in this definition. Specifically, \cite{AZ2} uses the formally weaker notion
that $V$ be flat in $\rQgr S_C$  rather than requiring that the $V_n$ be flat.
The two notions coincide for strongly noetherian algebras by
\cite[Lemma~E5.3]{AZ2}. Fortunately, the distinction is not significant since 
we will only need to make computations when $C$ is a noetherian domain. In this
case  both notions of flatness follow automatically from the fact that the
$V_n$ have constant rank.}

By \cite[Theorem~E5.1]{AZ2}, if $S$ satisfies the \emph{strong
$\chi$ condition} \cite[p.346]{AZ2}, then  the functor
$\mathcal{P}_h^{\mathrm{qgr}}$ is represented by a scheme $Y$ that is a
countable union of projective closed subschemes. Artin and Zhang
conjecture that $Y$ is actually a projective scheme and 
the first result of the section   shows that this conjecture is true for 
point modules. These are  defined as follows.
Suppose that  $S$ is a
connected graded ring, generated in degree one, $P=S$ and 
$h$ is the function $1$.  
Then we will write $\ptfn$ for $\mathcal{P}_h$ and $\qptfn$ 
for $\mathcal{P}^{\mathrm{qgr}}_h$.  Elements of the  
set $\ptfn(C)$ will be called
\emph{point modules over $ S_C$} while elements of 
 $\qptfn(C)$ will be called
\emph{point modules in $\rqgr S_C$.} 
 If $h$ is the function $h(t)=1$ for $0\leq
t\leq d-1$ but $h(t)=0$ for $t\geq d$, then the modules in   $\mathcal{P}_h(C)$
are defined to be the \emph{truncated point modules of length $d$} over $C$.

By \cite[Corollary~E4.4]{AZ2}, if $S$ is strongly 
noetherian then there exists an
integer $t$ (independent of the commutative $k$-algebra $C$) 
such that:
\begin{enumerate}
\item[(i)] \emph{For all $C$ and   $V\in \ptfn(C)$, 
$\mathrm{Ker}(S\to V)$ is generated in degrees $\leq t$.}

\item[(ii)] \emph{Conversely, given a truncated point 
module $V'$ of length $d\geq
t$, then $V'$ is the homomorphic image of a unique point module $V$.}
\end{enumerate}
The significance of these results is that it is easy to show that the truncated
point modules  of length $d$ are parametrized by a projective scheme $Y_d$ (see
\cite[Lemma~E4.6]{AZ2}). It then follows from (i) and (ii) that, for $d\geq t$,
the scheme $Y_d\cong Y_{d+1}$  parametrizes
the point modules.  As we next show, this scheme also represents $\qptfn$.

  This result was proved jointly by Michael
Artin and the third-named author and we are grateful to Artin for letting us
include it here.

\begin{proposition}\label{shift-aut}  Let $S$ be a  noetherian, connected
graded $k$-algebra that is generated in degree one. Assume that $S$ is strongly
noetherian or, more generally, that (i) and (ii) hold on both the left and
right  so that the left, respectively right, point modules are  
parametrized by a projective scheme $Y^\ell$, respectively $Y^r$. Then:
\begin{enumerate}
\item[(1)] For any commutative noetherian $k$-algebra $C$ the shift functor $s:
M\mapsto M[1]_{\geq 0}$ is an automorphism of the set 
of right $S_C$-point modules.

\item[(2)] The shift functor $s$ induces an automorphism
of both $Y^r$ and $Y^\ell$.

\item[(3)] The scheme $Y^r$ also represents the functor
$\qptfn$.
\end{enumerate}
\end{proposition}
\begin{proof} (1) As the base ring $C$ is fixed, 
we may write $S$ for $ S_C$ without confusion.  
By (i) and (ii),  every truncated right $S$-point
module $M=\bigoplus_{i=0}^{r-1}M_i$ of length $r\geq t$ is the factor
$N/N_{\geq r}$ for a unique point module $N$. Thus, to prove the result it
suffices to show that there exists a unique shifted truncated $S$-point
module $M' =\bigoplus_{i=-1}^{r-1}M'_i$ such that $M=M'_{\geq 0}$.

As $C$ is noetherian, each $M_i $ is a projective $C$-module of constant rank
one.  Consider the Matlis dual $M^\vee = \bigoplus_{i=1-r}^0M^\vee_i$ of $M$;
thus  $M^\vee_i=\hom_C(M_{-i},C)$ for each $i$. Clearly $M^\vee $ is a left
$S$-module for which each  $M^\vee_i$ is a projective $C$-module of constant
rank one. We claim that $M^\vee_{1-i} = S_1M^\vee_{-i}$ for  $1 \leq i \leq
r-1$.  It suffices to prove this locally, so assume that $C$ is local.  
Then each $M_i$ is free, say $M_i=m_iC\cong C$. As $S$ is generated in degree one,
$M_i=M_{i-1}S_1=m_{i-1}S_1$ and $m_i=m_{i-1}s$ for some  $s\in S_1$.   Let 
$\theta$ be the generator of $ M_{-i}^\vee$; thus $m_i^\theta=1$. Then
$\phi=s\theta\in M_{1-i}^\vee$ satisfies  $m_{i-1}^\phi=(m_{i-1}s)^\theta=1$.
In other words,  $\phi$ is the generator of $M_{1-i}^\vee $ and $\phi\in
S_1M_{i}^\vee$. This proves the claim.

By the claim $M^\vee = SM^\vee_{1-r}$ and so it is  the shift of a
truncated left point module. By  hypothesis (ii), $M^\vee$ is a
homomorphic image of a unique shifted point module and so  
there exists a unique shifted truncated   
point module $L=\bigoplus_{i=1-r}^1L_i$ such that
$M^\vee = L/L_1$. Taking Matlis duals, again, gives the 
required module $M'=L^\vee$.
The uniqueness of $L$ implies the uniqueness of $M'$.

(2) We prove the statement for $Y^r$ only; the proof for $Y^{\ell}$ is
symmetric.  For each finitely generated  commutative $k$-algebra $C$, in part~1
 we proved that $s$ induces a bijection  from the set $\ptfn(C)$ to
itself.  It is easy to check that these bijections  are functorial, so for any
commutative $k$-algebra $C$, we get an induced bijection from $\ptfn(C)$ to
itself  by taking limits over the finitely generated subalgebras $C'$ of $C$. 
Thus  we have actually defined a natural isomorphism from the functor $\ptfn$
to itself in  the category $(\fun)$.  Since Yoneda's lemma embeds $(\schemes)$
as a full subcategory of $(\fun)$, and $\ptfn$ is represented by the scheme
$Y^r$, we must have a scheme automorphism $\sigma: Y^r \to Y^r$ induced by $s$.

(3) Fix a commutative noetherian ring $C$  and   $\mc{M}=\pi(M)\in \qptfn(C)$.
 As $S_C$ is generated in degree one, we may choose a  tail $N=M_{\geq n}
= \bigoplus_{i\geq n}M_i$ of $M$  such that $N=N_{n}S_C$ and  $N_i$ is a
projective $C$-module locally of rank one for all $i\geq n$. By part~1, $N$ is
the tail $L_{\geq n}$ of a unique point  module $L\in \ptfn(C)$.  Thus the
sets $\ptfn(C)$ and $\qptfn(C)$ are in natural bijection for all finitely
generated $k$-algebras $C$ and the same holds for all $C \in (\rings)$ by
taking  limits.  Thus  $\ptfn(C)$ and $\qptfn(C)$ are naturally
isomorphic functors and   $\qptfn$ is also  represented by $Y^r$. \end{proof}

We now turn to the structure of the point modules  over
$R=R(X,c,\mathcal{L},\sigma)$ and show that $R$ satisfies none of the
conclusions of Theorem~\ref{AZ-param} or Proposition~\ref{shift-aut}. 
To do
this we either have to assume that $R$ is generated in degree one (since that
is required for  the definition of point modules) or to work with general $R$
and use a slightly  more artificial class of $R$-modules.  We will use the
second approach in a way that also includes the first case.

\begin{notation}\label{def-pseudopoint} Fix an open affine subset $U\subset X$
and recall from Theorem~\ref{not strong noeth} that $\mathcal{R}(U)$ is a
finitely generated right $R_{\OO_X(U)}$-module. Fix a finitely generated  
graded free $R$-module $P$ such that $\mathcal{R}(U)$ is a homomorphic image of
$P_{\OO_X(U)}=P\otimes_RR_{\OO_X(U)}$. Now take $h$ to be the constant 
function $1$ and let
$\mathcal{P}=\mathcal{P}_h$ be the corresponding functor.
\end{notation}

We would like to thank Brian Conrad for his help which was invaluable in the
proof of the next result.

\begin{theorem}\label{non-represent} Keep the hypotheses of 
\eqref{global-convention2} and \eqref{def-pseudopoint}.
Then:
\begin{enumerate} 
\item[(1)] The functor $\mathcal{P}$ is 
not represented by any scheme $Y$ of
locally  finite type.

\item[(2)] For $m\gg 0$, set $S=R(X,c,\mc{L}^{m},\sigma)$. 
Then $S$ is generated in degree one but $\ptfn$ is not represented
by  a scheme $Y$ of
locally  finite type.
\end{enumerate}
\end{theorem}

\begin{proof}  (1) 
Assume that $\mathcal{P}$ is represented by the scheme $Y$ of
locally finite type. The intuitive reason for our choice of $\mathcal{P}$ is that
$\mathcal{R}(U)$ is ``almost'' in   $\mathcal{P}(\OO_X(U))$. More precisely,
the choice of $P$ ensures that  $\mathcal{R}(U)$ is a homomorphic image of
$P_{\OO(U)}$. Moreover, if $p\in U\smallsetminus\{c_i\}$ is a closed point
then  $(\mathcal{I}_n\otimes_{\OO_X}\mathcal{L}_n)_p \cong
\OO_{X,p}$, for all $n\geq 0$. Thus, if $C=\OO_{X,p}$ for some such $p$, then 
$\mathcal{R}_p = \mathcal{R}(U)\otimes_{\OO(U)}C\cong
\bigoplus_{i\geq 0} C$ does belong to $\mathcal{P}(C)$. However,  
Theorem~\ref{not strong noeth} and Remark~\ref{not strong noeth1}(1) imply that 
$\mathcal{R}(U)$ is not a flat $\OO(U)$-module; indeed it fails to be flat at
precisely the points  in $U\cap\{c_i\}$.  Thus, $\mathcal{R}(U)$ does not lie
in $\mathcal{P}(\OO(U))$.

Fix a closed point $p\in U\smallsetminus\{c_i: i \in \ZZ \}$ and set
$C=\OO_{X,p}$.  The last
paragraph implies that there exists 
$\theta_p\in \mathcal{P}(C)=\Morph(\spec C, Y)$
corresponding to $\mathcal{R}_p$.
By the definition of locally finite type \cite[p.84]{Ha}, we may pick an open affine 
neighbourhood $V$ of $\theta_p(p)$ in $Y$ of finite type over $k$. 
Then  we get a map of algebras
$\theta_p' : \OO_Y(V)\to  \OO_{\spec C}(\theta_p^{-1}(V))$. 
Since $\theta_p^{-1}(V)$ is an
open set containing $p$, it is necessarily $\spec C$ and so 
 $\mathrm{Im}(\theta_p')\subseteq C$.
Since $\OO_Y(V)$ is a finitely generated $k$-algebra
  and $\OO_X(U)$ is a domain,
 $\theta_p'(\OO_Y(V))\subseteq \OO_X(U')$, for some open set $U'\subseteq U$.
 Since it does no harm to replace $U$ by a smaller open set containing $p$, 
 we may as well assume that $U=U'$.
 In other words, we have extended $\theta_p$ to a map 
 $\widetilde{\theta}_p \in \Morph(U,Y)$ such that $\theta_p=\widetilde{\theta}_p\circ \pi_p$,
  where $\pi_p : \spec C\to U$ is the natural
morphism.

By construction, $\widetilde{\theta}_p$ corresponds to a module 
$M_U\in \mathcal{P}(\OO(U))$ 
with the property that $M_U\otimes_{\OO(U)}C \cong \mathcal{R}_p$.
But $\mathcal{R}(U)$ is a second finitely generated $R_{\OO(U)}$-module 
with $\mathcal{R}(U)\otimes_{\OO(U)}C \cong \mathcal{R}_p$. This local
isomorphism of $R_C$-modules  lifts to an isomorphism
$M_W=M_U\otimes_{\OO(U)}\OO(W) \cong \mathcal{R}(W)$ of $R_{\OO(W)}$-modules,
 for some open affine set
$W\subseteq U$.  By the definition of $\mathcal{P}$, 
$(M_W)_n=(M_U)_n\otimes_{\OO(U)}\OO(W)$ is a flat $\OO(W)$-module for
all $n$. On the other hand, for $n\gg 0$,   Remark~\ref{not strong noeth1}(1)
implies that  $(M_W)_n\cong \mathcal{R}(W)_n$ is \emph{not} flat over
$\OO(W)$.  This contradiction proves (1). 

(2)  Pick an integer 
$M$ by Proposition~\ref{degree 1} and assume that $m\geq M$. 
 Then  $S$ is generated in degree one and   $P=S$ satisfies
the hypotheses of \eqref{def-pseudopoint}. Thus, part~2 follows from part~1.
\end{proof}

With minor changes the proof of Theorem~\ref{non-represent} also shows that the
point modules in $\rqgr R$ are  not parametrized by a scheme of locally
finite type.

\begin{corollary}\label{non-represent3} 
Keep the hypotheses of \eqref{global-convention2}
 and \eqref{def-pseudopoint}.   
Then the functor $\mathcal{P}^{\mathrm{qgr}} = \mathcal{P}^{\mathrm{qgr}}_1$ is 
not represented by any scheme $Y$ of locally  finite type.

Similarly, if $S=R(X,c,\mc{L}^m,\sigma)$  for $m\gg0$, then
 $S$ is generated in degree one but $\qptfn$ is not represented by 
a scheme $Y$ of locally  finite type.\end{corollary}

\begin{remark} Combined with 
Theorem~\ref{non-represent} and Remark~\ref{point-class12},
this proves  Theorem~\ref{mainthm}(5).  
\end{remark}

\begin{proof}
Consider the proof of Theorem~\ref{non-represent}(1). 
In the final paragraph, 
$M_W \in \mathcal{P}(\OO(W))$ and so it  certainly lies in 
$\mathcal{P}^{\mathrm{qgr}}(\OO(W))$. 
However, as $\mathcal{R}(W)_n$ is not flat 
as an $\OO(W)$-module for any $n\gg 0$, no tail $\mathcal{R}(W)_{\geq n}$
of $\mathcal{R}(W)$ is a flat $\OO(W)$-module and hence 
$\mathcal{R}(W)$ cannot belong to 
$\mathcal{P}^{\mathrm{qgr}}(\OO(W))$. 

Thus, the proof of Theorem~\ref{non-represent} can also be used to prove 
the corollary.
\end{proof}

Corollary~\ref{non-represent3}  is in stark contrast to Remark~\ref{GT equiv1}.
To see this, assume that $R=S$, so that 
the functor $\mathcal{P}=\qptfn$
determines point modules.
It  follows from Remark~\ref{point-class12} below that 
the points in $\rqgr R$ are simply the images of the point modules in $\rgr R$.
Thus Remark~\ref{GT equiv1} can be rephrased as saying that, in $\rqgr R$,
 the point modules are in (1-1) correspondence with
the points of $X$---which is definitely a scheme of finite type.
 The way to think of the difference is as follows: If the
point modules in qgr were indeed parametrized by $X$ then, in 
 $\rqgr R_{\OO(U)}$, one would need to find not only the point modules
induced from $R=R_k$ but also a module corresponding to
the immersion $U\subseteq X$. The proof of Corollary~\ref{non-represent3}
 can be interpreted as
saying that one does have indeed a  module corresponding to $U$. 
Unfortunately it is $\pi(\mathcal{R}(U))$ which, as 
 the proof also shows, is  not a point module 
in   $\rqgr R_{\OO(U)}$.

Although the last few results have shown that one cannot parametrize the point 
modules for $R$, they do not say anything about the point modules 
over $R=R_k$ itself. The next two results consider these modules and show that
they also have interesting properties. For any closed point 
$x\in X$, recall the definition of the modules
 $\widetilde{x}=\bigoplus_{i\geq 0} k(x)_{\sigma^i}\in \rqgr\mc{R}$
from \eqref{exceptional}.

\begin{proposition}\label{point-class} Keep the hypotheses of 
 \eqref{global-convention2}.
 Set $\dim X=d$ and write
 $R=R(X,c,\mc{L},\sigma)$. If $x$ is a closed point in $X$, then 
 the $R$-module of global sections $\Gamma_{AZ}(\widetilde{x})$ 
has Hilbert series
$$H(t) = \begin{cases}
\sum_{i=0}^{r} dt^i + \sum_{i=r+1}^\infty t^i & \text{ if }
\ x=c_r\ \text{ for }\ r\geq 0;\\
\noalign{\vskip 7pt}
1/(1-t) & \text{otherwise}.
\end{cases}$$
\end{proposition}
\begin{remark}\label{point-class12}
By  Remark~\ref{GT equiv1}
and the equivalence of categories Theorem~\ref{general-ampleness},
the simple objects in $\rQgr \mc{R}\simeq \rQgr R$ are precisely
  $\{\widetilde{x} : x \text{ a closed point in } X\}$.
Thus, if $R$ is generated in degree one, the proposition implies 
that these objects are precisely the point modules in $\rqgr R$.
\end{remark}

\begin{proof} We remind the reader that $\Gamma_{AZ}$ denotes the 
Artin-Zhang sections functor from \eqref{artin-zhang21}.
By Lemma~\ref{shifts}, $\widetilde{x}[m]=
\bigoplus_i (k(x)^{\sigma^{-m}})_{\sigma^i}
=\bigoplus_i k(\sigma^{m}(x))_{\sigma^i}.$
By Lemma~\ref{Hom in Qgr}(1) and 
Theorem~\ref{general-ampleness},
$$\Gamma_{AZ}(\widetilde{x})_m
= \lim_{n\to \infty}
\hom_{\OO_X}(\mc{I}_n,\, k(\sigma^m(x))) \qquad\text{for}\  m\geq 0.$$

There are two possibilities. If $x=c_j$ with 
$0\leq j-m\leq n-1$, then $\sigma^m(x)=c_{j-m}$
and 
$\hom(\mc{I}_n, k(c_{j-m})) 
\cong  \Gamma(\mc{I}_{c_{j-m}}/\mc{I}_{c_{j-m}}^2) 
\cong k^d$. For all other choices of $x$ one has 
$\hom(\mc{I}_n, k(x)) =k$. 
In other words, 
$$\Gamma_{AZ}(\widetilde{x})_m =\begin{cases}
k^d & \text{ if }
\ x=c_j\ \text{ with }\ j\geq m;\\
\noalign{\vskip 4pt}
k & \text{otherwise}\strut.
\end{cases}$$ 
This is equivalent to the assertion of the proposition.
\end{proof}

\begin{corollary}\label{point-class2}
Assume that  $R=R(X,c,\mc{L},\sigma)$ is generated in degree one.
Then the shift functor $s : P\mapsto P[1]_{\geq 0}$ does not 
induce an automorphism on 
the set of isomorphism classes of $R$-point modules.
\end{corollary}

\begin{proof} Set $M=\Gamma_{AZ}(\widetilde{c}_0)$. By
Proposition~\ref{point-class}, $\dim_kM_0=d\geq 2$ but $\dim M_j=1$ for $j>0$.
 Pick two linearly
independent elements $\alpha_i\in M_0$ and  write $M^i=\alpha_iR$. We will 
 show that the $M^i$ are non-isomorphic point modules for
$R$. 

By \cite[S2, p.252]{AZ1}, $M$ and
hence the $M^i$ are torsion-free, in the sense that they have no finite
dimensional submodules.   Since $R$ is generated in
degree one and $\dim_kM_i=1$, for $i>0$, this forces $(M^i)_{\geq 1} = M_{\geq
1}$ and so certainly the $M^i$ are point modules.  Suppose that there is an
isomorphism $\theta: M^1\to M^2$. Adjusting  $\theta$ by a scalar, we may
suppose that $\theta(\alpha_1)=\alpha_2$.  Now, $(M^1)_{\geq 1} =(M^2)_{\geq
1}$ is a shifted point module and so  its only automorphisms are given by
scalar multiplication. Hence the restriction of $\theta$ to $(M^1)_{\geq 1}$ is
given by multiplication by such  a scalar; say $\lambda$. Thus, for all $r\in R_1$,
$\lambda\alpha_1r=\theta(\alpha_1r)=\alpha_2r$. In other words,
$(\lambda\alpha_1-\alpha_2)R_1=0$,  contradicting the fact  that $M$ is
torsion-free. Thus the $M^i$ are nonisomorphic. Since 
$s(M^1) =(M^1)_{\geq 1}[1] =s(M^2)$ this implies that the shift functor $s$ 
is not an injection.
\end{proof}

The results of this section give several more 
ways in which the properties of $R$ differ from those of a 
 strongly noetherian $k$-algebra. Indeed, if $R$ were strongly
noetherian, then Theorem~\ref{non-represent} would contradict
\cite[Theorem~E4.3]{AZ2}, while, if $R$ were also generated in degree one,
then Corollary~\ref{non-represent3} would
 contradict Proposition~\ref{shift-aut}(3) and 
 Corollary~\ref{point-class2} would
 contradict Proposition~\ref{shift-aut}(1).

\section{Examples}
\label{section-critical}

In this last section we discuss the stringency of 
Assumptions~\ref{global-convention2} and give a number of examples where the
hypotheses are satisfied.  As is shown in \cite{Ke1},  projective schemes $X$
with  $\sigma$-ample  sheaves $\mc{L}$ exist in
abundance.  So the main issue is to determine 
the varieties $X$ and automorphisms  $\sigma$ for which 
there exists a critically dense orbit
  $\mc{C} = \{ \sigma^i(c) \}_{i \in \mb{Z}}$.

When  $X = \mb{P}^t$,   \cite[Section~14]{Ro} shows 
that $\mc{C}$ is critically dense  for generic choices
of $\sigma$ and $c$.  Below, we provide further examples of varieties
for which $\mc{C}$ is critically
dense for many choices of $\sigma, c$.  
The main technique is to find situations 
where one can reduce the problem of 
proving that $\mc{C}$ is critically dense to the problem 
of proving that $\mc{C}$ is dense.  
For this we use the following theorem of Cutkosky and Srinivas.

\begin{theorem} \cite[Theorem~7]{CS}
\label{alg group theorem} Let $G$ be a connected commutative algebraic group
 defined 
over an algebraically closed field $k$ of characteristic $0$. 
 Suppose that $g \in G$ 
is such that the cyclic subgroup $H = \langle g \rangle$ is dense 
in $G$.  Then any infinite subset of $H$ is dense in $G$.  \hfill $\Box$
\end{theorem}

Now let $X$ be any integral projective scheme with automorphism $\sigma$ and 
closed point $c$.  We think of $X$ as a variety in this section.

\begin{theorem}
\label{dense vs crit dense}
Assume that  $\cha k = 0$.   Suppose that there is an algebraic
group $G \subseteq \aut X$ with $\sigma \in G$
such that the map $\theta: G \times X \to X$ defined by $(g, x) \mapsto gx$ 
is regular.  
Then $\mc{C}=\{\sigma^i(c) : i\in {\mathbb Z}\}$ is dense in $X$ if and only if 
$\mc{C}$ is critically dense in $X$.
\end{theorem}

\begin{proof} 
Assume that $\mc{C}$ is dense in $X$.  Let $Z$ be the closure in $G$ of the subgroup
$H = \langle \sigma\rangle$ of $G$.  Then $Z$ is a subgroup of $G$ 
\cite[Proposition~7.4.A]{Hu}, which is abelian
since $H$ is. 

By definition, an algebraic group is of finite type, so 
$Z$ has finitely many connected components; let $Z_0$ be the connected
component containing $1$.   Since $\mc{C}$ is dense in $X$, the points $c_i =
\sigma^{-i}(c)$ are distinct (unless $X$ is a point, in which case the theorem
is trivial)  and so the automorphisms $\{\sigma^i\}_{i \in \mb{Z}}$ must be
distinct  points in $\aut X$.  Thus $H \cong \mb{Z}$ as groups and  it is easy
to see that $Z_0$ is the  closure of $H_0 = \langle\sigma^{e}\rangle$
for some $e \geq 1$.  Then the connected components of $Z$ are  the closures
$Z_j$ of the cosets $H_j = \sigma^jH_0$ for $0 \leq j \leq e-1$.  Now apply
Theorem~\ref{alg group  theorem} to $Z_0\supseteq H_0$ to show  that $H_0$
is critically dense in $Z_0$  and hence that each $H_j=\sigma^jH_0 $ is
critically dense in  $Z_j = \sigma^j Z_0$.

If $\mc{C}$  is not critically dense, pick $\{c_i\}_{i \in
I} \subseteq W \subsetneq X$,  where $W$ is closed and $I$ is infinite.
We can choose some $0 \leq j \leq e-1$ for which $I' = I \cap \{-j + e\mb{Z}\}$ 
is still infinite.  Now let $\phi: Z_j \times \{c\} \to X$ be the 
restriction of $\theta$, and let
$\rho_1: Z_j \times \{c\} \to Z_j$ be  the first projection.  Then
$\{ \sigma^{-i} \}_{i \in I'} \subseteq  \rho_1 (\phi^{-1}(W)) \subseteq
Z_j$.  Since $H_j$ is critically dense in $Z_j$, we must have 
$\rho_1 (\phi^{-1}(W)) = Z_j$.  But then
 $\{c_i \mid i \in -j + e\mb{Z} \} =\phi\rho_1^{-1}(H_j)\subseteq W$, 
and so $\mc{C} \subseteq \bigcup_{i = 0}^{e-1} \sigma^i(W) \subsetneq X$,
since $X$ is irreducible. This contradicts 
the density of $\mc{C}$ in $X$ and shows that $\mc{C}$ 
is critically dense. The other direction
is trivial.
\end{proof}

If  $\sigma, \tau$ are automorphisms of projective schemes $X$ and $ Y,$
respectively, and the hypotheses
of Theorem~\ref{dense vs crit dense} 
are satisfied for  $\sigma \in G \subseteq \aut X$ and 
$\tau \in G' \subseteq \aut Y$, 
then those hypotheses are also satisfied for 
$(\sigma, \tau) \in G \times G' \subseteq \aut (X \times Y)$.  
As an application of this   remark, we will
find critically dense sets in products of projective spaces.  Here we 
think of elements of $\pgl(t+1) = \aut \mb{P}^t$ as 
$(t+1) \times (t+1)$ matrices acting by left multiplication on
points of $\mb{P}^t$ written as 
column vectors of homogeneous coordinates.

\begin{example}\label{p-times-p}
Assume that $\cha k = 0$ and let $X = \mb{P}^{s_1} \times \mb{P}^{s_2}$.
Let $\tau_i\in \aut(\mb{P}^{s_i})$ be given by the diagonal matrix 
$\diag(1, p_{i1}, \dots, p_{is_i})$.  Assume that the 
$\{p_{ij}\}$ generate a multiplicative
 subgroup of $k$ isomorphic to $\mb{Z}^{s_1+s_2}$ and 
 write  $\sigma = (\tau_1, \tau_2)$.
 
If 
$c = \left((1: 1: \dots: 1), (1:1: \dots: 1)\right)\in X$, then  
$ \mc{C} = \{ \sigma^{i}(c) : i\in {\mathbb Z}\}$ is critically dense in $X$.
Moreover, if $\mc{L}$ is any very ample invertible sheaf on $X$ then $\mc{L}$ 
is $\sigma$-ample and so
\eqref{global-convention2} does hold for the data $(X, c, \mc{L}, \sigma)$.
\end{example}

\begin{proof} 
Since $\sigma \in \pgl(s_1+1) \times \pgl(s_2+1) \subseteq \aut X$,
by Theorem~\ref{dense vs crit dense} it is enough to  prove that $\mc{C}$ is
dense in $X$. If this is false,  there exists a proper closed set $Y$ with 
$\mc{C} \subseteq Y \subsetneq X$ and $\sigma(Y) = Y$. 
Let  $I$ be the defining ideal of $Y$ in the bigraded polynomial ring $U =
k[x_{10}, x_{11}, \dots, x_{1s_1}, x_{20}, x_{21}, \dots, x_{2s_2}]$. Thus $I$ 
 is bihomogeneous in the
$x$'s and satisfies $\phi(I) = I$, where $\phi$ is the automorphism
of $U$ corresponding to $\sigma$; explicitly, $\phi$ is defined
by $\phi(x_{ij}) = p_{ij} x_{ij}$.  Each
nonzero  bihomogeneous component $I_{uv}$ of 
$I$ is fixed by $\phi$ 
and so contains an  eigenvector, say $f=f_{uv}$,  for $\phi$.
  But the hypotheses on the $p_{ij}$ then
   force $f$ to be a single monomial  in the $x$'s, and
so $f(c) \neq 0$, contradicting $c \in Y$. Thus $\mc{C}$ is dense in $X$. 

The canonical sheaf on ${\mathbb P}^n$ is isomorphic to 
$\OO(-n-1)$ which is certainly minus ample. By 
\cite[Exercises~II.8.3 and II.5.11]{Ha}, the canonical sheaf on 
$ \mb{P}^{s_1} \times \mb{P}^{s_2}$ is therefore also minus ample.
Now apply \cite[Proposition~5.6]{Ke1} to see that $\mc{L}$ is $\sigma$-ample.
\end{proof}

Another large class of examples is provided by abelian varieties.
Here we recall some relevant definitions and refer the reader to 
 \cite{Ln} for the details.   An abelian variety $E$ is called \emph{simple}
if  the only (irreducible) abelian subvarieties of $E$ are itself and $0$.  
Two abelian varieties $E, E'$ are called \emph{isogenous}
 if there is a surjective
morphism of abelian  varieties  $E \to E'$ with
finite kernel. This is an equivalence relation on the set of abelian
varieties and every abelian variety $E$ is isogenous to a finite product of
simple abelian varieties \cite[Corollary, p.30]{Ln}. 

The following result is proved in \cite{RZ}.

\begin{proposition}
\label{prop-rz} 
Let $E  = E_1 \times E_2 \times \dots \times E_n$, where the $E_i$ 
are simple abelian
varieties.
  Let $a = (a_1, a_2, \dots, a_n) \in E$, and let $Z_a$ be the
Zariski closure in $E$ of $\{ia \}_{i \in \mb{Z}}$.   Then:
\begin{enumerate}
\item There is a countable set of closed subsets $Y_{\alpha} \subsetneq E$
 such that $Z_a = E$ for all $a \not \in \bigcup Y_{\alpha}$.
\item If the $E_i$ are pairwise
non-isogenous, then $Z_a = E$ if and only if each $a_i$ is a point of infinite
order in $E_i$. \qed
\end{enumerate}  
\end{proposition}

If $k$ is an uncountable field,  part~1 of the proposition shows that, for
a  sufficiently general point $a \in E$, we have $Z_a = E$. Part~2 shows that
in the special case where $E$ is a product of non-isogenous simples it is easy
to describe exactly for which $a \in E$ this happens. Thus the hypotheses of the
next result hold for generic $a\in E$.

\begin{example}\label{abelian-eg} Let $E$, $a$, and $Z_a$ be as in
Proposition~\ref{prop-rz}. Assume that $E$ is defined over a field of
characteristic zero and that $Z_a = E$.   Let $\sigma: E \to E$
be the translation automorphism defined by  $x \mapsto x + a$ and pick any
$c\in E$. Then $\mc{C} = \{ \sigma^i(c) \}_{i \in \mb{Z}}$ is critically dense
in $E$. If $\mc{L}$ is a very ample invertible sheaf on $E$ then $\mc{L}$  is
$\sigma$-ample and so \eqref{global-convention2} does hold for the data
$(E, c, \mc{L}, \sigma)$.
\end{example}

\begin{proof} 
The group of all translation automorphisms of $E$ is isomorphic
to $E$, so  it is an algebraic subgroup of $\aut E$ containing $\sigma$. 
By Theorem~\ref{dense vs crit dense} we just need to show that $\mc{C}$ is
dense.  This clearly does not depend on  the choice of $c$, so we may choose $c
= 0$.  By hypothesis, the closure of $\mc{C}$ is $Z_a = E$, as
required. The final assertion is proved in \cite{RZ}.
\end{proof}

When $k$ is a field of characteristic $p>0$,
there do exist examples of 
orbits $\mc{C}=\{\sigma^i(c_0)\}$ that are dense but not critically dense
\cite[Example~14.9]{Ro}. However, we cannot answer:

\begin{question}  If $\character k = 0$,   is every dense orbit 
$\{\sigma^i(c) : i\in \mathbb Z\}$ critically dense?
\end{question}

\bibliographystyle{amsplain}
\providecommand{\bysame}{\leavevmode\hbox to3em{\hrulefill}\thinspace}
\providecommand{\MR}{\relax\ifhmode\unskip\space\fi MR }
\providecommand{\MRhref}[2]{
  \href{http://www.ams.org/mathscinet-getitem?mr=#1}{#2}
}
\providecommand{\href}[2]{#2}

\end{document}